\documentclass[11pt]{article}
\usepackage{amsmath}
\usepackage{amsfonts}
\usepackage{color}

\def\U{{\cal U}}

\def\T{{\cal T}}

\def\U{{\cal U}}

\newfont{\Blackboard}{msbm10 scaled 1200}

\newfont{\roma}{cmr10 scaled 1200}

\providecommand{\bysame}{\leavevmode\hbox to3em{\hrulefill}\thinspace}
\providecommand{\MR}{\relax\ifhmode\unskip\space\fi MR }
\providecommand{\href}[2]{#2}

\newtheorem{theorem}{{}\hskip\parindent Theorem}[section]
\newtheorem{lem}{{}\hskip\parindent Lemma}[section]

\newtheorem{exl}{{}\hskip\parindent Example}[section]
\newtheorem{cor}{{}\hskip\parindent Corollary}[section]
\newtheorem{definition}{{}\hskip\parindent Definition}[section]
\newtheorem{rem}{{}\hskip\parindent Remark}[section]

\def\beq{\arraycolsep=1.5pt\begin{eqnarray}}
\def\eeq{\end{eqnarray}}

\large
\title{Second Order Necessary  Conditions for  Endpoints-Constrained  Optimal Control Problems on Riemannian manifolds \thanks{This work is  supported by the National Science
Foundation of China under grants 11401491 and 11931011, the Fundamental research funds for the Central Universities under grant 2682014CX052.
}
\date{}
\author{Li Deng\thanks{School of Mathematics,   Southwest Jiaotong University, Chengdu 611756, Sichuan Province, China. {\small\it E-mail:} {\small\tt
dengli@swjtu.edu.cn}.}, \ \ \
   and \ \ \ Xu Zhang\thanks{School of Mathematics, Sichuan University, Chengdu 610064, Sichuan Province, China. {\small\it E-mail:} {\small\tt
zhang$\_$xu@scu.edu.cn}.}}
}

\begin{document}
\maketitle

\begin{quote}
\begin{small}
{\bf Abstract} \,\,\,In this paper, we are concerned with optimal  control problems evolved on Riemannian manifolds, where the initial and final states satisfy some inequality and equality type constraints, and the control set is a separable metric space. We obtain the second order necessary conditions of integral and quasi-pointwise forms, both of which  work for Pontryagin type critical controls and involve  the curvature tensor. Also, we apply the  condition of integral form to the Bolza problem, where the initial and final states are subject to equality's type constraint.
\\[3mm]
\\[3mm]
{\bf MSC (2010) \,\,\,
49K15, 49K30, 93C15, 58E25, 70Q05} \\[3mm]
\end{small}
\end{quote}

\setcounter{equation}{0}

\section{Introduction }
\def\theequation{1.\arabic{equation}}
\hskip\parindent
 Let $n,j,k\in\mathbb N$ and $M$ be a complete simply connected, $n$-dimensional manifold with a Riemannian metric $g$. Let $\nabla$ be the Levi-Civita connection on $M$ related to $g$, $\rho(\cdot,\cdot)$ be the distance function on $M$, $T_xM$ be the tangent space of  $M$ at $x\in M$, and $T^*_x M$ be the  cotangent space. Denote by $\langle\cdot,\cdot\rangle$ and $|\cdot|$ the inner product and the norm over $T_xM$ related to $g$, respectively. Also, denote by  $ T M \equiv \bigcup\limits_{x\in M} T_xM$, $T^* M \equiv \bigcup\limits_{x\in M} T^*_x M$, $\mathcal X(M)$ and $C^\infty(M)$ the tangent bundle, the cotangent bundle, the set of smooth vector fields and the set of smooth functions on $M$, respectively. For $h\in C^\infty(M)$, we denote by $d h$  the differential of $h$.

Let $T>0$, $U$ be a metric  space, and the maps $f:[0,T]\times M\times U\to T M$,  $\phi_i: M\times M \to \mathbb R$ ($i=0,1,\cdots,j$) and $\psi: M\times M \to \mathbb R^k$  satisfy suitable assumptions to be given later.  We consider the following optimal control problem:
\begin{description}
\item[($P$)] Find a control $\bar{u}(\cdot)$ belonging to the set
\begin{equation}\label{436}
\U\equiv \big\{u(\cdot):[0,T]\to U;\ \ u(\cdot) \mbox{ is measurable}\big\},
\end{equation}
which minimizes the following cost functional:
$$
J(x(\cdot), u(\cdot)) \equiv \phi_{0}(x(0), x(T))
$$
subject to
\beq\label{25}
\dot{x}(t)=f(t,x(t),u(t)),\quad x(t)\in M\quad a.e. \,t\in[0,T],
\eeq
and
\begin{equation}\label{n22}
\left\{\begin{array}{l}{\phi_{i}(x(0), x(T)) \leq 0, i=1, \cdots, j,} \\ {\psi(x(0), x(T))=0,}\end{array}\right.
\end{equation}
where  $\dot{x}(t)=\frac{d }{dt}x(t)
$ for $t\in [0,T]$.
\end{description}

For the above problem, we call $\bar{u}(\cdot)$ an {\it optimal control}, the corresponding solution $\bar{x}(\cdot)$ to (\ref{25}) such that (\ref{n22}) holds an {\it optimal trajectory}, and $(\bar{x}(\cdot),\bar{u}(\cdot))$ an {\it optimal pair}.

Problem $(P)$ covers the optimal control problem that people usually consider:  the cost functional is of integral type (e.g., Problem $(P_1)$ in Section \ref{ec})  and the initial and final states satisfy some
inequality and equality type constraints. When $M$ is a Euclidean space, it is well-known that a necessary condition for optimal pairs is the classical Pontryagin type maximum principle. One way to derive this condition is to employ the first order needle variation of the control system, and then use the separation theorem for convex sets. As in calculus, Pontryagin type maximum principle is the first order necessary condition for optimal pairs in the sense of needle variation. A very natural question is, {\it what further necessary condition for optimal pairs can be obtained if the second order needle variation is introduced?} This sort of condition (if obtained) is called {\it second order necessary condition for optimal pairs in the sense of needle variation}.
The main purpose in the present work is to study such kind of second order necessary conditions for the above Problem $(P)$.

In the literatures, second order necessary conditions  in the sense of needle variation were studied for the case that  the state space is a Euclidean space or under some restrictive conditions. For example, Warga \cite{War78} considers the case that the initial and final states are subject to equality's type constraints. For Problem $(P)$ (with $M=\mathbb R^n$), Gilbert and Bernstein \cite{gb}  essentially require the control set to be a  subset of a compact set.  Lou \cite{l} and Cui, Deng and Zhang   \cite{cdz} consider the problems with free final states respectively on the Euclidean space  and Riemannian manifolds (while none of them needs the compact assumption on the control set). The paper \cite{cdz} also considered second order necessary conditions for optimal control problems on Riemannian
manifolds when the final state is fixed and the control set is an open set in a Euclidian space.

As mentioned above, one way to derive the Pontryagin's type maximum principle is  to use the separation theorem of convex sets, thanks to the very fact that the  set of all the first order variations of a control problem (in the sense of Ekeland's metric) is convex. However, no matter whether the control set $U$ is convex or not, the set of all the second order variations  may not  be convex, and consequently the same technique fails  in establishing the desired second order necessary condition.

To overcome the above difficulty, to the best of our knowledge,  most results in the literature are obtained by finding some convex set related to the high order variations of the  control problem.  For example, \cite[Theorem 3.1]{gb} and \cite[Theorem 6.1]{ls}  were focused on the case $M$ is a Euclidean space and $U$ is convex. Osmolovskii  \cite[Theorem 1]{O1} worked on the case that $M$ is a Euclidean space and the  pointwise control constraints are of inequality
type. P\'ales and  V. Zeidan \cite[Theorem 4.1]{pz}  concerned on the systems on a Euclidean space, with pure state and mixed control-state constraints.
It was observed in these results that, when fixing a critical variational direction, the set of the second order variations is convex, and the separation theorem of convex sets works.
 Sch\"attler and   Ledzewicz
\cite[Theorem 4.4.1]{sl} assumed $M$ is a differential manifold, and constructed an approximating cone related to the high order variations, and therefore for  a control-affine system on a Euclidean space with the control set $U$ being a closed ball, high order necessary conditions were obtained.
\cite[Theorem 20.6]{as} and \cite[Theorem 3.3]{cdz} were respectively concerned with the cases $M$ is a differential manifold and   $M$ is a Riemannian manifold, in both cases $U$ was required to open.
Warga
\cite[Theorem 2.2]{War78} considered the case that $M$ is a Euclidean space and $U$  is compact, and obtained the second order necessary condition by introducing relaxed controls.
%For more and early results, we refer the references where the mentioned papers cite.

Compared to the previous results, our results concern on a more general case: the state space is a manifold, the control set is neither compact nor convex, and the state is constrained both at the initial and final time.
We will encounter with two difficulties. The first one is, how to compute variations on manifolds. The second one is, with endpoint constraints, how to construct admissible trajectories around an optimal one. Actually, we use Riemannian geometric language to overcome the first difficulty. For the second one,
in order to obtain our main result (Theorem \ref{n295}), we first fix a critical variational direction, along which the Pontryagin's type maximum principle becomes trivial; then, we prove that along this direction, the set of all the second order variations of the control system is convex; finally, we apply the separation theorem of convex sets to this set. To obtain the quasi-pointwise second order necessary condition, we borrow the idea from \cite[Theorem 2.2(c)]{War78}.
It is notable that the curvature tensor enters explicitly our second order necessary condition.

This paper is organized  as follows: The main results are stated in Section 2, and we show their effectiveness by an example. In Section 3, we apply our main result to optimal control problems, where the initial and final states are subject to equality's type constraint, and the cost functional is of integral type. We also present a concrete example in this section. Sections 4 is devoted to the proof of our main results. In Section \ref{ap}, we list some notations, definitions and lemmas in Riemannian geometry, which are used in Sections 2-4.

\setcounter{equation}{0}
\hskip\parindent \section{Statement of the main results}
\def\theequation{2.\arabic{equation}}

\subsection{Notations and assumptions}
We begin with the following notions: Denote by $i(x)$, $|\T(x)|$, $\nabla \T$, $R$ , the injectivity radius (at the point  $x\in M$), the norm of the tensor field $\T$ (at the point $x\in M$), the covariant derivative of the tensor field $\T$ and the curvature tensor ( of $(M,g)$), respectively. For any $x,y\in M$ with $\rho(x,y)<\min\{i(x),i(y)\}$, there exists a unique shortest geodesic connecting $x$ and $y$. We denote   the parallel translation of a tensor from $x$ to $y$ along this geodesic  by   $L_{xy}$. For the definitions of the above notions, see Section \ref{ap}: Appendix.

Moreover, when a differentiable function $h: M\times M\to \mathbb R$ has two arguments, we denote by $\nabla_i h$ and $d_i h$ respectively the covariant derivative and differential of $h$ with respect to  the $i^{th}$ argument with $i=1,2$, i.e., for $X\in T M$ and $(x_1,x_2)\in M\times M$,
\begin{equation}\label{476}
 \nabla_i h(x_1,x_2)\left(X(x_i)\right)=d_ih(x_1,x_2)(X(x_i))= X(x_i)h(x_1,x_2),
\end{equation}
where we have used relation (\ref{ced}). When the vector-valued map $\eta=(\eta_1,\cdots,\eta_l)^\top: M\times M\to \mathbb R^l$ is differentiable, we denote by
\begin{align*}
\nabla_i\eta=(\nabla_i\eta_1,\cdots,\nabla_i\eta_l)^\top;\;d_i \eta=(d_i\eta_1,\cdots,d_i\eta_n)^\top,
\end{align*}
respectively the covariant derivative and differential of $\eta$ with respect to the $i^{th}$ argument with $i=1,2$.

The main assumptions are exhibited as follows:

\begin{description}
\item[$(C1)$] $(U,\tilde{d})$ is a separable metric space.
\item[$(C2)$] The map $f(=f(t,x,u)): [0, T]\times M\times U \to T M$ is measurable in $t$, continuous in $u$, and $C^1$ in $x$. The maps $\phi_i(=\phi_i(x_1,x_2)): M\times M \to \mathbb R$($i=0,\cdots,j$) and $\psi(=\psi(x_1,x_2))=(\psi_1,\cdots,\psi_k)^\top: M\times M
\to \mathbb R^k$ are $C^1$. Moreover,   there exist a constant $L>1$ and $x_0\in M$ such that,
\begin{equation}\begin{array}{l}\label{10}
|L_{x\hat x}f(s,x,u)-f(s,\hat x,u)|\leq L\rho(x,\hat x),
\\ |\phi_i(x_1,x_2)-\phi_i(\hat x_1,\hat x_2)|\leq L(\rho(x_1,\hat x_1)+\rho(x_2,\hat x_2)),\quad i=0,\cdots,j,
\\ |\psi(x_1,x_2)-\psi(\hat x_1,\hat x_2)|\leq L(\rho(x_1,\hat x_1)+\rho(x_2,\hat x_2)),
\\ |f(s,x_0,u)|\leq L,\end{array}
\end{equation} for all $s\in [0,T]$, $u\in U$, $x,\hat x\in M$ with $\rho(x,\hat x)\leq \min\{i(x),i(\hat x)\}$, and  $x_1, x_2,\hat x_1,\hat x_2\in M$.

\item[$(C3)$]The map $f(t,\cdot,u)$  are $C^2$  for all $(t,u)\in[0,T]\times U$.  The maps $\phi_i$ ($i=0,1,\cdots,j$) and $\psi$ are $C^2$. Moreover, for $f$ and $\varphi=\phi_0,\cdots,\phi_j,\psi$, there eixsts  a positive constant $L$ such that
\begin{align}\label{116}\begin{array}{l}
|\nabla_{x}f(t,x_1,u)-L_{\hat x_1 x_1}\nabla_x f(t,\hat x_1,u)|\leq L\rho(x_1,\hat x_1),
\\ |\nabla_1\varphi(x_1,x_2)-L_{\hat x_1 x_1}\nabla_1\varphi(\hat x_1, x_2)|\leq L\rho(x_1,\hat x_1),
\\ |\nabla_2\varphi(x_1,x_2)-L_{\hat x_2 x_2}\nabla_2\varphi( x_1, \hat x_2)|\leq L\rho(x_1,\hat x_1),
\end{array}\end{align}
for $t\in[0,T]$, $u\in U$, and  $x_1,\hat x_1,x_2,\hat x_2\in M$ with $\rho(x_i,\hat x_i)\leq \min\{i(x_i),i(\hat x_i)\}$ ($i=1,2$), where $\nabla_x f(t,x,u)$ is a tensor of type $(1,1)$ (see Section \ref{a1} for definition)  given by
\begin{align*}
\nabla_xf(t,x,u)(Y,X)=\nabla_Xf(t,\cdot,u)(Y),\quad\forall\, Y\in T_x^*M,\;X\in T_xM,
\end{align*}
and
\begin{align*}
\nabla_i\psi=(\nabla_i\psi_1,\cdots,\nabla_i\psi_k)^\top,\quad|\nabla_i\psi|=\sum_{\eta=1}^k|\nabla_i\psi_\eta|,\quad i=1, 2.
\end{align*}

\end{description}
It should be mentioned that, the first two lines of (\ref{10}) and (\ref{116}) are essentially Lipschitz  conditions, and they can be checked by \cite[Lemma 4.1]{cdz}.
In this paper, for $x\in M$, we denote by $\tilde X\in T_x^*M$ the dual covector of $X\in T_xM$, which is defined by
\begin{align*}
\tilde X(Y)=\langle X,Y\rangle,\quad\forall\,Y\in T_xM.
\end{align*}
Analogously, we denote by $\tilde \eta\in T_xM$ the dual vector of $\eta\in T_x^*M$, which is defined by $\langle\tilde\eta,Y\rangle=\eta(Y)$ for all $Y\in T_xM$.
Denote by
$H: [0,T]\times T^*M\times U\to \mathbb R$ the Hamiltonian function corresponding to Problem $(P)$, which is defined by
\begin{equation}\label{n30}H(t,x,p,\varphi,u)\equiv p(f(t,x,u)),
\end{equation}
for all $(t,x,p,u)\in[0,T]\times T^*M \times U$.

\subsection{Second order necessary condition of integral type}\label{m1}

In this section, we fix an optimal control $\bar{u}(\cdot)\in\mathcal U$.  Let $\bar x(\cdot)$ be a solution to (\ref{25}) associated to $\bar u(\cdot)$  such that (\ref{n22}) holds.
For abbreviation, we denote by
\begin{equation}\label{434}
[t]\equiv (t,\bar x(t),\bar u(t)),\quad\forall \,t\in[0,T].
\end{equation}
Set
\begin{equation}\label{n28}
\begin{array}{l}{I_{A O} \equiv\{0\}\cup \left\{i \in\{1, \cdots, j\} | \phi_{i}(\overline{x}(0), \overline{x}(T))=0\right\}} \\ {I_{N} \equiv\{0,1, \cdots, j\} \backslash I_{A O}.}\end{array}
\end{equation}
When $k>0,$ we introduce a Lagrange function $\mathcal{L} : M \times M \times \mathbb{R}^{1+j+k} \rightarrow \mathbb{R}$ defined by
\begin{equation}\label{n839}\mathcal{L}\left(y_{1}, y_{2}, \ell\right) \equiv \sum_{i=0}^{j} \ell_i \phi_{i}\left(y_{1}, y_{2}\right)+\ell_{\psi}^{\top} \psi\left(y_{1}, y_{2}\right),\end{equation}
 where $\ell=\left(\ell_0, \cdots,\ell_j, \ell_{\psi}^{\top}\right)^{\top} .$ When $k=0$ ,
we introduce a Lagrange function $\mathcal{L} : M \times M \times \mathbb{R}^{1+j} \rightarrow \mathbb{R}$ defined by $\mathcal{L}\left(y_{1}, y_{2}, \ell\right) \equiv$
$\sum_{i=0}^{j} \ell_i \phi_{i}\left(y_{1}, y_{2}\right),$ where $\ell=\left(\ell_{0}, \cdots, \ell_{j}\right)^{\top}.$  Either $k>0$ or $k=0$, we denote by $\nabla_i\mathcal L(y_1,y_2,\ell)$ and $d_i\mathcal L(y_1,y_2,\ell)$ respectively  the covariant derivative and the exterior derivative of $\mathcal L$ with respect to the variable $y_i$, where  $i=1,2$.

First, we shall introduce the Pontryagin's type maximum principle.

\begin{theorem}\label{f}
Assume conditions $(C1)-(C2)$  hold. If $(\bar x(\cdot), \bar u(\cdot))$ is optimal for Problem  (P), then there exists $\ell=\left(\ell_{\phi_{0}}, \ell_{\phi_{1}}, \cdots, \ell_{\phi_{j}}, \ell_{\psi}\right) \in \mathbb{R}^{1+j+k} \backslash\{0\}$ (if $k=0$, $\ell_\psi$ is omitted) satisfying
\begin{equation}\label{n33}\begin{array}{l}
\ell_{\phi_{i}} \in(-\infty, 0], \quad i=0, \cdots, j,
\\ \ell_{\phi_i}=0,\quad\textrm{if}\; i\in I_N,
\end{array}
\end{equation}
such that
\begin{equation}\label{n34}
H(t,\bar x(t),p^\ell(t),\bar u(t))=\max_{u\in U}H(t,\bar x(t),p^\ell(t),u),\quad a.e.\,t\in(0,T),
\end{equation}
where    $p^\ell(\cdot)$ is a covector field  along $\bar x(\cdot)$  verifying  the dual equation \begin{equation}\label{n35}\displaystyle\begin{cases}\nabla_{\dot{
 \bar{y}}(t)}p^\ell=-\nabla_xf[t]
 (p^\ell(t),\cdot),\quad \text{a.e.} \,t\in[0,T),\cr
p^\ell(T)=d_{2} \mathcal{L}(\bar{x}(0), \bar{x}(T), \ell),
 \end{cases}
 \end{equation}
 and the initial condition
 \begin{align}\label{n261}
p^\ell(0)=-d_{1} \mathcal{L}(\bar{x}(0), \bar x(T), \ell),
\end{align}
 and $\nabla_x f[t](p^\ell(t),\cdot)$ ($t\in[0,T]$) is a tensor given by
$$
\nabla_x f[t]\Big(p^\ell(t),X(\bar{x}(t))\Big)\equiv \nabla_{X(\bar{x}(t))}f(t,\cdot,\bar{u}(t))(p^\ell(t)),\qquad \forall\; X\in T M.
$$
\end{theorem}

Several remarks are in order.
\begin{rem}\label{n291}The initial and final conditions of the dual variable (see (\ref{n261}) and (\ref{n35})) is in fact the transversality condition. Actually, in \cite[Theorem 1.3, p. 132]{ly}, if the constraint set $S$ is $C^1$, it is just a special case of Theorem \ref{f}. The corresponding  transversatily condition  (see \cite[(1.9), p.131]{ly}) can be implied from (\ref{n261}) and (\ref{n35}).
\end{rem}
\begin{rem}\label{n290}It is obvious that the Pontryagin's type maximum principle (\ref{n34}) is equivalent to
\begin{align}\label{n260}
&\int_0^T\left(H(t,\bar x(t),p^\ell(t),\sigma(t))-H(t,\bar x(t),p^\ell(t),\bar u(t))\right)dt\nonumber
\\&+\left(\nabla_1\mathcal L(\bar x(0),\bar x(T), \ell)+p^\ell(0)\right)(W)\leq 0,\quad\forall \;(W, \sigma(\cdot))\in T_{\bar x(0)}M\times \mathcal U,
\end{align} where $p^\ell$  verifies (\ref{n35}). This condition is obtained  by computing the first order needle variation of the trajectory with respect to the initial state and the control, and by using the seperation theorem of convex sets. Therefore, Theorem \ref{f} can be viewed as the first order necessary condition of an optimal control.
\end{rem}

\begin{rem}\label{rem1} \cite[Theorem 12.15, p. 188]{as} shows the  Pontryagin's type maximum principle for a special case of Problem $(P)$: the state $x(\cdot)$ of system (\ref{25}) satisfies the endpoint constraint: $(x(0), x(T))$ belongs to a submanifold of $M\times M$, and the   cost functional is of integral type. Theorem \ref{f} is consistent with this result.

\end{rem}

As in calculus, when the first order necessary condition becomes trivial, we need to seek the second order necessary condition. Before doing this, we should  clarify what  ``the first order necessary condition is trivial'' means. To this end, we introduce the following definition.

\begin{definition}\label{n380}
A vector $\ell=\left(\ell_{\phi_{0}}, \ell_{\phi_{1}}, \cdots, \ell_{\phi_{j}}, \ell_{\psi}\right) \in \mathbb{R}^{1+j+k} \backslash\{0\}$ is called a Lagrange multiplier of an optimal pair $(\bar x(\cdot),\bar x(\cdot))$  for Problem $(P)$, if it satisfies (\ref{n33}), (\ref{n34}),  (\ref{n35}) and (\ref{n261}). A Lagrange multiplier $\ell$ is  normal, if $\ell_{\phi_0}<0$. Otherwise, it is called an abnormal Lagrange multiplier. For a Lagrange multiplier $\ell$, if there is a $u(\cdot)\in \mathcal U$ such that
\begin{align*}
H(t,\bar x(t),p^\ell(t),u(t))=H(t,\bar x(t),p^\ell(t),\bar u(t)),\quad a.e.\,t\in(0,T),
\end{align*}
we say that the Lagrange multiplier is trivial along the direction $u(\cdot)$.

\end{definition}

From the viewpoint of calculus, the first order necessary condition is trivial in direction $u(\cdot)\in\mathcal U$, if all the Lagrange multipliers are  trivial along $u(\cdot)$. In what follows, we introduce ``critical direction'', along which all the Lagrange multipliers are trivial.

\begin{definition}\label{n40}  A control $u(\cdot)\in  \mathcal U$ is called a Pontryagin's type critical direction, if there exists a $V\in T_{\bar x(0)}M$ such that
\begin{equation}\label{n53}
\begin{array}{l}{\nabla_{1} \phi_{i}(\bar{x}(0), \bar{x}(T))(V)+\nabla_{2} \phi_{i}(\bar{x}(0), \bar{x}(T))\left(X_{u,V}(T)\right) \leq 0, \forall \,i \in I_{A O}}, \\ {\nabla_{1} \psi(\bar{x}(0), \bar{x}(T))(V)+\nabla_{2} \psi(\bar{x}(0), \bar{x}(T))\left(X_{u,V}(T)\right)=0 \quad(\textrm {omit if } k=0)},\end{array}
\end{equation}
where $X_{u,V}(\cdot)$ is a vector field along $\bar x(\cdot)$ and satisfies
\begin{equation}\label{n14}\left\{\begin{array}{l}
\nabla_{\dot{\bar x}(t)} X_{u,V}=\nabla_xf[t](\cdot,X_{u,V}(t))+f(t,\bar x(t),u(t))-f[t],\;a.e.\, t\in(0,T),
\\ X_{u,V}(0)=V,
\end{array}\right.
\end{equation}
and
\begin{align}\label{pdv}
\nabla_i\psi (X)=\left(\nabla_i\psi_1(X),\cdots,\nabla_i\psi_k(X)\right)^\top,\quad\forall\,X\in T M,\,i=1,2.
\end{align}
\end{definition}

 Actually, if $u(\cdot)$ is a Pontryagin's type critical direction, then for any Lagrange multiplier $\ell$,  by using (\ref{n33}), (\ref{n35}), (\ref{n261}), (\ref{n260})  and integration by parts over $[0,T]$, we can obtain
\begin{align}\label{n850}
0\leq& \nabla_1\mathcal L(\bar x(0),\bar x(T),\ell)(V)+\nabla_2\mathcal L(\bar x(0),\bar x(T),\ell)(X_{u,V}(T))\nonumber
\\=&\int_0^T(H(t,\bar x(t),p^\ell(t),u(t))-H(t,\bar x(t),p^\ell(t),\bar u(t)))dt\nonumber
\\\leq &0,
\end{align}
which implies $\ell$ is trivial along $u(\cdot)$.

To introduce the second order necessary condition, we assume conditions $(C1)-(C3)$ hold, and adopt the following notaions.
\begin{equation}\label{524}\begin{array}{l}
\\ \nabla_x H(t,x,p,u)(X):=\nabla_x f(t,x,u)(p,X),
\\ \nabla_x^2 H(t,x,p,u)(X,Y):=
 \nabla_x^2f(t,x,u)(p,X,Y),
\end{array}\end{equation}
for all $(t,x,p,u)\in[0,T]\times T^*M \times U$ and  $X,Y\in T M$. For the definition of covariant derivative of tensors, we refer to Section \ref{a1}.

If $u(\cdot)$ is a  critical direction, and $V\in T_{\bar x(0)}M$ and $u(\cdot)$ satisfies (\ref{n53}), we set
\begin{equation}\label{n43}
\begin{array}{l}{I_{0}^{\prime} \equiv I_{N} \cup\left\{i \in I_{A O} | \nabla_{1} \phi_{i}(\bar{x}(0), \bar{x}(T))(V)+\nabla_{2} \phi_{i}(\bar{x}(0), \bar{x}(T))\left(X_{u, V}(T)\right)<0\right\}}, \\
 {I_{0}^{\prime \prime} \equiv\{0,1, \cdots, j\} \backslash I_{0}^{\prime}}.\end{array}
\end{equation}

Our result on the second order necessary condition of integral type can be stated as follows.

\begin{theorem}\label{n295}
Assume conditions $(C1)-(C3)$ hold. Let $(\bar x(\cdot),\bar u(\cdot))$ be an optimal pair for Problem $(P)$.  For any critical direction $u(\cdot)\in\mathcal U$ with $V\in T_{\bar x(0)}M$ such that (\ref{n53}) holds,
there exists another Lagrange multiplier $\hat\ell=(\hat\ell_0,\hat\ell_1,\cdots,\hat\ell_j,\hat\ell_\psi)\in \mathbb R^{1+j+k}\setminus\{0\}$ satisfying
\begin{align}\label{n230}
\hat\ell_i\leq 0,\quad i=0,1,\cdots,j,
\\\label{n231} \hat\ell_i=0,\quad\textrm{if}\;i\notin I_0^{\prime\prime},
\end{align}
such that
\begin{align}\label{n236}
&\int_0^T\Big(\nabla_x^2H\{ t\}^{\hat\ell}(X_{u,V}(t),X_{u,V}(t))+2(\nabla_xH(t,\bar x(t),p^{\hat\ell}(t),u(t))\nonumber
\\&-\nabla_xH\{ t\}^{\hat\ell})(X_{u,V}(t))
-R(\tilde p^{\hat\ell}(t),X_{u,V}(t),f[t],X_{u,V}(t))\Big)dt\nonumber
\\&+\nabla_1^2\mathcal L(\bar x(0),\bar x(T),\hat\ell)(V,V)
+2\nabla_2\nabla_1\mathcal L(\bar x(0),\bar x(T),\hat\ell)(V,X_{u,V}(T))\nonumber
\\&+\nabla_2^2\mathcal L(\bar x(0),\bar x(T),\hat\ell)(X_{u,V}(T),X_{u,V}(T))\leq 0,
\end{align}
where we adopt the notation
\begin{align}\label{n241}
\{ t\}^{\hat\ell}\stackrel{\triangle}{=}(t,\bar x(t),p^{\hat\ell}(t),\bar u(t)),\quad \forall\;t\in[0,T],
\end{align}
$p^{\hat\ell}(\cdot)$ is the solution to (\ref{n35}) with $\ell$ replaced by $\hat\ell$, and $\tilde p^{\hat\ell}(t)$ is the dual vector of $p^{\hat\ell}(t)$ ($t\in[0,T]$).
\end{theorem}
\begin{rem}\label{n360}
In \cite[Theorem 6.4]{gb} the second odder necessary conditions for Problem $(P)$ was considered, where the state space  is  a Euclidean space.  It needs the following assumptions: 1)  the control set is a subset of a compact set; 2) the control system is convex, or all the Lagrange multipliers are normal (see Definition \ref{n380}). While Theorem \ref{n295} does not need these conditions, and it considers a more general case: the state space is a  Riemannian manifold. What is new is that, the curvature tensor ``$R$'' appears in (\ref{n236}). It is necessary to mention that, when the state space is a Euclidean space,  the curvature tensor  is zero in (\ref{n236}), and the corresponding result is consistent with  \cite[Theorem 6.4]{gb}.
\end{rem}

We may apply Theorem \ref{n295} to the following example to check whether a control is optimal, while the same example is solved by the second order necessary condition of quasi-pointwise form (see \cite[Example II.]{War78}).
\begin{exl}\label{ei}
Minimize
\begin{align*}
\phi_0(x_1(0),x_2(0),
x_1(T),x_2(T))\stackrel{\triangle}{=}x_1(T),
\end{align*}
subject to
\begin{align*}
\begin{array}{l}
\left(\begin{array}{c}\dot x_1(t)
\\ \dot x_2(t)\end{array}\right)=\left(\begin{array}{c}x_2(t)(u_1(t)+u_2(t))\\ u_2(t)-x_1(t)
\end{array}\right),\; a.e. \,t\in(0,T),
\end{array}
\end{align*}
and
\begin{align*}
(u_1(t),u_2(t))^\top\in [0,1]\times[-1,1],\;a.e. \,t\in[0,T],
\\
\psi(x_1(0),x_2(0),
x_1(T),x_2(T))\stackrel{\triangle}{=}\left(x_1(0),x_2(0),x_2(T)\right)^\top=(0,0,0)^\top.
\end{align*}
If  control $\bar u(\cdot)=(0,0)^\top$ is an optimal control of the above problem, the corresponding trajectory is $\bar x(\cdot)=(0,0)^\top$. Then, there exists $\ell_0\leq 0$ and $\ell_{\psi}^1, \ell_\psi^2, \ell_\psi^2\in \mathbb R$ with $(\ell_0,\ell_\psi^1,\ell_\psi^2,\ell_\psi^3)^\top\neq 0$ such
that
\begin{align}\label{epmp}
\max_{-1\leq u_2\leq 1}(-\ell_\psi^2 u_2)=-\ell_\psi^2 0,
\end{align}
and
\begin{align*}
\left(\dot p_1,\dot p_2\right)(t)=(p_2(t), 0), \;a. e. \,t\in(0,T),
\\ p_1(0)=-\ell_\psi^1,\,p_2(0)=-\ell_\psi^2,\,p_1(T)=\ell_0,\,p_2(T)=\ell_\psi^3.
\end{align*}
 We obtain that $\ell_\psi^3=\ell_\psi^2=0$ and $-\ell_\psi^1=\ell_0$.
Then, the Lagrange multiplier $(\ell_0,\ell_\psi^1,\ell_\psi^2,\ell_\psi^3)^\top$ is unique up to a positive factor. We take $\ell_0=-1$ and consequently $\ell_\psi^1=1$. Then $p_1(t)\equiv -1$ and $p_2(t)\equiv 0$ for all $t\in[0,T]$.

Take $u_1(t)\equiv 1$ and $u_2(t)=-I_{[0,\frac{T}{2}]}(t)+I_{(\frac{T}{2},T]}(t)$ for $t\in[0,T]$, where $I_A(\cdot)$ is the indicator function of set $A$. The variational equatin along direction $(u_1(\cdot),u_2(\cdot)) $ is as follows:
\begin{align*}
\left\{\begin{array}{l}
\left(\begin{array}{c}
\dot X_1(t)\\ \dot X_2(t)
\end{array}\right)=\left(\begin{array}{c}
0\\ -X_1(t)
\end{array}\right)+\left(\begin{array}{c}
0\\ u_2(t)
\end{array}\right),\;a.e. \,t\in(0,T),
\\ X_1(0)=X_2(0)=0.
\end{array}\right.
\end{align*}
We can check that
$(u_1(\cdot),u_2(\cdot))$
is a critical direction.
In this case, the left hand side of (\ref{n236}) is reduced to
\begin{align*}
-\int_0^Tu_2(t)\int_0^tu_2(s)dsdt-\int_0^T\int_0^tu_2(s)dsdt=\frac{T^2}{4},
\end{align*}
which means that the second order necessary condition does not hold, and consequently control $(0,0)$ is not optimal.

\end{exl}

\subsection{Second order necessary condition of  quasi-pointwise form}\label{psnc}

In this subsection, we seek the second order necessary condition of  quasi-pointwise form, by borrowing some idea from
 \cite[Theorem 2.2 (c)]{War78}.

Without loss of generality, we shall consider a simpler case of Problem $(P)$:
\begin{description}
\item[($P_1$)]  Find $\bar u(\cdot)\in\mathcal U$, which minimizes $\phi_0(x(T))$
 subject to (\ref{25}), $x(0)=x_0$, $\phi_\eta(x(T))\leq 0$ $(\eta=1,\cdots,j)$ and $\psi(x(T))\stackrel{\triangle}{=}(\psi_1(x(T)),\cdots,\psi_k(x(T)))^\top
 =0 (\in\mathbb R^k)$.
\end{description}

Assume $(\bar x(\cdot), \bar u(\cdot))$ is an optimal pair of problem $(P_1)$,  and there is a unique Lagrange multiplier (up to a positive factor)
$(\ell_0,\ell_1,\cdots,\ell_j,\ell_\psi^\top)^\top\in((-\infty,0]^{j+1}\times\mathbb R^k )\setminus\{0\}$. We adopt notation $[t]$ in (\ref{434}), and
\begin{align}\label{pdv1}
\eta^\top d\psi(x)=\sum_{i=1}^k\eta_i d\psi_i(x),\;\forall\eta=(\eta_1,\cdots,\eta_k)^\top\in\mathbb R^k,\;\forall x\in M,
\end{align}
where $d$ is the exterior derivative.

 By Theorem \ref{f}, we have
\begin{align*}
H(t,\bar x(t), p(t),\bar u(t))=\max_{u\in U}H(t,\bar x(t), p(t),u),\quad a.e. \,t\in[0,T],
\end{align*}
where $H$ is defined in (\ref{n30}) and  $p(\cdot)$ is the covector along $\bar x(\cdot)$ satisfying
\begin{align}\label{n750}
\left\{\begin{array}{ll}
\nabla_{\dot{\bar x}(t)}p=-\nabla_xf[t](p(t),\cdot),&a.e. t\in(0,T),
\\ p(T)=\sum_{\eta
=0}^j\ell_\eta d \phi_\eta(\bar x(T))+\ell_\psi^\top d\psi(\bar x(T)).&
\end{array}\right.
\end{align}
For $t\in[0,T]$, set
\begin{align*}
U(t)=\{u\in U;\, H(t,\bar x(t),p(t),u)=H(t,\bar x(t),p(t),\bar u(t))\}.
\end{align*}

Let $\left\{e_{1}, \cdots, e_{n}\right\} \subset T_{\bar {x}(0)} M$ be an orthonormal basis.  Denote by $\{d_1,\cdots,d_n\}\subset T_{\bar x(0)}^*M$ the dual basis to $\{e_1,\cdots,e_n\} $, i.e. $d_i(e_j)=\delta_i^j$ for $i,j=1,\cdots,n$, where $\delta_i^j$ are the Kronecker delta symbols.  For $t \in(0, T],$ denote respectively by $e_i(t) \equiv L_{\bar{x}(0) \bar{x}(t)}^{\bar x(\cdot)}e_i$ and $d_i(t)\equiv L_{\bar{x}(0) \bar{x}(t)}^{\bar x(\cdot)}d_i (i=1, \cdots, n)$ the parallel translations of $e_{i}$ and $d_i$
from $\bar{x}(0)$ to $\bar{x}(t)$ along the curve $\bar{x}(\cdot)$.  Then, it follows from (\ref{r1}) and (\ref{r50}) that $\{e_1(t),\cdots,e_n(t)\}$ is an orthonormal basis at $T_{\bar x(t)}M$, and $\{d_1(t),\cdots,d_n(t)\}$ is the dual basis to it.
Consequently, for $(t,u) \in[0, T]\times U$, we can express tensors $\nabla_x^2f[t]$, $\nabla_{x} f(t,\bar x(t),u)$, $f(t,\bar x(t),u)$ and $p(t)$ respectively  by  $\nabla_x^2f[t]=\sum_{i,\xi,\zeta=1}^n B_{i\xi\zeta}(t)e_i(t)\times d_\xi(t)\times d_\zeta(t)$, $\nabla_{x} f(t,\bar x(t),u)=\sum_{i, j=1}^{n} A_{i j}(t,u) e_{i}(t) \otimes d_{j}(t)$, $f(t,\bar x(t),u)=\sum_{i=1}^n f^i(t,u)e_i(t)$ and $p(t)=\sum_{i=1}^n p_i(t)d_i(t)$, where
for $i,j, \xi,\zeta=1,\cdots, n$,
\begin{align}\label{n738}\begin{array}{lll}
&B_{i\xi\zeta}(t)=\nabla_x^2f[t](d_i(t),e_\xi(t),e_\zeta(t)), &p_i(t)=p(t)(e_i(t)),\quad\quad\quad\quad
\\ &A_{i j}(t,u)=\nabla_{x} f(t,\bar x(t),u)\left(d_{i}(t), e_{j}(t)\right),&  f^i(t,u)= f(t,\bar x(t),u)(d_i(t)).
\end{array}\end{align}
Denote by
\begin{align}\label{gcs20}
\begin{array}{ll}
\vec f(t,u)=(f^1(t,u),\cdots,f^n(t,u))^\top,\;\vec p(t)=(p_1(t),\cdots, p_n(t)),
\\ A(t,u)=(A_{ij}(t,u))_{i,j=1}^n.
\end{array}
\end{align}
 Denote by $Z:[0,T]\to \mathbb R^{n\times n}$ the solution to
\begin{align}\label{n830}
\left\{\begin{array}{l}\dot Z(t)=-Z(t)A(t,\bar u(t)),\;t\in[0,T),
\\ Z(0)=I_{n},
\end{array}\right.
\end{align}
where $I_n$ is the identity matrix in $\mathbb R^{n\times n}$. Set
\begin{align}\label{n745}
\mathcal A(t,u)=& Z(t)[\vec f(t,u)-\vec f(t,\bar u(t))];
\\
\Delta  H(t,u)=&\Big([\nabla_x H(t,\bar x(t),p(t),u)-\nabla_xH(t,\bar x(t),p(t),\bar u(t))](e_1(t)),\cdots,\nonumber
\\ &[\nabla_x H(t,\bar x(t),p(t),u)-\nabla_xH(t,\bar x(t),p(t),\bar u(t))](e_n(t))\Big)^\top,\nonumber
\end{align}
for $(t,u)\in[0,T]\times U$.

\begin{definition}\label{n741} Given a  map $\mathcal L:\mathbb R\to\mathbb R^N$ ($N\in\mathbb N$), we say $\mathcal L$ is approximately continuous at $t_0\in\mathbb R$, if for any $\epsilon>0$,  the relation
\begin{equation}\label{n718}\lim_{r\to 0^+}r^{-1}\left|\{t\in\mathbb R; |t-t_0|\leq r, |\mathcal L(t)-\mathcal L(t_0)|>\epsilon\}\right|=0
\end{equation}
holds.

\end{definition}

It follows  from \cite[Theorem 3, p. 47]{eg92} that, if $\mathcal L$ is measurable, then it is approximately continuous almost every.

\begin{theorem}\label{n770}
  Assume all the assumptions in Theorem \ref{n295} hold. Let $(\bar x(\cdot), \bar u(\cdot))$ be an optimal pair with a unique Lagrange multiplier (up to a positive factor) $(\ell_0,\cdots,\ell_j, \ell_\psi^1,\cdots,\\ \ell_\psi^k)^\top\in (-\infty,0)\times (-\infty,0]^{j}\times\mathbb R^k$.  Let $u(\cdot)\in\mathcal U$ be such that $u(t)\in U(t)$ a.e.  $t\in[0,T]$. Assume $\tau_0, \tau_1, \cdots, \tau_\ell\in (0,T)$  with $\ell\geq k+j$   satisfy the following properties: i) $0<\tau_0<\cdots<\tau_\ell<T$; ii)  $\mathcal A(\cdot,u(\cdot))$ and $Z(\cdot)\left(A(\cdot,u(\cdot))-A(\cdot,\bar u(\cdot))\right)Z^{-1}(\cdot)$ are both approximately continuous at $\tau_0, \cdots, \tau_\ell$ and
 $0^{j+k}\in Int \,co \{(\nabla\Phi_1^\top,\cdots,\nabla\Phi_j^\top, \nabla\Psi^\top)^\top\mathcal A(\tau_i,u(\tau_i))\}_{i=0}^\ell,$
where $0^{j+k}$ is the zero in space $\mathbb R^{j+k}$,  and  ``$Int A$'' and "$co\,A$" respectively denote the interior and the convex hull of set $A$;
 iii) There exist  $\beta_0, \beta_1, \cdots,\beta_\ell\in(0,+\infty)$ such that
\begin{align}\label{n800}
\nabla\Phi_i\sum_{\eta=0}^\ell\beta_\eta\mathcal A(\tau_\eta,u(\tau_\eta))=0,
\quad \nabla\Psi\sum_{\eta=0}^\ell\beta_\eta\mathcal A(\tau_\eta,u(\tau_\eta))=0,\;i=1,\cdots,j,
\end{align}
  where
\begin{align}\label{n780}\begin{array}{l}\nabla\Phi_i=[\nabla \phi_i(\bar x(T))(e_1(T)),\cdots,\nabla\phi_i(\bar x(T))(e_n(T))]Z^{-1}(T),\;i=1,\cdots,j,
\\[2mm]
\nabla\Psi=\left(\begin{array}{ccc}
\nabla\psi_1(\bar x(T))(e_1(T))&\cdots&\nabla\psi_1(\bar x(T))(e_n(T))
\\ \vdots&\vdots&\vdots
\\\nabla\psi_k(\bar x(T))(e_1(T))&\cdots&\nabla\psi_k(\bar x(T))(e_n(T))
\end{array}
\right)Z^{-1}(T).
\end{array}
\end{align}
Then, it holds that
 \begin{align}\label{n740}
\begin{array}{ll}&\sum_{i=0}^{\ell}\sum_{ \eta,\hat\eta=0}^i\beta_\eta \beta_{\hat\eta} \mathcal A(\tau_\eta, u(\tau_\eta))^\top \int_{\tau_i}^{\tau_{i+1}} \left(Z^{-1}(t)\right)^\top \Big(\nabla_x^2H\{t\}(e_\xi(t),e_\zeta(t))
\\&-R(\tilde p(t),e_\xi(t),f[t],e_\zeta(t))\Big)_{\xi,\zeta=1}^n
Z(t)^{-1}dt
  \mathcal A(\tau_{\hat\eta},u(\tau_{\hat\eta}))
\\&+\sum_{\eta=0}^\ell\Big(2\beta_\eta \Delta H(\tau_\eta,u(\tau_\eta)) Z^{-1}(\tau_\eta)\sum_{0\leq i<\eta}\beta_i\mathcal A(\tau_i,u(\tau_i))
\\&+(\beta_\eta)^2 \Delta H(\tau_\eta,u(\tau_\eta))Z^{-1}(\tau_\eta)\mathcal A(\tau_\eta,u(\tau_\eta))\Big)
+\sum_{\eta,\hat\eta=0}^\ell \beta_\eta\beta_{\hat\eta}\mathcal A(\tau_\eta,u(\tau_\eta))^\top
\\&Z^{-1}(T)^\top\left(\sum_{i=0}^j\ell_{\phi_i}\nabla^2\Phi_i+\sum_{\eta=1}^k\ell_\psi^\eta\nabla^2\Psi_\eta\right) Z^{-1}(T) \mathcal A(\tau_{\hat\eta},u(\tau_{\hat\eta}))\leq 0,
\end{array}\end{align}
where $\tau_{\ell+1}=T$, $\tilde p(t)$ is the dual vector of $p(t)$ for $t\in[0,T]$, and
\begin{align}\label{n759}\begin{array}{ll}
&\{t\}=(t,\bar x(t),p(t),\bar u(t)),\;\forall\, t\in[0, T];\\
[2mm]&\nabla^2\Phi_i=\Big(\nabla^2\phi_i(\bar x(T))(e_\xi(T),e_\zeta(T))\Big)_{\xi,\zeta=1}^n,\quad i=0,1,\cdots,j ;
\\[2mm]&\nabla^2\Psi_\eta=\Big(\nabla^2\psi_\eta(\bar x(T))(e_i(T),e_\xi(T))\Big)_{i,\xi=1}^n,\;\eta=1,\cdots,k.
\end{array}\end{align}

\end{theorem}

\begin{theorem}\label{n790} Assume all the assumptions in Theorem \ref{n770} hold and $U$ is compact.  Then, there exists a subset $\mathcal T\subset[0,T]$ with measure $T$ such that,
for any $\tau_0, \tau_1,\cdots, \tau_\ell\subset \mathcal  T$ with $0
<\tau_0<\cdots<\tau_\ell<T$ and $\ell\geq k+j$,  any $r_i\in U(\tau_i)$ $(i=0,\cdots,\ell)$ and $\beta_0,\cdots,\beta_\ell\in(0,+\infty)$ satisfying
\begin{align}\label{n791}
\nabla\Phi_i\sum_{\eta=0}^\ell\beta_\eta\mathcal A(\tau_\eta,r_\eta)=0,\;\nabla\Psi\sum_{\eta=0}^\ell\beta_\eta\mathcal A(\tau_\eta,r_\eta)=0,\; i=1,\cdots,\ell;
\\ \label{n792}
0^{j+k}\in Int\, co\{(\nabla\Phi_1^\top,\cdots,\nabla\Phi_j^\top,\nabla\Psi^\top)^\top\mathcal A(\tau_\eta,r_\eta)\}_{\eta=0}^\ell,
\end{align}
 it holds that
\begin{align}\label{n810}
\begin{array}{ll}
&\sum_{i=0}^{\ell}\sum_{ \eta,\hat\eta=0}^i\beta_\eta \beta_{\hat\eta} \mathcal A(\tau_\eta,r_\eta)^\top \int_{\tau_i}^{\tau_{i+1}} \Big(Z^{-1}(t)\Big)^\top \Big(\nabla_x^2H\{t\}(e_\xi(t),e_\zeta(t))
\\&-R(\tilde p(t),e_\xi(t),f[t],e_\zeta(t))\Big)_{\xi,\zeta=1}^n
Z(t)^{-1}dt
  \mathcal A(\tau_{\hat\eta},r_{\hat\eta})
\\&+\sum_{\eta=0}^\ell\Big(2\beta_\eta \Delta H(\tau_\eta,r_\eta) Z^{-1}(\tau_\eta)\sum_{0\leq i<\eta}\beta_i\mathcal A(\tau_i,r_i)
\\&+(\beta_\eta)^2 \Delta H(\tau_\eta,r_\eta)Z^{-1}(\tau_\eta)\mathcal A(\tau_\eta,r_\eta)\Big)
+\sum_{\eta,\hat\eta=0}^\ell \beta_\eta\beta_{\hat\eta}\mathcal A(\tau_\eta,r_\eta)^\top Z^{-1}(T)^\top
\\& \left(\sum_{i=0}^j\ell_{\phi_i}\nabla^2\Phi_i
+\sum_{\eta=1}^k\ell_\psi^\eta\nabla^2\Psi_\eta\right) Z^{-1}(T) \mathcal A(\tau_{\hat\eta},r_{\hat\eta})\leq 0,
\end{array}
\end{align}
where $\tau_{\ell+1}=T$.
\end{theorem}

\begin{rem}
Theorem \ref{n790} can be used to check Example \ref{ei}, see \cite[Example II]{War78} for details. \cite[Theorem 2.2(c)]{War78} considers problem $(P_1)$ when $M$ is a Euclidean space, and gives the quasi-pointwise second order necessary condition. Theorem  \ref{n790} extends this result to the case on manifolds.
\end{rem}
%
% Using the Pontryagin's type maximum principle, Sussmann (\cite[Example 5.10.2]{s2}) showed that the locally shortest curve connecting two fixed points on a Riemannian manifold must be a geodesic.   \cite[Example 4.2]{cdz} applies the normal case of  the second order necessary conditions of optimal controls  (  \cite[Theorem 3.3]{cdz}), and obtain the well known second variation of energy ( see Example \ref{510} below ). As mentioned in Remark \ref{rem.1},  Theorem \ref{430} is still applicable for the following  example, and we omit its proof.
%
%
%

%

%
%
%

 \setcounter{equation}{0}
\hskip\parindent \section{Applications}
\def\theequation{3.\arabic{equation}}

We shall apply  Theorem \ref{n295} to a special case of Problem $(P)$, and give an example as an application.

\subsection{Optimal control problems with endpoints constraints}\label{ec}

Given maps $f^0:[0,T]\times M\times U\to \mathbb R$, $\psi_1:[0,T]\to \mathbb R^{k_1}$,$ \psi_2:[0,T]\to\mathbb R^{k_2}$ and $ h:M\to \mathbb R$($k_1,k_2\in\mathbb N$), we  consider the following problem:
\begin{description}
\item[$(P_2)$] Minimize
\begin{align*}
J(x(\cdot), u(\cdot))\stackrel{\triangle}{=}\int_0^T f^0(t,x(t),u(t))dt+h(x(T)),
\end{align*}

which is subject to (\ref{25}), $u(\cdot)\in\mathcal U$, $\psi_1(x(0))=0$ and $\psi_2(x(T))=0.$
\end{description}

\begin{theorem}\label{n300}
Assume $(C1)$ holds, the maps $f:[0, T]\times M\times U\to  T M$ and $f^0:[0, T]\times M\times U\to \mathbb R$ are measurable in $t$, continuous in $u$, and $C^1$ in $x$. The maps $h, \psi_1,\psi_2$ are $C^1$. Moreover, there exists a constant $L>1$ such that the fisrt and last lines of (\ref{10}) hold both for $f$ and $f^0$, and for $\varphi=h,\psi_1,\psi_2$, it holds that
\begin{align*}
|\varphi(x)-\varphi(\hat x)|\leq L\rho(x,\hat x),\; i=1,2,
\end{align*}
where $x,\hat x\in M$ satisfy $\rho(x,\hat x)\leq\min\{i(x),i(\hat x)\}$. Then, if $(\bar x(\cdot),\bar u(\cdot))$ is an optimal pair for problem $(P_2)$,  there exists $\ell=(\ell_0,\ell_{\psi_1},\ell_{\psi_2})\in\mathbb R^{1+k_1+k_2}\setminus\{0\}$ such that
\begin{align}\label{n335}
&\ell_0\leq 0,
\\
&\label{n304}\max_{u\in U} H_2(t,\bar x(t),p_1^\ell(t),u,\ell_0)=H_2(t,\bar x(t),p_1^\ell(t),\bar u(t),\ell_0), \;a.e.\,t\in(0,T),
\end{align}
where $p_1^\ell$ is a covector field along $\bar x(\cdot)$ satisfying
\begin{align}\label{n305}
\left\{\begin{array}{ll}
\nabla_{\dot{\bar x}(t)}p_1^\ell=-\nabla_xf[t](p_1^\ell(t),\cdot)-\ell_0 d_xf^0[t],&a. e. \,t\in(0,T),
\\ p_1^\ell(0)=-\ell_{\psi_1}^\top d\psi_1(\bar x(0)),
\\ p_1^\ell(T)=\ell_0 d h(\bar x(T))+\ell_{\psi_2}^\top d\psi_2(\bar x(T)),
\end{array}
\right.
\end{align}
with $d_x f^0(t,x,u)$ being the exterior derivative of $f^0$ with respect to the variable $x$ and $\ell_{\psi_1}^\top d\psi_1(\bar x(\cdot))$ given by (\ref{pdv1}),
and the Hamiltonian function is given by
\begin{align}\label{n340}
H_2(t,x,p,u,l)=p(f(t,x,u))+lf^0(t,x,u),\quad\forall\, (t,x,p,u,l)\in[0,T]\times T^*M\times U\times \mathbb R.
\end{align}

\end{theorem}
The corresponding second order necessary condition is stated as follows.

\begin{theorem}\label{n310}
Assume all the assumptions in Theorem \ref{n300} hold. The maps $f(t,\cdot,u)$ and $f^0(t,\cdot,u)$ are $C^2$ for all $(t,u)\in[0,T]\times U$.  The maps $h,\psi_1$ and $\psi_2$ are $C^2$. Furthermore, there exists a positive constant $L$ such that the first line of (\ref{116}) holds for $f$ and  $f^0$,  and for $\varphi=h,\psi_1,\psi_2$, the following relation holds
\begin{align*}
|\nabla\varphi(x)-L_{\hat x x}\nabla\varphi(\hat x)|\leq L\rho(x,\hat x),
\end{align*}
 where $x,\hat x\in M$ satisfy $\rho(x,\hat x)\leq\min\{i(x),i(\hat x)\}$. Then, for any $(u(\cdot), V)\in\mathcal U\times T_{\bar x(0)}M$ satisfying
\begin{align}\label{n330}
\int_0^T\left(\nabla_xf^0[t]\left(X_{u,V}(t)\right)+f^0(t,\bar x(t), u(t))-f^0[t]\right)dt\leq 0
\\\label{n331} \nabla\psi_1(\bar x(0))(V)=\nabla\psi_2(\bar x(T))(X_{u,V}(T))=0,
\end{align}
where $X_{u,V}(\cdot)$ is the solution to (\ref{n14}),
there exists $\hat\ell=(\hat\ell_0,\hat\ell_{\psi_1},\hat\ell_{\psi_2})\in\mathbb R^{1+k_1+k_2}\setminus\{0\}$ satisfying (\ref{n335}), (\ref{n304}) and (\ref{n305}) with $\ell $ replaced by $\hat\ell$, and $\hat\ell_0=0$ if ``$\leq$'' in (\ref{n330}) is ``$<$'', such that
\begin{align}\label{n315}\begin{array}{ll}
&\int_0^T\left\{\nabla_x^2 H_2\{t\}^{\hat\ell}(X_{u,V}(t),X_{u,V}(t))+2\Big(\nabla_x H_2(t,\bar x(t),p^{\hat\ell}(t),u(t),\hat\ell_0)\right.
\\& \left.-\nabla_x H_2\{t\}^{\hat\ell}\Big)(X_{u,V}(t))-R(\tilde{p}_1^{\hat\ell}(t),X_{u,V}(t),f[t],X_{u,V}(t))\right\} dt
\\&+\Big(\hat\ell_0\nabla^2h(\bar x(T))+\hat\ell_{\psi_2}^\top\nabla^2\psi_2(\bar x(T))\Big)(X_{u,V}(T),X_{u,V}(T))
\\&+ \hat\ell_{\psi_1}^\top\nabla^2\psi_1(\bar x(0))(V,V)\leq 0,
\end{array}\end{align}
where $p_1^{\hat\ell}$ is the solution to (\ref{n305}) with $\ell$ replaced by $\hat\ell$, $\tilde{p}_1^{\hat\ell}$ is the dual vector of $p_1^{\hat\ell}$,  $H_2$ is defined in (\ref{n340}), and we use the notation
\begin{align*}
\{t\}^{\hat\ell}=(t,\bar x(t),p_1^{\hat\ell}(t),\bar u(t),\hat\ell_0)
\end{align*}
for abbreviation.
\end{theorem}

\begin{rem}\label{n360qqq}
In problem $(P_2)$, when $h\equiv0$, $\psi_1(x)=exp_{x_1}^{-1}x$ and $\psi_2(x)=exp_{x_2}^{-1}x$,
where $x_1, x_2\in M$ are fixed, and  $exp_{x_i}^{-1}$ ($i=1,2$) is the inverse of the expenential map  at $x_i\in M$ $exp_{x_i}$ (see Section \ref{a2}), problem $(P_2)$ is reduced to the case that the state is fixed at the initial and final time. Especially, when $M=\mathbb R^n$, $\exp_{x_i}^{-1}x$ is reduced to $x-x_i$ for each $x\in\mathbb R^n$.
\end{rem}

\medskip
{\it Proof of Theorem \ref{n300}.}
\quad First, we shall transform problem $(P_2)$ into the form of $(P)$. Given an admissible pair $(x(\cdot), u(\cdot))$ (i.e. it is subject to (\ref{25}), $u(\cdot)\in\mathcal U$, $\psi_1(x(0))=0$ and $\psi_2(x(T))=0$), we introduce another state variable $x^0(t)=\int_0^tf^0(s,x(s),u(s))ds$, the problem $(P_2)$ can be represented as
\begin{description}
\item[$(\tilde P_2)$]Minimize $x^0(T)+h(x(T))$ subject to
\begin{align}\label{n302}
\left\{\begin{array}{l}
\left(\begin{array}{c}\dot x^0(t)
\\\dot x(t)\end{array}\right)=\left(\begin{array}{c}f^0(t,x(t),u(t))\\ f(t,x(t),u(t))\end{array}\right),
\\[3mm](x^0(0),\psi_1(x(0)))=0,\quad \psi_2(x(T))=0.
\end{array}
\right.
\end{align}
\end{description}
Denote by $\bar x^0(t)=\int_0^tf(s,\bar x(s),\bar u(s))ds$. Then,  $(\bar x^0(\cdot),\bar x(\cdot),\bar u(\cdot))$ is an optimal pair for problem $(\tilde P_2).$
By Theorem \ref{f}, there exists $\ell=(\ell_0,\ell_1,\ell_{\psi_1},\ell_{\psi_2})\in\mathbb R^{1+1+k_1+k_2}\setminus\{0\}$ such that
$
\ell_0\leq 0
$,
and
\begin{align*}
\max_{u\in U}H^e(t,\bar x^0(t),\bar x(t),p^{0\ell}(t),p_1^\ell(t),u)
=H^e(t,\bar x^0(t),\bar x(t),p^{0\ell}(t),p_1^\ell(t),\bar u(t)),
\end{align*}
for almost all $t\in[0,T]$, where $(p^{0\ell},p^\ell)$ is the sulotion to
 the following
 dual system of (\ref{n302}):
\begin{align*}
\left\{\begin{array}{ll}\dot {p^{0\ell}}(t)=0,& a.e. \,t\in(0,T),
\\ \nabla_{\dot{\bar x}(t)}p_1^\ell=-\nabla_xf[t](p_1^\ell(t),\cdot)-p^{0\ell}(t)d_xf^0[t],&a .e. \,t\in (0,T),
\\  p^{0\ell}(0)=-\ell_1,\;p^\ell(0)=-\ell_{\psi_1}^\top d\psi_1(\bar x(0)),
\\ p^{0\ell}(T)=\ell_0,\;p_1^\ell(T)=\ell_0 d h(\bar x(T))+\ell_{\psi_2}^\top d\psi_2(\bar x(T)),
\end{array}\right.
\end{align*}
and the extended Hamiltonian function is defined by
\begin{align*}
H^e(t,x^0,x,p^0,p,u)
=p^0f^0(t,x,u)+p(f(t,x,u)),
\end{align*}
for all $(t,x^0,x,p^0,p,u)\in[0,T]\times T^*(\mathbb R\times M)\times U$. From the above relations we conclude the proof. $\Box$

\medskip
{\it Proof of Theorem \ref{n310}.}\quad  By applying Theorem  \ref{n295} to problem $(\tilde P_2)$ in the proof of Thoerm \ref{n300}, we conclude the proof.  $\Box$

\subsection{An Example}

In this subsection, we will consider the curves, which connect two fixed points on a Riemannian manifold $(M,g)$, and are subject to some restrictions. We would apply Theorem \ref{f} and Theorem \ref{n295} to characterise the shortest one among all this curves.

Given any two points $y_{0}, y_{1} \in M$ and a bounded domain $D\subset M$ such that $y_0,y_1\in D$. By the completeness of the Riemannian manifold $(M,g)$, there exist smooth vector fields $f_1,\cdots,f_m$ with compact support such that
\begin{equation}\label{e9}
span\left\{\left.f_{1}\right|_{\overline{D}}, \cdots,\left.f_{m}\right|_{\overline{D}}\right\}=\left\{\left.X\right|_{\overline{D}} ; \quad X \in T M\right\},
\end{equation}
where $\overline D$ is the closure of $D$.
For more details, please see \cite[Example 4.2]{cdz}. Denote by $U\equiv \{(u_1,\cdots,u_m)\in\mathbb R^m; u_1\geq 0\}$.

Consider the following control system
\begin{equation}\label{e1}
\begin{cases}\dot y_u(t)=\sum_{i=1}^m u_i(t)f_i(y_u(t)), \;a.e.\,t\in(0,T),
\cr y_u(0)=y_0,\;y_u(T)=y_1,
\end{cases}
\end{equation}
with the control restriction
\begin{align}\label{e2}
u_1(t)\geq 0,\;(u_2(t),\cdots,u_m(t))\in \mathbb R^{m-1}\quad a.e.\,t\in[0,T].
\end{align}
A control $u(\cdot)$ determines the direction of the corresponding curve $y_u(\cdot)$. Set $\mathcal U=\{u=(u_1,\cdots,u_m)^\top:[0,T]\to \mathbb R^m \;\textrm{is measurable;}\, u_1(t)\geq 0\;\textrm{a.e.}\,t\in[0,T]\}$. Denote the set of admissible controls by
$$\begin{array}{ll}\mathcal C_{ad}\equiv&\{u\in\mathcal U;
\textrm{corresponding to control } \, u(\cdot),
\textrm{(\ref{e1}) admits a solution } y_u(\cdot)\}.\end{array}$$
Given $u(\cdot)\in\mathcal C_{ad}$, the associated cost functional is given by
$$
J(u(\cdot))=\frac{1}{2}\int_0^T|\sum_{i=1}^mu_i(t)f_i(y_u(t))|^2dt.
$$
 Set
$\ell(u(\cdot))\equiv\int_0^T|\dot y_u(t)|dt$. Then $\ell(u(\cdot))$ is the length of the curve $y_u(\cdot)$. Analogous to \cite[Proposition 17, p.126]{p1}, we obtain that, if $\bar u(\cdot)\in\mathcal C_{ad}$  minimises $J$ over $\mathcal C_{ad}$, and the corresponding solution $\bar y(\cdot)$ has constant speed (i.e. $|\dot{\bar y}(t)|\equiv $ a positive constant, $\forall \,t\in[0,T]$),  then it also minimises $\ell(\cdot)$ over $\mathcal C_{ad}$. Thus, the problem
$\min_{u(\cdot)\in\mathcal C_{ad}}J(u(\cdot))$ is in fact to find the shortest curve, which is subject to restriction (\ref{e2}) and connects $y_0$ and $y_1$.

\begin{exl}\label{e3}
Assume that $\bar u(\cdot)=(\bar u_1(\cdot),\cdots,\bar u_m(\cdot))^\top\in\mathcal C_{ad}$ satisfies $J(\bar u(\cdot))=\\ \min_{u(\cdot)\in\mathcal C_{ad}}J(u(\cdot))$ and $|\dot{\bar y}(t)|\equiv $ a positive constant a.e.  $t\in[0,T]$, where $\bar y(\cdot)$ is the  corresponding solution to (\ref{e1}).  By Theorem \ref{n300},  there exists
$(\psi_0,\psi_1)\in\Big((-\infty,0]\times T_{y_1}^*M\Big)\setminus\{0\}$ such that
\begin{align}\label{n400}\begin{array}{ll}
&\sum_{i=1}^m\psi(t)(f_i(\bar y(t)))\bar u_i(t)+\frac{1}{2}\psi_0\Big|\sum_{i=1}^m\bar u_i(t)f_i(\bar y(t))\Big|^2
\\ =&\max\Big\{\sum_{i=1}^m\psi(f_i(\bar y(t))) u_i+\frac{1}{2}\psi_0\Big|\sum_{i=1}^m u_if_i(\bar y(t))\Big|^2; u_1\geq 0, u_2,\cdots,u_m\in\mathbb R\Big\},
\end{array}\end{align}
for almost all $t\in[0,T]$,
where $\psi(\cdot)$ is a covector field along $\bar y(\cdot)$, and satisfies
\begin{equation}\label{e13}
\begin{cases}\nabla_{\dot{\bar y}(t)}
\psi=-\sum_{i=1}^m\bar u_i(t)\nabla f_i(\bar y(t))(\psi(t),\cdot)-\psi_0\sum_{i=1}^m\bar u_i(t)\nabla f_i(\bar y(t))(\tilde{\dot{\bar y}}(t),\cdot),
\cr\quad\quad\quad\quad\quad\quad\quad\quad\quad\quad\quad\quad\quad\quad\quad\quad\quad\quad\quad\quad\quad\quad\quad\quad\quad a.e.\,t\in[0,T),
\cr \psi(T)=\psi_1,
\end{cases}
\end{equation}
where $\tilde{\dot{\bar y}}(t)$ is the dual covector of $\dot{\bar y}(t)$.
The maximum principle (\ref{n400}) implies
\begin{align}\label{n500}
&\psi(t)(f_i(\bar y(t)))+\psi_0\langle f_i(\bar y(t),\dot{\bar y}(t))\rangle=0,\quad a.e.\,t\in(0,T),\; i=2,\cdots,m.
\end{align}

To figure out what $\bar u_1(\cdot)$ is, for fixed $t\in[0,T]$, we set
\begin{align*}
&h_t(u_1)
\\=&\frac{1}{2}\psi_0|f_1(\bar y(t))|^2u_1^2+\Big(\psi_0\Big\langle f_1(\bar y(t)),\sum_{i=2}^m\bar u_i(t)f_i(\bar y(t))\Big\rangle+\psi(t)(f_1(\bar y(t)))\Big)u_1
\\&+\sum_{i=2}^m \psi(t)(f_i(\bar y(t)))\bar u_i(t)+\frac{1}{2}\psi_0\sum_{i,j=2}^m\bar u_i(t)\bar u_j(t)\langle f_i(\bar y(t)),f_j(\bar y(t))\rangle,
\end{align*} for all $u_1\geq 0$.  By (\ref{e13}), for almost every $t\in(0,T)$,  $h_t(\cdot)$ attains its maximum at $\bar u_1(t)$.

If $\psi_0=0$, we obtain from the nontriviality of $(\psi_0,\psi_1)$ that $\psi_1\neq 0$, and consequently $\psi(t)\neq 0$ for all $t\in[0,T]$. We obtan from (\ref{n500}) and (\ref{e9}) that $\psi(t)(f_1(\bar y(t)))\neq 0$ for all $t\in[0,T]$.
 Then we have
\begin{align*}
\bar u_1(t)=0\;a.e. \,t\in(0,T)\quad\textrm{and}\quad\psi(t)(f_1(\bar y(t)))<0\quad \forall\, t\in[0,T].
\end{align*}

If $\psi_0<0$, we have
\begin{align*}
\bar u_1(t)=\max\left\{0, \frac{-\psi_0\langle f_1(\bar y(t)),\sum_{i=2}^m\bar u_i(t)f_i(\bar y(t))\rangle-\psi(t)(f_1(\bar y(t)))}{\psi_0|f_1(\bar y(t))|^2}\right\},\quad a.e. \,t\in(0,T).
\end{align*}
For the case that
$E\stackrel{\triangle}{=}\{t\in[0,T]; \bar u_1(t)>0\}$ is of positive measure, we obtain (\ref{n500}) with $i=1,\cdots,n$, which implies that the Lagrange multipler $(\psi_0,\psi_1)$ ( satisfying (\ref{n400})) is normal (see Definition \ref{n380} ), and   unique ( up to a positive factor), and $\bar y(\cdot)$ satisfies $\nabla_{\dot{\bar y}(t)}\dot{\bar y}=0$ for almost all $t\in E$.
 For the detailed argument, please see \cite[Example 4.3]{cdz}.

Then, we are going to seek the second order necessary condition. By Theorem \ref{n310},
for any $u(\cdot)=(u_1(\cdot),\cdots, u_m(\cdot))^\top\in\mathcal U$ satisfying
\begin{equation}\label{e8}
\left\{\begin{array}{ll}\nabla_{\dot{\bar y}(t)}X_u=\sum_{i=1}^m\bar u_i(t)\nabla_{X_u(t)}f_i+\sum_{i=1}^m(u_i(t)-\bar u_i(t))f_i(\bar y(t)),& a.e.\,t\in(0,T),
\\ X_u(0)=0,\qquad X_u(T)=0,
\end{array}\right.
\end{equation}
and
\begin{align}\label{n710}\begin{array}{ll}
&\int_0^T\Big(-\langle X_u(t),\nabla_{\dot{\bar y}(t)}\dot{\bar y}\rangle-\sum_{i=1}^m(u_i(t)-\bar u_i(t))\langle f_i(\bar y(t)),\dot{\bar y}(t)\rangle
\\&+\frac{1}{2}|\sum_{i=1}^mu_i(t)f_i(\bar y(t))|^2-\frac{1}{2}|\sum_{i=1}^m\bar u_i(t)f_i(\bar y(t))|^2\Big)dt\leq 0,
\end{array}\end{align}
there exists another $(\hat\psi_0,\hat\psi_1)\in (-\infty,0]\times T_{\bar y(T)}^*M\setminus\{0\}$ satisfying (\ref{n400}) and (\ref{e13}), with $\psi_0,\psi_1,\psi(\cdot)$ replaced respectively by $\hat\psi_0,\hat\psi_1$ and $\hat\psi(\cdot)$, and $\hat\psi_0=0$ when ``$\leq$'' in (\ref{n710}) is ``$<$'',  such that the following inequality holds:

\begin{align}\label{e15}\begin{array}{ll}
&\int_0^T\left\{\sum_{i=1}^m\bar u_i(t)\nabla^2f_i(\bar y(t))(\hat\psi(t),X_u(t),X_u(t))\right.
\\&+\hat\psi_0\sum_{i,j=1}^m \bar u_i(t)\bar u_j(t)\nabla^2f_i(\bar y(t))(\tilde f_j(\bar y(t)),X_u(t),X_u(t))
\\&+\hat\psi_0\sum_{i,j=1}^m\bar u_i(t)\bar u_j(t)\langle\nabla_{X_u(t)}f_i,\nabla_{X_u(t)}f_j\rangle
\\&+2\Big(\sum_{i=1}^m(u_i(t)-\bar u_i(t))\nabla f_i(\bar y(t))(\hat\psi(t),X_u(t))
\\&+\hat\psi_0\sum_{i,j=1}^m(u_i(t)u_j(t)-\bar u_i(t)\bar u_j(t))\nabla f_i(\bar y(t))(\tilde f_j(\bar y(t)),X_u(t))\Big)
\\& \left.-R(\tilde{\hat\psi}(t),X_u(t),\dot{\bar y}(t),X_u(t))\right\}dt\leq 0,\end{array}
\end{align}
where $\tilde{\hat\psi}(t)$ is the dual vector of $\hat\psi(t)$ ($t\in[0,T]$), and $\tilde f_j(\bar y(t))$
 is the dual covector of $ f_j(\bar y(t))$.
\end{exl}

\setcounter{equation}{0}
\section{Proof of the Main Results}
\def\theequation{4.\arabic{equation}}
This section is split into three parts. In the first subsection, we give some lemmas related to Liapounoff's Theorem. In the second subsection, we shall prove Theorem \ref{n295} first, and show the sketch of the proof of Theorem \ref{f}, and the proofs of Theorem \ref{n770} and Theorem \ref{n790} are given in the last subsection.

\subsection{Some lemmas}
By making a little revision to the proofs of \cite[Lemma 3.7, p. 143]{ly} and \cite[Corollary 3.8, p. 144]{ly}, we have the following  Liapounoff's type lemma.
\begin{lem}\label{n90}
Asume $h\in C([0,T];L^1(0,T;\mathbb R^k))$ ($k\in\mathbb N$). Fix $\epsilon>0$. For any $\rho\in[0,1]$, there exist measurable subset $E_\rho\subset[0,T]$ and $R:[0,T]\times [0,1]\to\mathbb R^k$   such that
\begin{align}\label{n91}
&\rho\int_0^T h(t,s)ds=\int_{E_\rho}h(t,s)ds+R(t,\rho),\quad\forall \,t\in[0,T],
\\
\label{n94} &|R(t,\rho)|\leq\epsilon,\quad\forall \,t\in[0,T],
\\ \label{n95}& |E_\rho|=\rho T,  E_\rho\subseteq E_{\hat\rho},\quad\textrm{if}\;0\leq \rho\leq \hat\rho\leq 1,
\end{align}
and $R(t, \cdot)$ is continuous on $[0,1]$ for each $t\in[0,T]$. Furthermore, fix  measurable  subset $E\subseteq [0,T]$, for any $\rho\in[0,1]$, there exist measurable subset $A_\rho$ and $B(\cdot,\rho, E):[0,T]\to\mathbb R^k$  such that
\begin{align}\label{n97}
\rho\int_E h(t,s)ds=\int_{A_\rho}h(t,s)ds +B(t,\rho,E),\quad\forall\, t\in[0,T],
\\\label{n98}  |B(t,\rho,E)|\leq \epsilon,\quad\forall\, t\in[0,T],
\\ \label{n99} A_\rho\subseteq E,\;|A_\rho|=\rho |E|,\quad A_\rho\subseteq A_{\hat\rho},\quad\textrm{if }\; 0\leq \rho\leq\hat\rho\leq 1.
\end{align}
\end{lem}

{\it Proof.}\quad First, we shall prove (\ref{n91}) - (\ref{n95}). Fix any $\epsilon>0$, there exists $\delta>0$ such that
\begin{align*}
\int_0^T|h(t,s)-h(\hat t,s)|ds\leq \frac{1}{5}\epsilon,\quad|t-\hat t|\leq \delta,\;t,\hat t\in[0,T].
\end{align*}
Then, there exist $0=t_0\leq t_1<\cdots<t_l=T$ $(l\in\mathbb N)$ such that $|t_i-t_{i+1}|\leq \delta$ for $i=0,1,\cdots,l-1$. Deonte by $g(s)=(h(t_0,s),h(t_1,s),\cdots,h(t_l,s))^\top$ for $s\in[0,T]$. Then $g(\cdot)\in L^1(0,T;\mathbb R^{ (l+1)\times k})$, and consequently  there exists a $\mathbb R^{ (l+1)\times k}$ valued  simple function $S_g(\cdot)=\sum_{j=1}^p I_{F_j}(\cdot)g_j$ with $p\in\mathbb N$, $\cup_{j=1}^pF_j=[0,T]$, $F_i\cap F_j=\emptyset$ when $i\neq j$, and $g_j\in\mathbb R^{(l+1)\times k}$, such that
\begin{align*}
\int_0^T|g(s)-S_g(s)|ds\leq \frac{1}{5}\epsilon.
\end{align*}
For any $\rho\in(0,1]$,  there exist $E_\rho^1,\cdots, E_\rho^p$  satisfying
\begin{align*}
E_\rho^j\subseteq F_j,\quad |E_\rho^j|=\rho|F_j|,&\quad j=1,\cdots,p,
\\ E_{\hat\rho} ^j\subseteq E_\rho^j,\;|E_{\hat \rho}^j|=\hat\rho|F_j|,\quad\textrm{if}\;0\leq \hat\rho\leq \rho,&\quad j=1,\cdots,p.
\end{align*}
Set $E_\rho=\cup_{j=1}^p E_\rho^j $. Then $E_\rho$ fulfills  (\ref{n95}).
Consequently, we have
\begin{align*}
&\rho\int_0^Tg(s)ds\nonumber
\\=&\rho\int_0^TS_g(s)ds+\rho\int_0^T(g(s)-S_g(s))ds\nonumber
%\\=&\rho\sum_{j=1}^p|F_j|g_j+\rho\int_0^T(g(s)-S_g(s))ds\nonumber
%\\=&\sum_{j=1}^p\int_{E_\rho^j}g_jds+\rho\int_0^T(g(s)-S_g(s))ds\nonumber
%\\=&\int_{E_\rho} S_g(s)ds+\rho\int_0^T(g(s)-S_g(s))ds\nonumber
\\=& \int_{E_\rho}g(s)ds+R(\rho,g)
\end{align*}
where $R(\rho,g)=(R_0(\rho,g),\cdots,R_l(\rho,g))^\top=\int_{E_\rho}(S_g(s)-g(s))+\rho\int_0^T(g(s)-S_g(s))ds$ is continuous  with respect to $\rho$, and satisfies
$|R(\rho,g)|\leq \|S_g-g\|_{L^1(0,T;\mathbb R^{(l+1)\times k})}(1+\rho)\leq\frac{2}{5}\epsilon$.

For any $t\in[0,T]$, there exists $t_j$ ($j\in\{0,\cdots, l\}$)such that $|t-t_j|\leq \delta$. Then, we have
\begin{align*}
&\rho\int_0^T h(t,s)ds
\\=&\rho\int_0^T h(t_j,s)ds+\rho\int_0^T (h(t,s)-h(t_j,s))ds
\\=&\int_{E_\rho}h(t_j,s)ds+R_j(\rho,g)+\rho\int_0^T (h(t,s)-h(t_j,s))ds
\\=&\int_{E_\rho}h(t,s)ds+R(t,\rho),
\end{align*}
where $R(t,\rho)=\int_{E_\rho}(h(t_j,s)-h(t,s))ds+R_j(\rho,g)+\rho\int_0^T (h(t,s)-h(t_j,s))ds$ is continunous with respect to $\rho$ when $t$ is fixed, and satisfies (\ref{n94}).  (\ref{n91}) follows immediately.

Then, we shall prove (\ref{n97}) - (\ref{n99}). For each $j=1,\cdots,p$, there exist a measurable set $A_\rho^j\subset E\cap F_j $ satisfying $|A
_\rho^j|=\rho |F_j\cap E|$ and $A_\rho^j\subseteq A_{\hat\rho}^j$ when $0\leq\rho\leq\hat\rho\leq 1$. Set $A_\rho=\cup_{j=1}^p A_\rho^j$.  Then $A_\rho$ satisfies (\ref{n99}), and we have
\begin{align*}
&\rho\int_E g(s)ds
%\\=&\rho\int_E S_g(s)ds+\rho\int_E (g(s)-S_g(s))ds
\\= &\sum_{j=1}^p\int_{A_\rho^j}I_{F_j}(s) g_jds+\rho\int_E (g(s)-S_g(s))ds
\\=&\int_{A_\rho} S_g(s)ds+\rho\int_E (g(s)-S_g(s))ds
\\=&\int_{A_\rho} g(s)ds+B(\rho,E),
\end{align*}
where $B(\rho,E)=(B_0(\rho,E),\cdots, B_l(\rho,E))^\top=\int_{A_\rho}(S_g(s)-g(s))ds+\rho\int_E (g(s)-S_g(s))ds$ is continuous with respect to $\rho$.

For any $t\in[0,T]$, there exists $t_j$ such that $|t-t_j|\leq \delta$. Then, we have
\begin{align*}
\rho\int_Eh(t,s)ds
=\int_{A_\rho}h(t,s)ds+B(t,\rho,E),
\end{align*}
where $B(t,\rho,E)=\int_{A_\rho}(h(t_j,s)-h(t,s))ds+B_j(\rho,E)+\rho\int_E(h(t,s)-h(t_j,s))ds$ satisfies (\ref{n98}). The proof is concluded.

$\Box$

By the induction argument, we immediately obtain the following corollary.
\begin{cor}\label{n130}Assume $h_1,\cdots,h_l\in C([0,T];L^1(0,T;\mathbb R^k))$ $(k,l\in\mathbb N)$. Fix  $\epsilon>0$. Then,  for any $\vec\rho=(\rho_1,\cdots,\rho_l)\in\mathbb R^l$ with
\begin{align}\label{n102}
\sum_{j=1}^l\rho_j=1\quad\textrm{ and} \quad \rho_j\geq 0\quad\textrm{for} \;j=1,\cdots,l,
\end{align} there exist measurable subsets $E_{\vec\rho}^1,\cdots,E_{\vec\rho}^l$ of $[0,T]$ such that $|E_{\vec\rho}^i|=\rho_i T$ for $i=1,\cdots,l$, $E_{\vec\rho}^i\cap E_{\vec\rho}^j=\emptyset$ when $i\neq j$, $\cup_{i=1}^lE_{\vec\rho}^i=[0,T]$, and
\begin{align}\label{n771}
\sum_{i=1}^l\rho_i\int_0^Th_i(t,s)ds=\sum_{i=1}^l\int_{E_{\vec\rho}^i}h_i(t,s)ds+R(t,\vec\rho),\quad\forall\, t\in[0,T],
\end{align}
where $|R(t,\vec\rho
)|\leq \epsilon$ for all $t\in[0,T]$. Furthermore, there exists a positive constant $C$ such that $\sum_{i=1}^lH(E_{\vec\rho}^i,E_{\vec{\hat\rho}}^i)\leq C\sum_{i=1}^l|\rho_i-\hat\rho_i|$ for all $\vec\rho=(\rho_1,\cdots,\rho_l)$ and $\vec{\hat\rho}=(\hat\rho_1,\cdots,\hat\rho_l)$ satisfying  (\ref{n102}), where $H$ is the Hausdorff metric.
\end{cor}

For the sets $E_{\vec \rho}^1,\cdots, E_{\vec \rho}^l$ in Corollary \ref{n130}, we have another choice, such that  the rest term $R$ in (\ref{n771}) has different properties:

\begin{lem}\label{l6} Assume $h_i\in C([0,T]; L^1(0,T;\mathbb R^k))$ ($i=1,\cdots,l$) with $k,l\in \mathbb N$. Fix $\vec\rho=(\rho_1,\cdots,\rho_l)^\top\in\mathbb R^l$ satisfying (\ref{n102}). Then, given any $\epsilon>0$, there exist disjoint subsets $E_1,\cdots, E_l$ of $[0,T]$ such that
\begin{align*}
&\sum_{i=1}^l\rho_i\int_0^Th_i(t,s)ds=\sum_{i=1}^l\int_{E_i}h_i(t,s)ds+R(t,\epsilon),\;\forall\,t\in[0,T],
\\ & \cup_{i=1}^l E_i=[0,T],\quad |E_i|=\rho_i T,\;i=1,\cdots,l,
\end{align*}
where the term $R$ satisfies
$|R(t,\epsilon)|<\epsilon$ for all    $ t\in[0,T]$, and
$R(T,\epsilon)=0$.
\end{lem}

{\it Proof.} For any $\epsilon>0$, there exists $\delta>0$ such that
\begin{align*}
\int_0^T |h_i(t,s)-h_i(\hat t,s)|ds<\frac{\epsilon}{2(l+1)},\quad t,\hat t\in[0,T],\;|t-\hat t|\leq \delta,\;i=1,\cdots,l.
\end{align*}
Choose $0=t_0<t_1<\cdots<t_N=T$ such that $|t_i-t_{i-1}|\leq\delta$ for $i=1,\cdots,N.$
For $i=1,\cdots,l$, set $g_i(s)=(h_i(t_0,s),\cdots,h_i(t_N,s))^\top$ for $s\in[0,T]$. Then, $g_i(\cdot)\in L^1(0,T;\mathbb R^{(N+1)\times k})$. Applying  Liapounoff's Theorem to $g_1(\cdot),\cdots, g_l(\cdot)$, there exist mutual disjoint subsets $E_1,\cdots, E_l$ of $[0,T]$ such that
\begin{align}\label{l8}
\sum_{i=1}^l\rho_i\int_0^Tg_i(s)ds=\sum_{i=1}^l\int_{E_i}g_i(s)ds,
\end{align}
$|E_i|=\rho_iT$ for $i=1,\cdots,l$ and $\cup_{i=1}^lE_i=[0,T]$.
For any $t\in[0,T]$, there eixsts $t_i$ ($i=0,1,\cdots,N-1$) such that $t\in[t_i,t_{i+1}]$. Then, we have
\begin{align*}
&\Big|\sum_{j=1}^l\rho_j\int_0^Th_j(t,s)ds-\sum_{j=1}^l\int_{E_j}h_j(t,s)ds\Big|
\\\leq&\Big|\sum_{j=1}^l\rho_j\int_0^Th_j(t_i,s)ds-\sum_{j=1}^l\int_{E_j}h_j(t_i,s)ds\Big|+\sum_{j=1}^l\rho_j\int_0^T |h_j(t,s)-h_j(t_i,s)|ds
\\&+\sum_{j=1}^l\int_{E_j} |h_j(t,s)-h_j(t_i,s)|ds
\\<&\epsilon,
\end{align*}
which, together with (\ref{l8}), completes the proof.   $\Box$

\subsection{Proofs of Theorems \ref{f}--\ref{n295}}
Since the idea of the  proof of Theorem \ref{f} is similar to that of Theorem \ref{n295}, we will prove Theorem \ref{n295} in detail, and give the sketch of the proof of Theorem \ref{f} at the end of this subsection.
\begin{lem}\label{n65}
 Given any index set $I \subseteq\{0,1, \cdots, j\},$ we denote by $\Phi_{I}=(\bar{\phi}_{0}, \bar{\phi}_{1}, \cdots, \bar{\phi}_{j}, \\
\psi)^\top$ if $k>0$, and by $\Phi_I=(\bar\phi_0,\bar \phi_1,\cdots,\bar\phi_j)^\top$ if $k=0$, where $\bar{\phi}_{i}=\phi_{i}$
if $i \in I,$ and $\bar{\phi}_{i}=0$ if $i \notin I .$ Assume conditions $(C1)-(C3)$ hold and $(\bar x(\cdot), \bar u(\cdot))$ is an optimal pair for Problem $(P)$.
We also assume $u(\cdot)\in\mathcal U$ is a  Pontryagin's type critical  direction (see Definition \ref{n40}), and  $(u(\cdot),V)\in \mathcal U\times T_{\bar x(0)}M$ satisfies (\ref{n53}). For $\varphi=\phi_0,\cdots,\phi_j,\psi$, denote by
\begin{equation}\label{n55}\begin{array}{ll}
D^2\varphi(u,V)\equiv&\nabla_1^2\varphi(\bar x(0),\bar x(T))(V,V)+2\nabla_2\nabla_1\varphi(\bar x(0),\bar x(T))(V,X_{u,V}(T))
\\&+\nabla_2^2\varphi(\bar x(0),\bar x(T))(X_{u,V}(T),X_{u,V}(T)),
\end{array}\end{equation}
 where for the vector valued function $\psi$,
\begin{align*}
&\nabla_i^2\psi(\bar x(0),\bar x(T))(X,Y)
\\=&[\nabla_i^2\psi_1(\bar x(0),\bar x(T))
(X,Y), \cdots, \nabla_i^2\psi_k(\bar x(0),\bar x(T))
(X,Y)]^\top, \,i=1,2,
\\& \nabla_2\nabla_1\psi(\bar x(0),\bar x(T))(X,
Y)
\\=&[\nabla_2\nabla_1\psi_1(\bar x(0),\bar x(T))(X,
Y),\cdots,\nabla_2\nabla_1\psi_k(\bar x(0),\bar x(T))(X,
Y)]^\top,
\end{align*}
for all $X,Y\in\mathcal X(M)$,
and $X_{u,V}(\cdot)$ is the solution to  (\ref{n14}).
Set
\begin{equation}\label{n51}
\begin{array}{ll}\mathcal K_{u,V}\equiv& \{\nabla_1\Phi_{I_0^{\prime\prime}}(\bar x(0),\bar x(T))(W)+\nabla_2\Phi_{I_0^{\prime\prime}}(\bar x(0),\bar x(T))(Y_{u,\sigma}^{\lambda,W}(T))
\\&+\frac{1}{2}D^2\Phi_{I_0^{\prime\prime}}(u,V)|
 (\sigma(\cdot),\lambda,W)\in\mathcal U\times (0,+\infty)\times T_{\bar x(0)}M\},\end{array}
\end{equation}
where the covariant derivative of a vector-valued function is given by (\ref{pdv}), and   $Y_{u,\sigma}^{\lambda,W}(\cdot)$ is the  solution to
\begin{equation}\label{n20}
\left\{\begin{array}{lll}
\nabla_{\dot{\bar x}(t)}Y_{u,\sigma}^{\lambda,W}&=&\nabla_xf[t](\cdot,Y_{u,\sigma}^{\lambda,W}(t))+\lambda (f(t,\bar x(t),\sigma(t))-f[t])
\\&&+\frac{1}{2}\nabla_x^2f[t](\cdot,X_{u,V}(t),X_{u,V}(t))+\nabla_xf(t,\bar x(t),u(t))(\cdot,X_{u,V}(t))
\\&&-\nabla_xf[t](\cdot,X_{u,V}(t))-\frac{1}{2}R(\cdot,X_{u,V}(t),f[t],X_{u,V}(t)),\;a.e. t\in(0,T),
\\ Y_{u,\sigma}^{\lambda,W}(0)&=&W.
\end{array}
\right.
\end{equation}
Then, $\mathcal K_{u,V}$ is  a convex subset of $\mathbb R^{1+j+k}$.

\end{lem}

{\it Proof.}    In fact, we only have to show that $\{(W,Y_{u,\sigma}^{\lambda,W}(T))| \sigma\in\mathcal U,\lambda>0,W\in T_{\bar x(0)}M\}\subset T_{\bar x(0)}M\times T_{\bar x(T)}M $ is convex.

Let $\left\{e_{1}, \cdots, e_{n}\right\} \subset T_{\bar {x}(0)} M$ be an orthonormal basis.  For each $t \in(0, T]$, let $\{e_1(t),\cdots,$ $e_n(t)\}$ and $\{d_1(t),\cdots,d_n(t)\}$ be given in Section \ref{psnc}.    Then, $\{e_1(t),\cdots,e_n(t)\}$ is an orthonormal basis at $T_{\bar x(t)}M$, and $\{d_1(t),\cdots,d_n(t)\}$ is the dual  basis to it.
Consequently, for $t \in[0, T]$ we can express tensors $\nabla_{x} f[t]$ and  $f(t,\bar x(t),u(t))$  respectively  by  $\nabla_{x} f[t]=\sum_{i, j=1}^{n} A_{i j}(t,\bar u(t)) e_{i}(t) \otimes d_{j}(t)$ and $f(t,\bar x(t),u(t))=\sum_{i=1}^n f^i(t,u(t))e_i(t)$, where $A_{ij}(t,\bar u(t))$ and $f^i(t,u(t))$ ($i,j=1,\cdots,n$) are defined by (\ref{n738}).
Set
$$\Theta(t;u,V)=(\theta^1(t;u,V),\cdots,\theta^n(t;u,V))^\top$$ with
$$\begin{array}{ll}&\theta^i(t;u,V)
\\=&\frac{1}{2} \nabla_{x}^{2} f[t]\left(d_i(t), X_{u, V}(t), X_{u, V}(t)\right)+\nabla_{x} f(t, \overline{x}(t), u(t))\left(d_i(t), X_{u, V}(t)\right)
\\&-\nabla_{x} f[t]\left(d_i(t), X_{u, V}(t)\right)-\frac{1}{2} R\left(e_i(t), X_{u, V}(t), f[t], X_{u, V}(t)\right),\qquad i=1,\cdots,n.\end{array}$$
Given $\sigma\in\mathcal U$ and $\lambda>0$, denote by $Y_{u,\sigma}^{\lambda,W}(t)=\sum_{i=1}^ny^{i,W}_{u,\sigma,\lambda}(t) e_i(t)$. Then, $\vec Y_{u,\sigma}^{\lambda,W}(t)\equiv(y^{1,W}_{u,\sigma,\lambda}(t),\cdots,y^{n,W}_{u,\sigma,\lambda}(t))^
\top$ ($t\in[0,T]$) solves
$$\left\{\begin{array}{l}
\dot{\vec{Y}}_{u,\sigma}^{\lambda,W}(t)=A(t,\bar u(t))\vec{Y}_{u,\sigma}^{\lambda,W}(t)+\lambda (\vec f(t,\sigma(t))-\vec f(t,\bar u(t)))+\Theta(t;u,V),\quad a.e.\,t\in(0,T],
\\ \vec Y_{u,\sigma}^{\lambda,W}(0)=\vec W,
\end{array}
\right.$$
where $\vec W\stackrel{\triangle}{=}(w_1,\cdots,w_n)^\top$ with $w_i=\langle W,e_i\rangle$ ($i=1,\cdots,n$), and $A(t,\bar u(t))$ and $\vec f(t,\sigma(t))$ are defined by (\ref{gcs20}).
Assume $\eta:[0,+\infty)\to \mathbb R^{n\times n}$ solves
$$\left\{\begin{array}{l}
\dot\eta(t)=A(t,\bar u(t))\eta(t),\quad t>0,
\\ \eta(0)=I,
\end{array}\right.
$$
 where $I\in \mathbb R^{n\times n}$ is the identity matrix. Then, we have
 \begin{align*}
\vec Y_{u,\sigma}^{\lambda,W}(t)=&\eta(t)\vec W+\lambda\int_0^t\eta(t)\eta(s)^{-1}(\vec f(s,\sigma(s))-\vec f(s,\bar u(s)))ds&
\\&+\int_0^t\eta(t)\eta(s)^{-1}\Theta(s;u,V)ds,& \forall\,t\in[0,T].
\end{align*}

 Fix any $\sigma_1,\sigma_2\in\mathcal U$, $\lambda_1,\lambda_2\in(0,+\infty)$, $W^1, W^2\in T_{\bar x(0)}M$ and $ \nu\in(0,1)$. Denote by $W^j=\sum_{i=1}^n w^j_i e_i$ and $\vec W^j=(w_1^j,\cdots,w_n^j)^\top$ for $j=1,2$. By Liapounoff's convexity theorem (see \cite[Lemma 4.2]{y} ) one can find measurable subset $E\subset[0,T] $ with measure $\frac{\nu\lambda_1}{\nu\lambda_1+(1-\nu)\lambda_2}T$  such that
 $$\begin{array}{l}
 \frac{\nu\lambda_1}{\nu\lambda_1+(1-\nu)\lambda_2}\int_0^T \eta(T)\eta(s)^{-1}\vec f(s,\sigma_1(s))ds
 \\ +\frac{(1-\nu)\lambda_2}{\nu\lambda_1+(1-\nu)\lambda_2}\int_0^T \eta(T)\eta(s)^{-1}\vec f(s,\sigma_2(s))ds
 \\=\int_0^T \eta(T)\eta(s)^{-1}\vec f(s,\sigma_{12}(s))ds,
 \end{array}$$
 where
 \begin{equation}\label{n71}
 \sigma_{12}(s)=\left\{\begin{array}{l}
 \sigma_1(s),\quad\textrm{if}\,s\in E,
 \\\sigma_2(s),\quad\textrm{if}\,s\in [0,T]\setminus E.
 \end{array}\right.
 \end{equation}
 Thus we have
 \begin{equation}\label{n72}\begin{array}{ll}
 &\nu \vec Y_{u,\sigma_1}^{\lambda_1,W^1}(T)+(1-\nu)\vec Y_{u,\sigma_2}^{\lambda_2,W^2}(T)
 \\  =&(\nu\lambda_1+(1-\nu)\lambda_2)\int_0^T\eta(T)\eta(s)^{-1}(\vec f(s,\sigma_{12}(s))-\vec{f}(s,\bar u(s)))ds
 \\& +\eta(T)(\nu \vec W^1+(1-\nu)\vec W^2)+\int_0^T\eta(T)\eta(s)^{-1}\Theta(s;u,V)ds
 \\=&\vec Y_{u,\sigma_{12}}^{\nu\lambda_1+(1-\nu)\lambda_2,\nu W^1+(1-\nu)W^2}(T),
 \end{array}\end{equation}
 which implies the convexity of $\mathcal K_{u,V}$. $\Box$

\begin{lem}\label{n66} Assume all the assumptions in Lemma \ref{n65} hold,  $k>0$, and $(u(\cdot), V)\in \mathcal U\times T_{\bar x(0)}M$ satisfies (\ref{n53}).
 Set
$$
Y\equiv (Y_0,Y_1,\cdots,Y_j)^\top,
$$
where
$$
Y_i=\begin{cases}\nabla_1\phi_i(\bar x(0),\bar x(T))(V)+\nabla_2\phi_i(\bar x(0),\bar x(T))(X_{u,V}(T)),\;i\in  I_{AO},\cr
0,\quad \quad \quad\quad\quad\quad\quad\quad\quad\quad\quad\quad\quad\quad\quad\quad\quad\quad\quad\quad\quad\;\;i\notin  I_{AO}.
\end{cases}
$$
 Set  $Z=(-\infty,0)^{1+j}-cone\{\phi(\bar x(0),\bar x(T))+Y\}$,
 where $\phi=(\phi_0,\phi_1,\cdots,\phi_j)^\top$, and $cone\,S$ is the convex cone generated by set $S$ (see \cite[p. 14]{rock} for the definition). Assume $\phi_0(\bar x(0),\bar x(T))=0$.
Then, the dimension  of the affine hull of the following set (see \cite[p. 4]{rock})
\begin{align*}
\mathcal K_{u,V}^\psi=\left\{\begin{array}{l}\nabla_2\psi(\bar x(0),\bar x(T))(Y_{u,\sigma}^{\lambda,W}(T))
\\+\nabla_1\psi(\bar x(0),\bar x(T))(W)+\frac{1}{2}D^2\psi(u,V)\end{array}\Big| \sigma\in\mathcal U, \lambda>0, W\in T_{\bar x(0)}M.\right\}
\end{align*} is bigger than or equal to one, which is denoted by $l$. Moreover,
 if there does not exist $\ell=(\ell_0,\cdots,\ell_j,
 \ell_\psi^\top)^\top\in\mathbb R^{1+j+k}\setminus\{0\}$ such that
\begin{align}\label{n70}\begin{array}{ll}
&\ell^\top\Big(\nabla_1\Phi_{I_0^{\prime\prime}}(\bar x(0),\bar x(T))(W)+\nabla_2\Phi_{I_0^{\prime\prime}}(\bar x(0),\bar x(T))(Y_{u,\sigma}^{\lambda,W}(T))+\frac{1}{2}D^2\Phi_{I_0^{\prime\prime}}(u,V)\Big)
\\ &\leq \ell^\top (z^\top\quad 0)^\top,
\end{array}\end{align}
 for all $(\sigma,\lambda,W,z)\in \mathcal U\times (0,+\infty)\times T_{\bar x(0)}M\times Z$, then there exist $(\sigma_1,\lambda_1,W^1), \cdots, \\(\sigma_{l+1},\lambda_{l+1},W^{l+1})\in \mathcal U\times (0,+\infty)\times T_{\bar x(0)}M$ and $\delta_0>0$,
 such that
 \begin{align}\label{n60}\begin{array}{l}
B_{\mathbb R^l}(\delta_0)\\\subset co\{\nabla_1\psi(\bar x(0),\bar x(T))(W^\eta)+\nabla_2\psi(\bar x(0),\bar x(T))(Y_{u,\sigma_\eta}^{\lambda_\eta,W^\eta}(T))+\frac{1}{2}D^2\psi(u,V)\}_{\eta=1}^{l+1},\end{array}
\end{align}
where $co A$ denotes the convex hull of set $A$, and $B_{\mathbb R^l}(\delta_0)$ is the closed ball in $\mathbb R^l$ with center at the origin and of radius $\delta_0$,
and if $i\in  I_0^{\prime\prime}$,
\begin{align}
\label{n61}&\nabla_1\phi_i(\bar x(0),\bar x(T))(W^\eta)+\nabla_2\phi_i(\bar x(0),\bar x(T))(Y_{u,\sigma_\eta}^{\lambda_\eta,W^\eta}(T))+\frac{1}{2}D^2\phi_i(u,V)<0,
\end{align}
 for $\eta=1,\cdots,l+1$.
\end{lem}

{\it Proof.} First,  we claim
$0\in \textrm{ri}\, \mathcal K_{u,V}^\psi$, where
ri $ A$ is the interior of set $A$ relative to its affine hull of $ A$ (see \cite[p. 44]{rock} for its detailed definition).
By contradiction, we assume it was not true.
Since $\mathcal K_{u,V}^\psi$ is convex (by Lemma \ref{n65}),  the affine hull of $\mathcal K_{u,V}^\psi$ is closed (see \cite[p. 44]{rock}), and $ri\,\mathcal K_{u,V}^\psi\neq \emptyset$ (by \cite[Theorem 6.2, p. 45]{rock}),  we obtain from  \cite[Lemma 3.1]{bern} or \cite[Theorem 11.1, p.95 \& Theorem 11.3, p.97]{rock} that,
there exists $\xi\in\mathbb R^k\setminus\{0\}$ such that
\begin{align}\label{gcs27}
\xi^\top\Big(\nabla_1\psi(\bar x(0),\bar x(T))(W)+\nabla_2\psi(\bar x(0),\bar x(T))(Y_{u,\sigma}^{\lambda,W}(T))+\frac{1}{2}D^2\psi(u,V)\Big)\leq 0,
\end{align}
for all $\sigma\in\mathcal U$, $\lambda>0$ and $W\in T_{\bar x(0)}M$.
Consequently, we have
\begin{align}\label{gcs28}\begin{array}{ll}
&(0,\xi^\top)\Big(\nabla_1\Phi_{I_0^{\prime\prime}}(\bar x(0),\bar x(T))(W)+\nabla_2\Phi_{I_0^{\prime\prime}}(\bar x(0),\bar x(T))(Y_{u,\sigma}^{\lambda,W}(T))+\frac{1}{2}D^2\Phi_{I_0^{\prime\prime}}(u,V)\Big)
\\ &\leq (0,\xi^\top) (z^\top\; 0)^\top,\end{array}
\end{align}
for all $\sigma\in\mathcal U, \lambda>0$, $W\in T_{\bar x(0)}M$ and $z\in Z$, which leads to a  contradiction.

Second, we claim that $l\geq 1$. If this assertion were not true, we have $l=0$, because $\mathcal K_{u,V}^\psi\neq\emptyset$. Consequently $\mathcal K_{u,V}^\psi=\{0\}$. Thus, for any $\alpha\in\mathbb R^k\setminus\{0\}$, (\ref{gcs27}) holds with $\xi$ replaced by $\alpha$, and then (\ref{gcs28}) holds with $\xi$ replaced by $\alpha$, which leads to a contradiction.

Third, there exist $(\tilde \sigma_1,\tilde\lambda_1,\tilde W^1), \cdots, (\tilde\sigma_{l+1},\tilde\lambda_{l+1},\tilde W^{l+1})\in \mathcal U\times (0,+\infty)\times T_{\bar x(0)}M$ and $\tilde\delta_0>0$ such that (\ref{n60}) holds with $\delta_0$ and $\{(\sigma_\eta,\lambda_\eta, W^\eta)\}_{\eta=1}^{l+1}$ replaced respectively by $\tilde\delta_0$ and $\{(\tilde\sigma_\eta,\tilde\lambda_\eta,\tilde W^\eta)\}_{\eta=1}^{l+1}$.

According to the assumption and \cite[Theorem 11.3, p.97]{rock}, we have $\mathcal K_{u,V}\cap (Z\times\{0\})\neq\emptyset$. Then, there exist $(\sigma_0,\lambda_0,W^0)\in \mathcal U\times (0,+\infty)\times  T_{\bar x(0)}M$, $\theta_0>0$ and $(z_0,z_1,\cdots,z_j)\in(-\infty,0)^{1+j}$ such that
\begin{align}\label{gcs30}\begin{array}{ll}
\nabla_1\phi_i(\bar x(0),\bar x(T))(W^0)+\nabla_2\phi_i (\bar x(0),\bar x(T))(Y_{u,\sigma_0}^{\lambda_0,W^0}(T))&
\\+\frac{1}{2}D^2\phi_i (u,V)=z_i,& \textrm{if}\; i\in I_0^{\prime\prime},
\\
\nabla_1\phi_i(\bar x(0),\bar x(T))(W^0)+\nabla_2\phi_i (\bar x(0),\bar x(T))(Y_{u,\sigma_0}^{\lambda_0,W^0}(T))+\frac{1}{2}D^2\phi_i (u,V)&
\\=z_i-\theta_0(\nabla_1\phi_i(\bar x(0),\bar x(T))(V)+\nabla_2\phi_i(\bar x(0),\bar x(T))(X_{u,V}(T))),&\textrm{if}\;i\in I_{AO}\setminus I_0^{\prime\prime},
\\ \nabla_1\phi_i(\bar x(0),\bar x(T))(W^0)+\nabla_2\phi_i (\bar x(0),\bar x(T))(Y_{u,\sigma_0}^{\lambda_0,W^0}(T))+\frac{1}{2}D^2\phi_i (u,V)
\\=z_i-\theta_0\phi_i(\bar x(0),\bar x(T)),&\textrm{if}\;i\in I_N,
\end{array}\end{align}
and
\begin{align*}
\nabla_1\psi(\bar x(0),\bar x(T))(W^0)+\nabla_2\psi(\bar x(0),\bar x(T))(Y_{u,\sigma_0}^{\lambda_0,W^0}(T))+\frac{1}{2}D^2\psi (u,V)=0.
\end{align*}
When $i\in  I_0^{\prime\prime}$, one can find $\theta_1\in(0,1)$ such that
\begin{align*}
(1-\theta_1)\Big(\nabla_1\phi_i(\bar x(0),\bar x(T))(\tilde W^\eta)+\nabla_2\phi_i (\bar x(0),\bar x(T))(Y_{u,\tilde\sigma_\eta}^{\tilde\lambda_\eta,\tilde W^\eta}(T))+\frac{1}{2}D^2\phi_i (u,V)\Big)
\\+\theta_1 \Big(\nabla_1\phi_i(\bar x(0),\bar x(T))(W^0)+\nabla_2\phi_i (\bar x(0),\bar x(T))(Y_{u,\sigma_0}^{\lambda_0,W^0}(T))+\frac{1}{2}D^2\phi_i (u,V)\Big)<0,
\end{align*}
for all $ \eta=1,\cdots,l+1$. Since $\mathcal K_{u, V}$ is convex, there exist $\{(\sigma_\eta,\lambda_\eta, W^\eta)\}_{\eta=1}^{l+1}\in \mathcal U\times (0,+\infty)\times T_{\bar x(0)}M$ such that
\begin{align*}
&(1-\theta_1)\Big(\nabla_1\Phi_{I_0^{\prime\prime}} (\bar x(0),\bar x(T))(\tilde W^\eta)+\nabla_2\Phi_{I_0^{\prime\prime}} (\bar x(0),\bar x(T))(Y_{u,\tilde\sigma_\eta}^{\tilde\lambda_\eta, \tilde W^\eta}(T))+\frac{1}{2}D^2\Phi_{I_0^{\prime\prime}} (u,V)\Big)
\\&+\theta_1 \Big(\nabla_1\Phi_{I_0^{\prime\prime}} (\bar x(0),\bar x(T))(W^0)+\nabla_2\Phi_{I_0^{\prime\prime}} (\bar x(0),\bar x(T))(Y_{u,\sigma_0}^{\lambda_0, W^0}(T))+\frac{1}{2}D^2\Phi_{I_0^{\prime\prime}} (u,V)\Big)
\\=&\nabla_1\Phi_{I_0^{\prime\prime}} (\bar x(0),\bar x(T))(W^\eta)+\nabla_2\Phi_{I
_0^{\prime\prime}} (\bar x(0),\bar x(T))(Y_{u,\sigma_\eta}^{\lambda_\eta, W^\eta}(T))+\frac{1}{2}D^2\Phi_{I_0^{\prime\prime}} (u,V)
\end{align*}
holds
for $\eta=1,\cdots,l+1$. Set  $\delta_0=(1-\theta_1)\tilde\delta_0$. Then, (\ref{n60}) and (\ref{n61}) follow.
$\Box$

\begin{lem}\label{n190}
 Assume assumptions $(C1)-(C3)$ hold,  $(\bar x(\cdot),\bar u(\cdot))$ is an optimal pair of Problem $(P)$, $k>0$,  and $\phi_0(\bar x(0),\bar x(T))=0$. Assume $(u(\cdot), V)\in \mathcal U\times T_{\bar x(0)}M$ satisfies (\ref{n53}).  Then there exists $(\ell_0,\ell_1,\cdots,\ell_j,\ell_\psi)\in\mathbb R^{1+j+k}\setminus\{0\}$ such that (\ref{n70}) holds for all $(\sigma,\lambda,z)\in \mathcal U\times (0,+\infty)\times Z$, where  $Z$ is given in Lemma \ref{n66}.
\end{lem}

{\it Proof.} By contradiction, it follows from Lemma \ref{n66} that  there exist $(\sigma_1,\lambda_1, W^1), \cdots,\\ (\sigma_{l+1},\lambda_{l+1}, W^{l+1})\in \mathcal U\times (0,+\infty)\times T_{\bar x(0)}M$ and $\delta_0>0$ such that (\ref{n60}) and (\ref{n61}) hold.
Without loss of generality, we assume that $\lambda_1=\max\{\lambda_1,\cdots,\lambda_{l+1}\}>0$. The following argument is split into three steps.

{\it Step 1.} We claim that, given any small $\alpha>0$, there exists $\epsilon_0>0$ such that, for all $\epsilon\in[0,\epsilon_0]$, there exist measurable subset $A_\eta\subseteq[0,T]$ ($\eta=1,\cdots,l+1$) with measure $\frac{\lambda_\eta}{\lambda_1}T$,  and   $ E_\epsilon,F_\epsilon\subset [0,T]$ with $|F_\epsilon|=\lambda_1\epsilon^2T$ and $|E_\epsilon|=\epsilon T$ such that the following properties hold:
\begin{description}
\item[(i)]For any $\nu=(\nu_1,\cdots,\nu_{l+1})$ satisfying
\begin{equation}\label{n80}\nu_\eta\geq0,\quad \eta=1,\cdots,l+1;\quad\sum_{\eta=1}^{l+1}\nu_\eta=1,
\end{equation}
there exist measurable subsets $E_\nu^1,\cdots, E_\nu^{l+1}$ of $[0,T]$ and a positive constant $C$ such that
\begin{align}\label{ns1900}
&\cup_{\eta=1}^{l+1}E_\nu^\eta=[0,T],
\\ \label{n191}&E_\nu^i\cap E_\nu^\eta=\emptyset,\quad\textrm{if}\;i\neq \eta,
\\\label{n192}& |E_\nu^\eta|=\nu_\eta T,\quad \eta=1,\cdots,l+1,
\\\label{n193}& \sum_{\eta=1}^{l+1}H(E_\nu^\eta,E_{\hat\nu}^\eta)\leq C\sum_{\eta=1}^{l+1}|\nu_\eta-\hat\nu_\eta|,
\end{align}
where $\hat\nu=(\hat\nu_1,\cdots,\hat\nu_{l+1})$ satisfies (\ref{n80}).
\item[(ii)] Set
\begin{align}\label{n200}
u_\nu^\epsilon(t)=I_{(F_\epsilon\cup E_\epsilon)^c}(t)\bar u(t)+I_{E_\epsilon\setminus F_\epsilon}(t)u(t)+\sum_{\eta=1}^{l+1}I_{F_\epsilon\cap E_\nu^\eta\cap A_\eta}(t)\sigma_\eta(t),\quad a.e.\,t\in[0,T],
\end{align}
where $A^c$ is the complement of set $A$, and $I_A(\cdot)$ is the indicator funtion of set $A$.
Denote by $x(\cdot;u_\nu^\epsilon)$ the solution to (\ref{25})  corresponding to the initial state $exp_{\bar x(0)}(\epsilon V+\epsilon^2\sum_{\eta=1}^{l+1}\nu_\eta W^\eta)$ and control $u_\nu^\epsilon(\cdot)$, where $exp_{\bar x(0)}$ is the exponential map at $\bar x(0)$ (see Section \ref{a2}). It holds that
\begin{align}\label{n201}
|V_\nu^\epsilon(t)-\epsilon X_{u,V}(t)-\epsilon^2 \sum_{\eta=1}^{l+1}\nu_\eta Y_{u,\sigma_\eta}^{\lambda_\eta, W^\eta}(t)|\leq \alpha\epsilon^2,\quad\forall\;t\in[0,T],
\end{align}
where
\begin{align}\label{rv}
V_\nu^\epsilon(t)\stackrel{\triangle}{=}exp_{\bar x(t)}^{-1}x(t;u_\nu^\epsilon).
\end{align}

\end{description}

To show this,
we adopt the notation (\ref{434}) for abbreviation.
Let $\{e_1,\cdots,e_n\}\subset T_{\bar x(0)}M$ be an orthonormal basis. For $t\in(0,T]$, let $\{e_i(t)\}_{i=1}^n$ and $\{d_i(t)\}_{i=1}^n$ be given in the proof of Lemma \ref{n65}.
 Fix  $\epsilon>0$. By Lemma \ref{n90}, for $\eta=1,\cdots,l+1$, there exist measurable subset $A_\eta\subset[0,T]$ with measure $\frac{\lambda_\eta}{\lambda_1}T$ and $R_\eta=(R_\eta^1,\cdots,R_\eta^n)^\top: [0,T]\to \mathbb R^n$ with $|R_\eta(t)|\leq \epsilon^2$ for all $t\in[0,T]$, such that
\begin{align}\label{n175}\begin{array}{ll}
&\frac{\lambda_\eta}{\lambda_1}\displaystyle\int_0^t\left(\begin{array}{c}\langle e_1(s),f(s,\bar x(s),\sigma_\eta(s))-f[s]\rangle
\\ \vdots
\\ \langle e_n(s),f(s,\bar x(s),\sigma_\eta(s))-f[s]\rangle
\end{array}\right)ds
\\=&\displaystyle\int_{A_\eta\cap [0,t]}\left(\begin{array}{c}\langle e_1(s),f(s,\bar x(s),\sigma_\eta(s))-f[s]\rangle
\\ \vdots
\\ \langle e_n(s),f(s,\bar x(s),\sigma_\eta(s))-f[s]\rangle
\end{array}\right)ds+R_\eta(t),
\end{array}\end{align}
for all $t\in[0,T]$. For $t\in[0,T]$, denote by $\hat \sigma_\eta(t)=I_{A_\eta}(t)\sigma_\eta(t)+I_{A_\eta^c}(t)\bar u(t)$ if $\eta=2,\cdots,l+1$, and by $\hat\sigma_1(t)=\sigma_1(t)$.
We still obtain by Lemma \ref{n90} that, there exist measurable subset $F_\epsilon\subset [0,T]$ with measure $\lambda_1\epsilon^2 T$, $S_\eta=(S_\eta^1,\cdots,S_\eta^n)^\top: [0,T]\to \mathbb R^n$ with $\eta=1,\cdots,l+1$  and $Q=(Q^1,\cdots,Q^n)^\top: [0,T]\to\mathbb R^n$ such that
\begin{align}\label{n171}\begin{array}{ll}
&\lambda_1\epsilon^2\displaystyle\int_0^t\left(\begin{array}{c}
\langle e_1(s),f(s,\bar x(s),\hat\sigma_\eta(s))-f[s]\rangle
\\\vdots
\\ \langle e_n(s),f(s,\bar x(s),\hat\sigma_\eta(s))-f[s]\rangle
\end{array}
\right) ds
\\[2mm]=&\displaystyle\int_{F_\epsilon\cap [0,t]}\left(\begin{array}{c}
\langle e_1(s),f(s,\bar x(s),\hat\sigma_\eta(s))-f[s]\rangle
\\\vdots
\\ \langle e_n(s),f(s,\bar x(s),\hat\sigma_\eta(s))-f[s]\rangle
\end{array}
\right) ds+S_\eta(t),\end{array}
\end{align}
\begin{align}\label{n170}\begin{array}{ll}
&\lambda_1\epsilon^2\displaystyle\int_0^t\left(\begin{array}{c}
\langle e_1(s),f(s,\bar x(s),u(s))-f[s]\rangle
\\\vdots
\\ \langle e_n(s),f(s,\bar x(s),u(s))-f[s]\rangle
\end{array}
\right) ds
\\=&\displaystyle\int_{F_\epsilon\cap [0,t]}\left(\begin{array}{c}
\langle e_1(s),f(s,\bar x(s),u(s))-f[s]\rangle
\\\vdots
\\ \langle e_n(s),f(s,\bar x(s),u(s))-f[s]\rangle
\end{array}
\right) ds+Q(t),
\end{array}\end{align}
and
\begin{align*}
|S_\eta(t)|\leq \epsilon^3,\quad|Q(t)|\leq\epsilon^3,\quad\forall\, t\in[0,T],\;\eta=1,\cdots,l+1.
\end{align*}
Also, there exists a measurable subset $E_\epsilon\subset[0,T]$ with measure $\epsilon T$, $G=(G^1,\cdots,G^n)^\top:    [0,T]\to\mathbb R^n$ and $D=(D^1,\cdots,D^n)^\top: [0,T]\to\mathbb R^n$ such that
\begin{align}\label{n160}\begin{array}{ll}
&
\epsilon\displaystyle\int_0^t\left(\begin{array}{c}
\langle e_1(s),f(s,\bar x(s),u(s))-f[s]\rangle
\\\vdots
\\\langle e_n(s),f(s,\bar x(s),u(s))-f[s]\rangle
\end{array}\right) I_{F_\epsilon^c}(s)ds
\\=&\displaystyle\int_{[0,t]\cap E_\epsilon}\left(\begin{array}{c}
\langle e_1(s),f(s,\bar x(s),u(s))-f[s]\rangle
\\\vdots
\\\langle e_n(s),f(s,\bar x(s),u(s))-f[s]\rangle
\end{array}\right) I_{F_\epsilon^c}(s)ds+G(t),
\end{array}\end{align}
\begin{align}\label{n188}
\begin{array}{ll}
&
\epsilon\displaystyle\int_0^t\left(\begin{array}{c}
(\nabla_xf(s,\bar x(s),u(s))-\nabla_xf[s])(d_1(s),X_{u,V}(s))
\\\vdots
\\(\nabla_xf(s,\bar x(s),u(s))-\nabla_xf[s])(d_n(s),X_{u,V}(s))
\end{array}\right)I_{F_\epsilon^c}(s) ds
\\=&\displaystyle\int_{[0,t]\cap E_\epsilon}\left(\begin{array}{c}(\nabla_xf(s,\bar x(s),u(s))-\nabla_xf[s])(d_1(s),X_{u,V}(s))
\\\vdots
\\(\nabla_xf(s,\bar x(s),u(s))-\nabla_xf[s])(d_n(s),X_{u,V}(s))
\end{array}\right) I_{F_\epsilon^c}(s)ds+D(t),
\end{array}\end{align}
and
\begin{align*}
|G(t)|\leq \epsilon^3,\quad|D(t)|\leq \epsilon^3,\quad\forall\, t\in[0,T].
\end{align*}
For any $\nu=(\nu_1,\cdots,\nu_{l+1})$ satisfying (\ref{n80}),
it follows from Corollary \ref{n130} that, there exist measurable subsets $E_\nu^1,\cdots, E_\nu^{l+1}$ of $[0,T]$ such that (\ref{ns1900})-(\ref{n193}) hold, and
\begin{align}\label{n161}\begin{array}{ll}
&\sum_{\eta=1}^{l+1}\nu_\eta\displaystyle\int_0^t\left(\begin{array}{c}
\langle f(s,\bar x(s),\hat\sigma_\eta(s))-f[s],e_1(s)\rangle
\\ \vdots
\\ \langle f(s,\bar x(s),\hat\sigma_\eta(s))-f[s],e_n(s)\rangle
\end{array}
\right) I_{F_\epsilon}(s)ds
\\=& \sum_{\eta=1}^{l+1}\displaystyle\int_{[0,t]\cap E_\nu^\eta}\left(\begin{array}{c}
\langle f(s,\bar x(s),\hat\sigma_\eta(s))-f[s],e_1(s)\rangle
\\ \vdots
\\ \langle f(s,\bar x(s),\hat\sigma_\eta(s))-f[s],e_n(s)\rangle
\end{array}
\right) I_{F_\epsilon}(s)ds+C(t),
\end{array}\end{align}
where $C=(C^1,\cdots,C^n)^\top: [0,T]\to \mathbb R^n$ satisfies $|C(t)|\leq \epsilon^3$.

Recall $u_\nu^\epsilon(\cdot)$ defined in (\ref{n200}). Denote by $\hat x^\epsilon(\cdot)$ the solution to (\ref{25}) with initial state $\exp_{\bar x(0)}\left(\epsilon V+\epsilon^2\sum_{\eta=1}^{l+1}\nu_\eta W^\eta\right)$ and control $\bar u(\cdot)$.
Then, by \cite[Proposition 4.2]{d} and \cite[Lemma 5.2]{cdz}, there exists $\epsilon_1>0$ depending on $|V|$ and $\bar x(0)$ such that
\begin{align}\label{n155}\begin{array}{ll}
&\rho(x(t;u_\nu^\epsilon),\bar x(t))
\\\leq&\rho(x(t;u_\nu^\epsilon),\hat x^\epsilon(t))+\rho(\hat x^\epsilon(t),\bar x(t))
\\\leq& 2L(1+\rho(x_0,exp_{\bar x(0)}(\epsilon V+\epsilon^2\sum_{\eta=1}^{l+1}\nu_\eta W^\eta)))e^{Lt}(\epsilon+\lambda_1\epsilon^2)T
\\&+C_{\bar x(0)}\rho(\bar x(0),exp_{\bar x(0)}(\epsilon V+\epsilon^2\sum_{\eta=1}^{l+1}\nu_\eta W^\eta))
\\\leq & 2L\Big(1+\rho(x_0,\bar x(0))+\epsilon|V|+\epsilon^2\sum_{\eta=1}^{l+1}|W^\eta|\Big)e^{Lt}(\epsilon+\lambda_1\epsilon^2)T
\\&+C_{\bar x(0)}(\epsilon|V|+\epsilon^2\sum_{\eta=1}^{l+1}|W^\eta|),
\end{array}
\end{align}
for $t\in[0,T]$ and $\epsilon\in[0,\epsilon_1]$, where $C_{\bar x(0)}$ is the positive constant depending on $\bar x(0)$, and we have used (\ref{80}) and condition $(C2)$. Thus, for $\epsilon>0$ small enough, we can define $V^\epsilon_\nu(\cdot)$  given by (\ref{rv}).
Fix $t\in[0,T]$.
 Denote by
 \begin{equation}\label{n150}\beta^\nu_\epsilon(\theta;t)\equiv exp_{\bar x(t)}\theta  V^\epsilon_\nu(t),\;\theta\in[0,1].
 \end{equation}
 Then, $\beta^\nu_\epsilon(\cdot;t)$ is a geodesic  starting from $\bar x(t)$ and ending at $x(t;u_\nu^\epsilon)$.  Applying Lemma \ref{317} to (\ref{n150}), we obtain
 \begin{equation}\label{n157}\frac{\partial}{\partial\theta}\Big|_{\theta=0}\beta^\nu_\epsilon(\theta;t)= V_\nu^\epsilon(t).\end{equation}
  By Lemma \ref{317}, (\ref{80}), (\ref{81}),  (\ref{18}), (\ref{96}), (\ref{260}),  (\ref{n7}), (\ref{r1}), (\ref{n155}), (\ref{n157}), \cite[Lemma 2.3]{cdz},  and Newton-Leibniz formula, we have, for each $i=1,\cdots,n$,
 \begin{equation}\label{n181}\begin{array}{ll}
&\langle V^\epsilon_\nu(t),e_i(t)\rangle-\langle\epsilon V+\epsilon^2\sum_{\eta=1}^{l+1}\nu_\eta W^\eta, e_i\rangle
\\=&-\frac{1}{2}\int_0^t\frac{\partial}{\partial s} \nabla_1\rho^2(\bar x(s),x(s,u_\nu^\epsilon))(e_i(s)) ds
\\=&\int_0^t\left\{-\frac{1}{2} \nabla_{2} \nabla_{1} \rho^{2}\left(\bar x(s), x(s,u_\nu^\epsilon)\right)\left(e_i(s), f\left(s, x(s,u_\nu^\epsilon), u_\nu^\epsilon(s)\right)\right)\right.
\\&+\frac{1}{2} \nabla_{2} \nabla_{1} \rho^{2}(\bar x(s), \bar x(s))(e_i(s), f(s, \bar x(s), u_\nu^\epsilon(s)))
\\&-\frac{1}{2}\nabla_1^2\rho^2(\bar x(s),x(s,u_\nu^\epsilon))(e_i(s),f[s])+\frac{1}{2}\nabla_1^2\rho^2(\bar x(s),\bar x(s))(e_i(s),f[s])
\\&\left.+\langle e_i(s),f(s,\bar x(s),u_\nu^\epsilon(s))-f[s]\rangle\right\}ds
\\ =&-\frac{1}{2}\int_0^t\int_0^1\frac{\partial}{\partial\theta}\nabla_2\nabla_1\rho^2(\bar x(s),\beta_\epsilon^\nu(\theta;s))(e_i(s),f(s,\beta_\epsilon^\nu(\theta;s),u_\nu^\epsilon(s)))d\theta ds
\\&-\frac{1}{2}\int_0^t\int_0^1\frac{\partial}{\partial\theta}\nabla_1^2\rho^2(\bar x(s),\beta_\epsilon^\nu(\theta;s))(e_i(s),f[s])d\theta ds
\\&+\int_0^t \langle e_i(s),f(s,\bar x(s),u_\nu^\epsilon(s))-f[s]\rangle ds
 \\=&\int_0^t\left\{\nabla_xf[s](d_i(s),V^\epsilon_\nu(s))+\langle e_i(s),f(s,\bar x(s),u_\nu^\epsilon(s))-f[s]\rangle\right.
 \\&+\nabla_xf(s,\bar x(s),u_\nu^\epsilon(s))(d_i(s),V^\epsilon_\nu(s))-\nabla_xf[s](d_i(s),V^\epsilon_\nu(s))
 \\&-\frac{1}{2}R(e_i(s),V^\epsilon_\nu(s),f[s],V^\epsilon_\nu(s))
 +\frac{1}{2}\nabla_x^2f[s](d_i(s),V^\epsilon_\nu(s),V^\epsilon_\nu(s))\}ds
 \\&+A_i^\epsilon(t)+B_i^\epsilon(t),
 \end{array}
 \end{equation}
where
\begin{align*}
A_i^\epsilon(t)=&-\frac{1}{2}\int_0^t\int_0^1\Big[\nabla_2\nabla_1
\rho^2(\bar x(s),\beta_\epsilon^\nu(\theta;s))[e_i(s),\nabla_{\frac{\partial}{\partial\theta}\beta_\epsilon^\nu(\theta;s)}f(s,\cdot,u_\nu^\epsilon(s))]
\\&-\nabla_2\nabla_1\rho^2(\bar x(s),\bar x(s))[e_i(s),\nabla_{V_\nu^\epsilon(s)}f(s,\cdot,u_\nu^\epsilon(s))]\Big]d \theta ds
\\&-\frac{1}{2}\int_0^t\nabla_x^2f[s][d_i(s),V_\nu^\epsilon(s),V_\nu^\epsilon(s)]ds,
\\ B_i^\epsilon(t)=&-\frac{1}{2}\int_0^t\int_0^1\Big[\nabla_2^2\nabla_1
\rho^2(\bar x(s),\beta_\epsilon^\nu(\theta;s))[e_i(s),f(s,\beta_\epsilon^\nu(\theta;s),\bar u(s)),\frac{\partial}{\partial\theta}\beta_\epsilon^\nu(\theta;s)]
\\&-\nabla_2^2\nabla_1
\rho^2(\bar x(s),\bar x(s))[e_i(s),f(s,\bar x(s),\bar u(s)),V_\nu^\epsilon(s)]
\\&+\nabla_2\nabla_1^2
\rho^2(\bar x(s),\beta_\epsilon^\nu(\theta;s))[e_i(s),f[s],\frac{\partial}{\partial\theta}\beta_\epsilon^\nu(\theta;s)]
\\&-\nabla_2\nabla_1^2
\rho^2(\bar x(s),\bar x(s))[e_i(s),f[s],V_\nu^\epsilon(s)]
-R[e_i(s),V_\nu^\epsilon(s),f[s],,V_\nu^\epsilon(s)]
\\&\nabla_2^2\nabla_1\rho^2(\bar x(s),\beta_\epsilon^\nu(\theta;s))[e_i(s),f(s,\beta_\epsilon^\nu(\theta;s),u_\nu^\epsilon(s))-f(s,\beta_\epsilon^\nu(\theta;s),\bar u(s)),
\\&\frac{\partial}{\partial \theta}\beta_\epsilon^\nu(\theta;s)]\Big]d\theta ds.
\end{align*}
Set
\begin{align}\label{gcs10}
V_\nu^\epsilon(s)=\sum_{k=1}^n a_k^\epsilon (s)e_k(s),\quad\forall\,s\in[0,T].
\end{align}
Since $\beta_\epsilon^\nu(\cdot;s)$ is a geodesic, we have
\begin{align}\label{gcs}
\frac{\partial}{\partial \tau}\beta_\epsilon^\nu(\tau;s)=L_{\bar x(s)\beta_\epsilon^\nu(\tau;s)}V_\nu^\epsilon(s)=\sum_{k=1}^na_k^\epsilon(s) L_{\bar x(s)\beta_\epsilon^\nu(\tau;s)} e_k(s),\quad s\in[0,T],\;\tau\in[0,1],
\end{align}
and
\begin{align}\label{gcs1}
|\frac{\partial}{\partial \tau}\beta_\epsilon^\nu(\tau;s)|^2=|V_\nu^\epsilon(s)|^2=\sum_{k=1}^n a_k^\epsilon(s)^2.
\end{align}
We obtain from \cite[(2.17)]{cdz} and (\ref{n155}) that
\begin{align}\label{gcs7}
\sup_{s\in[0,T]}\sum_{k=1}^n a_k^\epsilon(s)^2=O(\epsilon^2),
\end{align}
where $O(\alpha)$ is a number satisfying $|\lim_{\alpha\to 0^+}O(\alpha)|<\infty$.

By \cite[(2.20)]{cdz}, (\ref{gcs}) and Newton-Leibniz formula, we have
\begin{align}\label{gcs5}
\begin{array}{ll}&A_i^\epsilon(t)
\\=&-\frac{1}{2}\int_0^t\int_0^1\int_0^\theta\Big[\sum_{k,l=1}^n a_k^\epsilon(s)a_l^\epsilon(s)\nabla_2^2\nabla_1\rho^2(\bar x(s),\beta_\epsilon^\nu(\tau;s))(e_i(s),
\\&\nabla_{L_{\bar x(s)\beta_\epsilon^\nu(\tau;s)} e_k(s)}f(s,\cdot,u_\nu^\epsilon(s)),
L_{\bar x(s)\beta_\epsilon^\nu(\tau;s)} e_l(s))
\\&+\nabla_2\nabla_1\rho^2(\bar x(s),\beta_\epsilon^\nu(\tau;s))(e_i(s),
\\&\nabla_{\frac{\partial}{\partial\tau}\beta_\epsilon^\nu(\tau;s)}\nabla_{\frac{\partial}{\partial\tau}\beta_\epsilon^\nu(\tau;s)}f(s,\cdot,u_\nu^\epsilon(s))\Big]d\tau d\theta ds
\\&-\frac{1}{2}\int_0^t\nabla_x^2f[s](d_i(s),V_\nu^\epsilon(s),V_\nu^\epsilon(s))ds
\\=& -\frac{1}{2}\int_0^t\int_0^1\int_0^\theta\sum_{k,l=1}^n a_k^\epsilon(s)a_l^\epsilon(s)\nabla_2^2\nabla_1\rho^2(\bar x(s),\beta_\epsilon^\nu(\tau;s))\Big(e_i(s),
\\&\nabla_{L_{\bar x(s)\beta_\epsilon^\nu(\tau;s)} e_k(s)}f(s,\cdot,u_\nu^\epsilon(s)),
L_{\bar x(s)\beta_\epsilon^\nu(\tau;s)} e_l(s)\Big)d\tau d\theta ds
\\&+\frac{1}{2}\int_0^t\Big[\nabla_x^2f(s,\bar x(s),u_\nu^\epsilon(s))-\nabla_x^2f[s]\Big](d_i(s),V_\nu^\epsilon(s),V_\nu^\epsilon(s))ds
\\&+\int_0^t\int_0^1\sum_{k,l=1}^n\Big[\nabla_x^2f(s,\beta_\epsilon^\nu(\tau;s),u_\nu^\epsilon(s))\Big(\widetilde{d\exp_{\beta_\epsilon^\nu(\tau;s)}^{-1}}\Big|_{\bar x(s)}e_i(s),
\\&L_{\bar x(s)\beta_\epsilon^\nu(\tau;s)} e_k(s),L_{\bar x(s)\beta_\epsilon^\nu(\tau;s)} e_l(s)\Big)
\\&-\nabla_x^2f(s,\bar x(s),u_\nu^\epsilon(s))\Big(d_i(s),e_k(s),e_l(s)\Big)\Big](1-\tau)a_k^\epsilon(s)a_l^\epsilon(s)d\tau ds,
\end{array}\end{align}
where $\widetilde{d\exp_{\beta_\epsilon^\nu(\tau;s)}^{-1}}\Big|_{\bar x(s)}e_i(s)$ is the dual covector of $d\exp_{\beta_\epsilon^\nu(\tau;s)}^{-1}\Big|_{\bar x(s)}e_i(s)$.
Similarly we have
\begin{align}\label{gcs6}
\begin{array}{ll}&B_i^\epsilon(t)
\\=&-\frac{1}{2}\int_0^t\int_0^1\int_0^\theta\sum_{k,l=1}^n\Big[\nabla_2^3\nabla_1\rho^2(
\bar x(s),\beta_\epsilon^\nu(\tau;s))\Big(e_i(s),f(s,\beta_\epsilon^\nu(\tau;s),\bar u(s)),
\\&L_{\bar x(s)\beta_\epsilon^\nu(\tau;s)} e_k(s), L_{\bar x(s)\beta_\epsilon^\nu(\tau;s)} e_l(s)\Big)+\nabla_2^2\nabla_1^2\rho^2(
\bar x(s),\beta_\epsilon^\nu(\tau;s))\Big(e_i(s),f[s],
\\& L_{\bar x(s)\beta_\epsilon^\nu(\tau;s)} e_k(s), L_{\bar x(s)\beta_\epsilon^\nu(\tau;s)} e_l(s)\Big)-2R(e_i(s),e_k(s),f[s],e_l(s
))
\\&+\nabla_2^2\nabla_1\rho^2(\bar x(s),\beta_\epsilon^\nu(\tau;s))
\Big( e_i(s),\nabla_{L_{\bar x(s)\beta_\epsilon^\nu(\tau;s)}e_k(s)}f(s,\cdot,\bar u(s)),L_{\bar x(s)\beta_\epsilon^\nu(\tau;s)}e_l(s)\Big)\Big]
\\&\cdot a_k^\epsilon(s)a_l^\epsilon(s)d\tau d\theta ds+o(\epsilon^2),
\end{array}\end{align}
where we have used  (\ref{gcs1}), (\ref{gcs7}) and (\ref{472}), and $o(\alpha)$ is a tensor of proper type and satisfies  $\lim_{\alpha\to0^+}\frac{o(\alpha)}{\alpha}=0$.

  We obtain from (\ref{n160}), (\ref{n161}), (\ref{n171}) and (\ref{n175}), that
 \begin{align}\label{n183}\begin{array}{ll}
&\int_0^t\langle e_i(s),f(s,\bar x(s),u_\nu^\epsilon(s))-f[s]\rangle ds
\\=&\epsilon\int_0^t\langle e_i(s),f(s,\bar x(s),u(s))-f[s]\rangle ds
\\&+\epsilon^2\sum_{\eta=1}^{l+1}\nu_\eta\lambda_\eta\int_0^t\langle e_i(s),f(s,\bar x(s),\sigma_\eta(s))-f[s]\rangle ds
\\&-\Big(G^i(t)+\sum_{\eta=1}^{l+1}\nu_\eta (S_\eta^i(t)+R_\eta^i(t)\lambda_1\epsilon^2)+C_i(t)\Big)+o(\epsilon^2),\quad \forall\;t\in[0,T].
\end{array}\end{align}
Recall that $X_{u,V}(\cdot)$ is the solution to (\ref{n14}). We obtain that
\begin{align}\label{n180}\begin{array}{ll}
&\epsilon\langle X_{u,V}(t), e_i(t)\rangle-\epsilon\langle V, e_i\rangle
\\=&\epsilon\int_0^t\Big(\nabla_xf[s](d_i(s),X_{u,V}(s))+\langle e_i(s),f(s,\bar x(s),u(s))-f[s]\rangle\Big)ds.
\end{array}\end{align}

By subtracting (\ref{n180}) from (\ref{n181}), we obtain from  \cite[Lemma 4.1]{cdz}, ($C2$), ($C3$), (\ref{n183}), (\ref{gcs5}), (\ref{gcs6}), and (\ref{gcs7}) that
$$
\langle V_\nu^\epsilon(t)-\epsilon X_{u,V}(t),e_i(t)\rangle=\int_0^t\nabla_xf[s](d_i(s),V_\nu^\epsilon(s)-\epsilon X_{u,V}(s))ds+o(\epsilon).
$$
 Applying the Gronwall's inequality to the above inequality, we obtain
$$
|\langle V_\nu^\epsilon(t)-\epsilon X_{u,V}(t),e_i(t)\rangle|\leq C \cdot o(\epsilon),\quad \forall \,t\in[0,T],
$$
for some positive constant $C$. Consequently, we have
\begin{equation}\label{n187}
V_\nu^\epsilon(t)-\epsilon X_{u,V}(t)=o(\epsilon),\quad\forall \,t\in[0,T].
\end{equation}
It follows from  (\ref{n187}), (\ref{n188}) and ($C2$) that
\begin{align}\label{n189}\begin{array}{ll}
&\int_0^t\Big(\nabla_xf(s,\bar x(s),u_\nu^\epsilon(s))(d_i(s),V_\nu^\epsilon(s))-\nabla_xf[s](d_i(s),V_\nu^\epsilon(s))\Big)ds
\\=&\epsilon^2\int_0^t(\nabla_xf(s,\bar x(s),u(s))-\nabla_xf[s])(d_i(s),X_{u,V}(s))ds+o(\epsilon^2).
\end{array}\end{align}

Denote by $Y^\nu(\cdot)$ the solution to  the following equation
\begin{equation}\label{n1900}
\left\{\begin{array}{lll}
\nabla_{\dot{\bar x}(t)}Y^\nu&=&\nabla_xf[t](\cdot,Y^\nu(t))+\sum_{\eta=1}^{l+1}\lambda_\eta \nu_\eta (f(t,\bar x(t),\sigma_\eta(t))-f[t])
\\&&+\frac{1}{2}\nabla_x^2f[t](\cdot,X_{u,V}(t),X_{u,V}(t))+\nabla_xf(t,\bar x(t),u(t))(\cdot,X_{u,V}(t))
\\&&-\nabla_xf[t](\cdot,X_{u,V}(t))-\frac{1}{2}R(\cdot,X_{u,V}(t),f[t],X_{u,V}(t)),\;a.e. \,t\in(0,T),
\\ Y^\nu (0)&=&\sum_{\eta=1}^{l+1}\nu_\eta W^\eta.
\end{array}
\right.
\end{equation}
It is easy to check that
\begin{align*}
\sum_{\eta=1}^{l+1}\nu_\eta Y_{u,\sigma_\eta}^{\lambda_\eta, W^\eta}(t)=Y^\nu(t),\quad\forall\;t\in[0,T],
\end{align*}
where $ Y_{u,\sigma_\eta}^{\lambda_\eta, W^\eta}(t)$ is defined by (\ref{n20}).

By employing (\ref{n181}), (\ref{n183}), (\ref{n180}), (\ref{n187}), (\ref{n189}), (\ref{n190}), we  derive
\begin{align*}
&\langle V_\nu^\epsilon(t)-\epsilon X_{u,V}(t)-\epsilon^2 Y^\nu(t),e_i(t)\rangle
\\=&\int_0^t\nabla_xf[s](d_i(s),V_\nu^\epsilon(s)-\epsilon X_{u,V}(s)-\epsilon^2 Y^\nu(s))ds+A_i^\epsilon(t)+B_i^\epsilon(t)+o(\epsilon^2),\quad\forall \;t\in[0,T].
\end{align*}
Applying the Gronwall's inequality and \cite[Lemma 4.1]{cdz} to the above identity, we obtain that
\begin{align*}
|V_\nu^\epsilon(t)-\epsilon X_{u,V}(t)-\epsilon^2Y^\nu(t)|\leq  2e^{LT}\Big[\sum_{i=1}^n\max_{s\in[0,T]}(|A_i^\epsilon(s)|+|B_i^\epsilon(s)|)+o(\epsilon^2 )\Big],\;\forall\,t\in[0,T].
\end{align*}
Fix any small $\alpha>0$. By Lebesgue's dominated convergence theorem, \cite[(2.24)]{cdz}, \cite[Lemma 2.3]{cdz} (\ref{gcs7}),   (\ref{n155}),  ($C2$) and ($C3$), there exists $\epsilon_0\in(0,\epsilon_1]$ such that, for all $\epsilon\in[0,\epsilon_0]$,
\begin{align*}
2e^{LT}\Big[\sum_{i=1}^n\max_{s\in[0,T]}(|A_i^\epsilon(s)|+|B_i^\epsilon(s)|)+o(\epsilon^2)\Big]\leq \alpha \epsilon^2,
\end{align*}
and (\ref{n201}) follows.

 \medskip

{\it Step 2.} Given $\alpha>0$, we claim that, there exists $\hat\epsilon_0\in(0,\epsilon_0]$ such that, for all $\epsilon\in[0,\hat\epsilon_0]$ and any $ \nu=(\nu_1,\cdots,\nu_{l+1})$ satisfying (\ref{n80}), it holds for $\varphi=\phi_i, \psi$ with $i=0,1,\cdots,j$ that
\begin{equation}\label{n69}\begin{array}{ll}
&\Big|\varphi(x(0;u_\nu^\epsilon),x(T;u_\nu^\epsilon))-\varphi(\bar x(0),\bar x(T))
\\&-\epsilon\Big(\nabla_1\varphi(\bar x(0),\bar x(T))(V)+\nabla_2\varphi(\bar x(0),\bar x(T))(X_{u,V}(T))\Big)
\\ &-\epsilon^2\sum_{\eta=1}^{l+1}\nu_\eta\Big(\nabla_1\varphi(\bar x(0),\bar x(T))(W^\eta)+\nabla_2\varphi(\bar x(0),\bar x(T))(Y_{u,\sigma_\eta}^{\lambda_\eta, W^\eta}(T))
\\&+\frac{1}{2} D^2\varphi(u,V)\Big)\Big|
\\\leq& \alpha \epsilon^2,
\end{array}\end{equation}
where
$D^2\varphi$ is defined in (\ref{n55}).

In fact, recalling (\ref{n150}), (\ref{gcs}) and (\ref{n157}), and following the same argument as that in (\ref{n181}) we obtain that
\begin{align}\label{gcs70}
\begin{array}{ll}
&\varphi(x(0;u_\nu^\epsilon),x(T;u_\nu^\epsilon))-\varphi(\bar x(0),\bar x(T))
\\=&\nabla_1\varphi(\bar x(0),\bar x(T))(\epsilon V+\epsilon^2\sum_{\eta=1}^{l+1}\nu_\eta W^\eta)+\nabla_2\varphi(\bar x(0),\bar x(T))(V_\nu^\epsilon(T))
\\&+\frac{1}{2}\nabla_1^2\varphi(\bar x(0),\bar x(T))(\epsilon V+\epsilon^2\sum_{\eta=1}^{l+1}\nu_\eta W^\eta,\epsilon V+\epsilon^2\sum_{\eta=1}^{l+1}\nu_\eta W^\eta)
\\&+\nabla_2\nabla_1\varphi(\bar x(0),\bar x(T))(\epsilon V+\epsilon^2\sum_{\eta=1}^{l+1}\nu_\eta W^\eta,V_\nu^\epsilon(T))
\\&+\frac{1}{2}\nabla_2^2\varphi(\bar x(0),\bar x(T))(V_\nu^\epsilon(T),V_\nu^\epsilon(T))+C^\epsilon,
\end{array}
\end{align}
where
\begin{align*}
&C^\epsilon
\\=&\int_0^1 \sum_{\eta,\xi=1}^n \Big\{ \Big[\nabla_1^2\varphi(\beta_\epsilon^\nu(\tau;0),\beta_\epsilon^\nu(\tau;T))(L_{\bar x(0)\beta_\epsilon^\nu(\tau;0)}e_\xi,L_{\bar x(0)\beta_\epsilon^\nu(\tau;0)}e_\eta)
\\&-\nabla_1^2\varphi(\bar x(0),\bar x(T))(e_\xi,e_\eta)\Big]
a_\xi^\epsilon(0)a_\eta^\epsilon(0)
\\&+2 a_\eta^\epsilon(0) a_\xi^\epsilon(T)\Big[\nabla_2\nabla_1\varphi(
\beta_\epsilon^\nu(\tau;0),\beta_\epsilon^\nu(\tau;T))(L_{\bar x(0)\beta_\epsilon^\nu(\tau;0)}e_\eta,
\\&L_{\bar x(T)\beta_\epsilon^\nu(\tau;T)}e_\xi(T))-\nabla_2\nabla_1\varphi(\bar x(0),\bar x(T))(e_\eta,e_\xi(T))\Big]+ a_\eta^\epsilon(T) a_\xi^\epsilon(T)\\&\Big[\nabla_2^2\varphi(
\beta_\epsilon^\nu(\tau;0),\beta_\epsilon^\nu(\tau;T))(L_{\bar x(T)\beta_\epsilon^\nu(\tau;T)}e_\eta(T),L_{\bar x(T)\beta_\epsilon^\nu(\tau;T)}e_\xi(T))
\\&-\nabla_2^2\varphi(\bar x(0),\bar x(T))(e_\eta(T),e_\xi(T))\Big]\Big\}(1-\tau)d\tau,
\end{align*}
with $a_\eta^\epsilon(\cdot)$  defined by (\ref{gcs10}). Applying Lesbegue's dominated convergence theorem, (\ref{n155}) and
(\ref{gcs7}) to (\ref{gcs70}), we obtain via (\ref{n201}) that, there exists $\hat\epsilon_0\in(0,\epsilon_0]$ such that (\ref{n69}) holds for all $\epsilon\in[0,\hat\epsilon_0]$.

 \medskip
{ \it Step 3.}  According to (\ref{n69}), (\ref{n60}), (\ref{n53}) and (\ref{n61}), there exists $\tilde\epsilon_0>0$ such that, for $\epsilon\in[0,\tilde\epsilon_0]$ and $\nu$ satisfying
 (\ref{n80}),  the following relations hold:
 \begin{align*}
&\Big| \epsilon^{-2}\psi(x(0;u_\nu^\epsilon),x(T;u_\nu^\epsilon))-\sum_{\eta=1}^{l+1}\nu_\eta\Big(\nabla_1\psi(\bar x(0),\bar x(T))(W^\eta)
\\&+\nabla_2\psi(\bar x(0),\bar x(T))(Y_{u,\sigma_\eta}^{\lambda_\eta, W^\eta}(T))
+\frac{1}{2} D^2\psi(u,V)\Big)\Big|\leq\delta_0;
\end{align*}
 If $i\in I_N$,
 \begin{align}\label{n210}
\phi_i(x(0;u_\nu^\epsilon),x(T;u_\nu^\epsilon))
=\phi_i(\bar x(0),\bar x(T))+O(\epsilon)
<0;
\end{align}
If $i\in I_{AO}\setminus I_0^{\prime\prime}$,
 \begin{align}\label{n211}
&\phi_i(x(0;u_\nu^\epsilon),x(T;u_\nu^\epsilon))\nonumber
\\=&\epsilon(\nabla_1\phi_i(\bar x(0),\bar x(T))(V)+\nabla_2\phi_i(\bar x(0),\bar x(T))(X_{u,V}(T)))+o(\epsilon)<0;
\end{align}
If $i\in I_0^{\prime\prime}$,
\begin{align}\label{n212}
\begin{array}{ll}&\phi_i(x(0;u_\nu^\epsilon),x(T;u_\nu^\epsilon))
\\=&\epsilon^2\Big(\nabla_1\phi_i(\bar x(0),\bar x(T))(\sum_{\eta=1}^{l+1}\nu_\eta W^\eta)+\nabla_2\phi_i(\bar x(0),\bar x(T))\sum_{\eta=1}^{l+1}\nu_\eta Y_{u,\sigma_\eta}^{\lambda_\eta, W^\eta}(T)
\\&+\frac{1}{2}D^2\phi_i(u,V)\Big)+\Big[\phi_i(x(0;u_\nu^\epsilon),x(T;u_\nu^\epsilon))
\\&-\epsilon^2\Big(\nabla_1\phi_i(\bar x(0),\bar x(T))(\sum_{\eta=1}^{l+1}\nu_\eta W^\eta)+\nabla_2\phi_i(\bar x(0),\bar x(T))\sum_{\eta=1}^{l+1}\nu_\eta Y_{u,\sigma_\eta}^{\lambda_\eta, W^\eta}(T)
\\&+\frac{1}{2}D^2\phi_i(u,V)\Big)\Big]
\\<&0.
\end{array}\end{align}

Then, recalling (\ref{n60}), we can define a map
\begin{align*}
F:&co\{\nabla_1\psi(\bar x(0),\bar x(T))(W^\eta)+\nabla_2\psi(\bar x(0),\bar x(T))(Y_{u,\sigma_\eta}^{\lambda_\eta, W^\eta}(T))+\frac{1}{2}D^2\psi(u,V)\}_{\eta=1}^{l+1}\\&\to co\{\nabla_1\psi(\bar x(0),\bar x(T))(W^\eta)+\nabla_2\psi(\bar x(0),\bar x(T))(Y_{u,\sigma_\eta}^{\lambda_\eta, W^\eta}(T))+\frac{1}{2}D^2\psi(u,V)\}_{\eta=1}^{l+1}
\end{align*}
as follows:
\begin{align*}
&F\Big(\sum_{\eta=1}^{l+1}\nu_\eta\Big(\nabla_1\psi(\bar x(0),\bar x(T))(W^\eta)+\nabla_2\psi(\bar x(0),\bar x(T))(Y_{u,\sigma_\eta}^{\lambda_\eta, W^\eta}(T))+\frac{1}{2}D^2\psi(u,V)\Big)\Big)
\\=& -\tilde\epsilon_0^{-2}\psi(\bar x(0;u_\nu^{\tilde\epsilon_0}),x(T;u_\nu^{\tilde\epsilon_0}))+\sum_{\eta=1}^{l+1}\nu_\eta\Big(\nabla_1\psi(\bar x(0),\bar x(T))(W^\eta)
\\&+\nabla_2\psi(\bar x(0),\bar x(T))(Y_{u,\sigma_\eta}^{\lambda_\eta, W^\eta}(T))+\frac{1}{2}D^2\psi(u,V)\Big),
\end{align*}
where $u_\nu^{\tilde\epsilon_0}(\cdot)$ is given in Step 1. We obtain from \cite[Lemma 5.1]{cdz} and (\ref{n193}) that $F$ is continuous.
By Browner's fixed point theorem, there exists $\nu^0=(\nu_1^0,\cdots,\nu_{l+1}^0)$ satisfying (\ref{n80}) such that
\begin{align*}
F\Big(\sum_{\eta=1}^{l+1}\nu_\eta^0\Big(\nabla_1\psi(\bar x(0),\bar x(T))(W^\eta)+\nabla_2\psi(\bar x(0),\bar x(T))(Y_{u,\sigma_\eta}^{\lambda_\eta, W^\eta}(T))+\frac{1}{2}D^2\psi(u,V)\Big)\Big)
\\=\sum_{\eta=1}^{l+1}\nu_\eta^0\Big(\nabla_1\psi(\bar x(0),\bar x(T))(W^\eta)+\nabla_2\psi(\bar x(0),\bar x(T))(Y_{u,\sigma_\eta}^{\lambda_\eta, W^\eta}(T))+\frac{1}{2}D^2\psi(u,V)\Big),
\end{align*}
 which implies that
\begin{align}\label{n213}
\psi(x(0,u_{\nu^0}^
{\tilde\epsilon_0}),x(T,u_{\nu^0}^
{\tilde\epsilon_0}))=0.
\end{align}
Inequalities (\ref{n210})-(\ref{n212}) imply that $\phi_i(x(0,u_{\nu^0}^
{\tilde\epsilon_0}),x(T,u_{\nu^0}^
{\tilde\epsilon_0}))<0$ for all $i=0,1,\cdots,j$, which together with (\ref{n213}), contradicts the optimality of $(\bar x(\cdot),\bar u(\cdot))$.
The proof is concluded. $\Box$

\medskip
{\it Proof of Theorem \ref{n295}.}  \quad Without loss of generality, we assume $\phi_0(\bar x(0),\bar x(T))=0$. First, we shall prove the  case $k>0$.
It follows from Lemma \ref{n190} that, there eixsts
\\$\hat\ell=(\hat\ell_0,\cdots,\hat\ell_j,\hat\ell_\psi)\in\mathbb R^{1+j+k}\setminus\{0\}$ such that
(\ref{n230}) and (\ref{n231}) hold,
and the following inequality
\begin{align}\label{n232}\begin{array}{ll}
&\sum_{i=0}^j\hat\ell_i\Big(\nabla_1\phi_i(\bar x(0),\bar x(T))(W)+\nabla_2\phi_i(\bar x(0),\bar x(T))Y_{u,\sigma}^{\lambda, W}(T)
\\&+\frac{1}{2}D^2\phi_i(\bar x(0),\bar x(T))(u,V)\Big)
+\hat\ell_\psi^\top
\Big(\nabla_1\psi(\bar x(0),\bar x(T))(W)
\\&+\nabla_2\psi(\bar x(0),\bar x(T))Y_{u,\sigma}^{\lambda, W}(T)+\frac{1}{2}D^2\psi(\bar x(0),\bar x(T))(u,V)\Big)\leq 0\end{array}\end{align}
 holds for all $(\sigma(\cdot),
\lambda, W)\in\mathcal U\times (0,+\infty)\times T_{\bar x(0)}M$.
Recall (\ref{n20}), (\ref{n35}) and (\ref{n241}). Applying Newton-Leibniz formula to   (\ref{n232}), we obtain
\begin{align}\label{n240}
&(d_1\mathcal L(\bar x(0),\bar x(T),\hat\ell)+p^{\hat\ell}(0))(W)+\lambda\int_0^T \Big(H(t,\bar x(t),p^{\hat\ell}(t),\sigma(t))-H\{ t\}^{\hat\ell}\Big)dt\nonumber
\\&+\frac{1}{2}\int_0^T\Big(\nabla_x^2H\{t\}^{\hat\ell}(X_{u,V}(t),X_{u,V}(t))\nonumber
+2(\nabla_xH(t,\bar x(t),p^{\hat\ell}(t),u(t))\nonumber
\\&-\nabla_xH\{t\}^{\hat\ell})(X_{u,V}(t))-R(\tilde p^{\hat\ell}(t),X_{u,V}(t),f[t],X_{u,V}(t))\Big)dt\nonumber
\\&+\frac{1}{2}\nabla_1^2\mathcal L(\bar x(0),\bar x(T),\hat\ell)(V,V)
+\nabla_2\nabla_1\mathcal L(\bar x(0),\bar x(T),\hat\ell)(V,X_{u,V}(T))\nonumber
\\&+\frac{1}{2}\nabla_2^2\mathcal L(\bar x(0),\bar x(T),\hat\ell)(X_{u,V}(T),X_{u,V}(T))\leq 0,
\end{align}
where $p^{\hat\ell}$ solves (\ref{n35}) with $\ell$ replaced by $\hat\ell$,  $\{t\}^{\hat\ell}$ is given in (\ref{n241}), and $\tilde p^{\hat\ell}(t)$ is the dual vector of $p^{\hat\ell}(t)$.
From the above relation, one can easily obtain by contradiction argument that
\begin{align*}
(\nabla_1\mathcal L(\bar x(0),\bar x(T),\hat\ell)+p^{\hat\ell}(0))(W)+\int_0^T \Big(H(t,\bar x(t),p^{\hat\ell}(t),\sigma(t))-H\{\hat t\}^{\hat\ell}\Big)dt\leq 0,
\end{align*}
for all $(W, \sigma(\cdot))\in T_{\bar x(0)}M\times \mathcal U$. If follows from Remark \ref{n290} that $\hat\ell$ is a Lagrange multiplier, and (\ref{n236}) follows.

Then, for the case $k=0$, we claim that there exists $(\ell_0,\cdots,\ell_j)^\top\in \mathbb R^{1+j}\setminus\{0\}$ such that
\begin{align}\label{gcs34}
\sum_{\eta=0}^{j}\ell_\eta\beta_\eta\leq \sum_{\eta=0}^{j}\ell_\eta z_\eta,\;\forall\,(\beta_0,\cdots,\beta_j)^\top\in\mathcal K_{u,V},\,(z_0,\cdots,z_j)^\top \in Z.
\end{align}
If it were not true,  it follows from \cite[Theorem 11.3, p.97]{rock} that $\mathcal K_{u,V}\cap Z\neq \emptyset$. Then, there exists $(\sigma_0,\lambda_0,W^0)\in \mathcal U\times (0,+\infty)\times T_{\bar x(0)}M$, $\theta_0>0$ and $(z_0,z_1,\cdots,z_j)\in(-\infty,0)^{1+j}$ such that (\ref{gcs30}) holds. Recalling the proof of  Lemma \ref{n190}, we can show similarly that, for any small $\alpha>0$, there exists $\epsilon_0>0$ such that, for all $\epsilon\in[0,\epsilon_0]$, there eixst measurable  subsets $E_\epsilon, F_\epsilon\subset[0,T]$  with $|E_\epsilon|=\epsilon T$ and $|F_\epsilon|=\lambda_0\epsilon^2 T$ such that the following relations hold:
Set by $u^\epsilon(t)= I_{
(F_\epsilon\cup E_\epsilon)^c}(t)\bar u(t)+I_{E_\epsilon\setminus F_\epsilon}(t)u(t)+I_{F_\epsilon}(t)\sigma_0(t)$ for almost all $t\in[0,T]$. Denote by $x(\cdot;u^\epsilon)$ the solution to  (\ref{25})  corresponding to the initial state $exp_{\bar x(0)}(\epsilon V+\epsilon^2 W^0)$ and control $u^\epsilon(\cdot)$. It holds that $|V^\epsilon(t)-\epsilon X_{u,V}(t)-\epsilon^2 Y_{u,\sigma_0}^{\lambda_0,W^0}(t)|\leq \alpha\epsilon^2$ for all $t\in[0,T]$, where $V^\epsilon(t)=\exp_{\bar x(t)}^{-1}x(t;u^\epsilon)$. Then, following the same argument as that in (\ref{n210})-(\ref{n212}), we obtain that  $\phi_i(x(0;u^\epsilon),x(T;u^\epsilon))-\phi_i(\bar x(0),\bar x(T))<0$ for $i=0,1,\cdots,j$, when $\epsilon>0$ is small enough, and a contradiction follows. $\Box$

\medskip
{\it A sketch of proof of Theorem \ref{f}.}
First, we shall show that the set
\begin{align*}
\hat K\stackrel{\triangle}{=}\{ \nabla_1\Phi_{I_{AO}}(\bar x(0),\bar x(T))(W)+\nabla_2\Phi_{I_{AO}}(\bar x(0),\bar x(T))(X_{\sigma,W}(T))| \sigma\in\mathcal U, W\in T_{\bar x(0)}M\}
\end{align*}
is convex.
Then, without loss of generaty, we assume $\phi_0(\bar x(0),\bar x(T))=0$,
and denote by
\begin{align*}
\hat Z\stackrel{\triangle}{=}(-\infty,0)^{j+1}-cone \phi(\bar x(0),\bar x(T)).
\end{align*}
Using the same argument as that in Lemma \ref{n190} that, there exists $\ell=(\ell_0,\ell_1,\cdots,\ell_j,\ell_\psi)\in \mathbb R^{1+j+k}\setminus\{0\}$ such that
\begin{align*}
&\sum_{i=0}^j\ell_i\Big(\nabla_1\phi_i(\bar x(0),\bar x(T))(W)+\nabla_2\phi_i(\bar x(0),\bar x(T))(X_{\sigma, W}(T))\Big)
\\&+\ell_\psi^\top\Big(\nabla_1\psi(\bar x(0),\bar x(T))(W)+\nabla_2\psi(\bar x(0),\bar x(T))(X_{\sigma, W}(T))\Big)\leq \sum_{i=0}^j\ell_i\hat z_i,
\end{align*}
for all $\sigma\in\mathcal U$, $W\in T_{\bar x(0)}M$ and $\hat z=(\hat z_0,\cdots,\hat z_j)^\top\in\hat Z$.
By integration by parts over $[0,T]$, we obtain (\ref{n260}) and (\ref{n33}) from the above relation.
$\Box$

\subsection{Proofs of Theorems \ref{n770}--\ref{n790}}

We first prove Theorem \ref{n770}.
\\

{\it Proof of Theorem \ref{n770}.}\quad
 Set  $\tau_{-1}=0$. Fix  $\epsilon\in(0, \min_{0\leq i\leq \ell}\{\tau_{i}-\tau_{i-1}\})$. Then, it follows from (\ref{n718}) and condition $(C2)$ that, there eixsts $r_\epsilon\in(0,\epsilon)$ such that
\begin{align}\label{n719}
\Big|\frac{1}{r_\epsilon}\int_{\tau_i-r_\epsilon}^{\tau_i}\mathcal A(s,u(s))ds-\mathcal A(\tau_i,u(\tau_i))\Big|<\epsilon,
\\\label{n901}\Big|\frac{1}{r_\epsilon}\int_{\tau_i-r_\epsilon}^{\tau_i}Z(s)[A(s,u(s))-A(s,\bar u(s))]Z^{-1}(s)ds\nonumber
\\-Z(\tau_i)[A(\tau_i,u(\tau_i))-A(\tau_i,\bar u(\tau_i))]Z^{-1}(\tau_i)\Big|<\epsilon,
\\\label{n720} r_\epsilon^{-1}\left|\{\tau;\,|\tau-\tau_i|\leq r_\epsilon, |\mathcal A(\tau,u(\tau))-\mathcal A(\tau_i,u(\tau_i))|>\epsilon\}\right|<\epsilon,
\\ \label{n903}
 r_\epsilon^{-1}\left|\left\{ \tau;|\tau-\tau_i|\leq r_\epsilon, |Z(\tau)[ A(\tau,u(\tau))-A(\tau,\bar u(\tau))]Z^{-1}(\tau)\right.\right.\nonumber
\\\left.\left.-Z(\tau_i)[A(\tau_i,u(\tau_i))-A(\tau_i,\bar u(\tau_i))]Z^{-1}(\tau_i)|>\epsilon \right\} \right|<\epsilon,
\end{align}
for $i=0,\cdots,\ell$.

 By \cite[Lemma]{War78}, there exists $\{\beta^\epsilon_i\}_{i=0}^\ell\subset(0,+\infty)$ such that $\beta^\epsilon_i\to \beta_i$ as $\epsilon\to 0$ for $i=0,\cdots,\ell$, and
\begin{align}\label{n700}
\nabla\Phi_\eta\sum_{i=0}^\ell \frac{\beta_i^\epsilon}{r_\epsilon}\int_{\tau_{i}-r_{\epsilon}}^{\tau_{i}} \mathcal A(s,u(s))d s=0,\;\eta=1,\cdots,j,
\\\label{n701}\nabla\Psi\sum_{i=0}^\ell \frac{\beta_i^\epsilon}{r_\epsilon}\int_{\tau_{i}-r_{\epsilon}}^{\tau_{i}} \mathcal A(s,u(s)) d s=0.
\end{align}
Set
\begin{align*}
u_i(s)=I_{[\tau_i-r_\epsilon,\tau_i]}(s) u(s)+I_{[\tau_i-r_\epsilon,\tau_i]^c}(s)\bar u(s),\quad s\in[0,T],\;i=0,\cdots,\ell.
\end{align*}
Then,  (\ref{n700}) and (\ref{n701}) are respectively equivalent to
\begin{align}\label{n705}
\frac{\sum_{p=0}^\ell\beta_p^\epsilon }{ r_{\epsilon} }\nabla \Phi_\eta \sum_{i=0}^{\ell}\hat{\beta}_{i}^{\epsilon}  \int_0^T \mathcal A(s,u_i(s)) d s=0, \;\eta=1,\cdots,j,
\end{align}
and
\begin{align}\label{n706}
 \frac{\sum_{p=0}^\ell\beta_p^\epsilon }{ r_{\epsilon}  } \nabla \Psi \sum_{i=0}^{\ell}\hat{\beta}_{i}^{\epsilon} \int_0^T\mathcal A(s,u_i(s))d s=0,
\end{align}
where $\hat{\beta}_{i}^{\epsilon}=\frac{\beta_i^\epsilon}{\sum_{p=0}^\ell\beta_p^\epsilon}$ for $i=0,\cdots,\ell$. By Lemma \ref{l6}, there exists muturally disjoint subsets $F_0^\epsilon,\cdots, F_\ell^\epsilon$ of $[0,T] $ with $|F_i^\epsilon|=\hat{\beta}_i^\epsilon T$ for $i=0,\cdots,\ell$, such that
\begin{align}\label{n713}
&\sum_{i=0}^\ell\hat\beta_i^\epsilon\int_0^t\mathcal A(s,u_i(s)) d s=\sum_{i=0}^\ell\int_{[0,t]\cap F_i^\epsilon}\mathcal A(s,u_i(s))d s+R_{-1}(t,\epsilon),
\\
&\sum_{i=0}^\ell\hat\beta_i^\epsilon\int_0^t\vec p(s)[A(s,u_i(s))-A(s,\bar u(s))]Z^{-1}(s)\int_0^s \mathcal A(\tau,u_\eta(\tau))d\tau ds\nonumber
\\\label{n715}&=\sum_{i=0}^\ell\int_{[0,t]\cap F_i^\epsilon}\vec p(s)[A(s,u_i(s))-A(s,\bar u(s))]Z^{-1}(s)\int_0^s \mathcal A(\tau,u_\eta(\tau))d\tau ds+R_\eta(t,\epsilon),
\end{align}
for all $t\in[0,T]$, where $R_\eta(T,\epsilon)=0$ and $|R_\eta(t,\epsilon)|\leq r_\epsilon^3$ for all $t\in[0,T]$, and $\eta=-1, 0,\cdots,\ell$.
 Set $v^\epsilon(\cdot)=\sum_{i=0}^\ell I_{F_i^\epsilon}(\cdot)u_i(\cdot)$. Then, (\ref{n705}) and (\ref{n706}) are respectively  reduced to
 \begin{align}\label{n752}
\nabla\Phi_\eta\int_0^T Z(s)(\vec f(s,v^\epsilon(s))-\vec f(s,\bar u(s)))ds=0,\;\eta=1,\cdots,j,
\\\label{n753} \nabla\Psi\int_0^T Z(s)(\vec f(s,v^\epsilon(t))-\vec f(s,\bar u(s)))ds=0.
\end{align}
Denote by $X_{v^\epsilon}(\cdot)$ the solution to (\ref{n14}) with $(u(\cdot),V)$ replaced by $(v^\epsilon(\cdot),0)$. Then, we can express $X_{v^\epsilon}(\cdot)$ by
 \begin{align*}
 X_{v^\epsilon}(t)=\sum_{i=1}^n X_{v^\epsilon}^i(t)e_i(t),\quad \vec X_{v^\epsilon}(t)=(X_{v^\epsilon}^1(t),\cdots,X_{v^\epsilon}^n(t))^\top,\;\forall\, t\in[0,T].
\end{align*}
If follows from (\ref{n752}) and (\ref{n753})  that
\begin{equation}\label{n754}\nabla\phi_\eta(\bar x(T))X_{v^\epsilon}(T)=0,\quad\nabla\psi(\bar x(T))(X_{v^\epsilon}(T))=0,\;\eta=1,\cdots,j.\end{equation}

By integrating by parts over $[0,T]$, we obtain
\begin{align*}
&\int_0^T [H(t,\bar x(t),p(t), v^\epsilon(t))-H(t,\bar x(t),p(t), \bar u(t))]dt
\\=&[\ell_0\nabla\phi_0(\bar x(T))+\sum_{i=1}^j\ell_i\nabla\phi_i(\bar x(T))+\ell_\psi^\top\nabla\psi(\bar x(T))](X_{v^\epsilon}(T)),
\end{align*}
which together with (\ref{n754}) and $v^\epsilon(t)\in U(t)$ a.e. $t\in[0,T]$, implies that 
$$\nabla\phi_0(\bar x(T))(X_{v^\epsilon}(T))=0.
$$ 
Therefore, $v^\epsilon(\cdot)$ is a Pontryagin's type critical direction.

Set $\vec p(t)=(p_1(t),\cdots, p_n(t))$. Recalling (\ref{n750}) and (\ref{n830}), we obtain
 \begin{align*}
\left\{\begin{array}{l}\dot{\vec p}(t)=-\vec p(t)A(t, \bar u(t)),\;a.e. \,t\in[0,T),
\\ \vec p(T)=(\sum_{i=0}^j\ell_i\nabla \Phi_i+\ell_\psi^\top\nabla \Psi) Z(T),\end{array}\right.
\end{align*}
where $\nabla\Phi_0$ is defined by (\ref{n780}) with $i=0$.
Recalling (\ref{n830}), we have
 \begin{align}\label{n725}
\vec p(t)=(\sum_{i=0}^j\ell_i\nabla \Phi_i+\ell_\psi^\top\nabla \Psi)Z(t),\quad\forall \,t\in[0,T];
\\\label{n726}
\vec X_{v^\epsilon}(t)=Z^{-1}(t)\int_0^t  Z(s)\Big(\vec f(s,v^\epsilon(s))-\vec f(s,\bar u(s))\Big)ds, \quad\forall \,t\in[0,T].
\end{align}

By Theorem \ref{n295}, (\ref{n725}) and (\ref{n726}), we derive
\begin{align}\label{n761}\begin{array}{ll}
&\int_0^T \left\{ \vec X_{v^\epsilon}(t)^\top\Big(\nabla_x^2f[t](p(t),e_\xi(t),e_\zeta(t))\Big)_{\xi,\zeta=1}^n \vec X_{v^\epsilon}(t)+2\vec p(t)\Big( A(t,v^\epsilon(t))\right.
\\&\left.-A(t,\bar u(t))\Big)\vec X_{v^\epsilon}(t)\right\}dt
-\int_0^T \vec X_{v^\epsilon}(t)^\top\left(R(\tilde p(t),e_i(t),f[t],e_\nu(t))\right)_{ i,\nu=1}^n
\\&\vec X_{v^\epsilon}(t)dt+\vec X_{v^\epsilon}(T)^\top\Big(\sum_{i=0}^j\ell_i \nabla^2\Phi_i+\sum_{\eta=1}^k\ell_\psi^\eta\nabla^2\Psi_\eta\Big)\vec X_{v^\epsilon}(T)\leq 0,
\end{array}\end{align}
where $\nabla^2\Phi_0,\cdots, \nabla^2\Phi_j$ and $\nabla^2\Psi_\eta$ $(\eta=1,\cdots, k)$ are defined in (\ref{n759}).
It follows from (\ref{n713}) and (\ref{n726}) that
\begin{align}\label{n716}
\vec X_{v^\epsilon}(t)=&Z^{-1}(t)\sum_{i=0}^\ell \hat\beta_i^\epsilon\int_0^t Z(s)\Big(\vec f(s,u_i(s))-\vec f(s,\bar u(s))\Big)ds\nonumber
\\&-Z^{-1}(t)R_{-1}(t,\epsilon),\quad\quad\quad\quad\quad\quad\quad\quad\quad\quad\quad\quad\quad\quad\forall \,t\in[0,T].
\end{align}
From the definition of $u_i(\cdot)$ we derive that
\begin{align}\label{n735}
\vec X_{v^\epsilon}(t)=&Z^{-1}(t)\sum_{0\leq\eta<i}\hat\beta_\eta^\epsilon\int_{\tau_\eta-r_\epsilon}^{\tau_\eta}\mathcal A(s,u(s))ds\nonumber
\\&+Z(t)^{-1}\hat\beta_i^\epsilon\int_{\tau_i-r_\epsilon}^t\mathcal A(s,u(s))I_{[\tau_i-r_\epsilon, \tau_i]}(s)ds-Z^{-1}(t)R_{-1}(t,\epsilon),
\end{align}
where $t\in[\tau_{i-1},\tau_i]$ with $i=0,\cdots,\ell+1$, and we set   $\hat\beta_{\ell+1}^\epsilon=0$.
Then, we obtain from (\ref{n720}) that
\begin{align}\label{n730}\begin{array}{ll}
&\frac{1}{r_\epsilon^2}\int_0^T\vec X_{v^\epsilon}(t)^\top\Big(\nabla_x^2f[t](p(t),e_\xi(t),e_\zeta(t))\Big)_{\xi,\zeta=1}^n \vec X_{v^\epsilon}(t)dt
\\[2mm]=& \frac{1}{r_\epsilon^2}\int_0^{\tau_0}\vec X_{v^\epsilon}(t)^\top\Big(\nabla_x^2f[t](p(t),e_\xi(t),e_\zeta(t))\Big)_{\xi,\zeta=1}^n \vec X_{v^\epsilon}(t)dt
\\[2mm]&+\frac{1}{r_\epsilon^2}\sum_{i=0}^\ell\int_{\tau_i}^{\tau_{i+1}}\vec X_{v^\epsilon}(t)^\top\Big(\nabla_x^2f[t](p(t),e_\xi(t),e_\zeta(t))\Big)_{\xi,\zeta=1}^n \vec X_{v^\epsilon}(t)dt
\\[2mm]=&\frac{1}{r_\epsilon^2}\sum_{i=0}^{\ell}\Big(\sum_{0\leq \eta\leq i}\hat\beta_\eta^\epsilon\int_{\tau_\eta-r_\epsilon}^{\tau_\eta}\mathcal A(s,u(s))ds\Big)^\top\int_{\tau_i}^{\tau_{i+1}} \Big(Z^{-1}(t)\Big)^\top
\\[2mm]&\Big(\nabla_x^2f[t](p(t),e_\xi(t),e_\zeta(t))\Big)_{\xi,\zeta=1}^n Z(t)^{-1}dt \sum_{0\leq \hat\eta\leq i}\hat\beta_{\hat\eta}^\epsilon\int_{\tau_{\hat\eta}-r_\epsilon}^{\tau_{\hat\eta}}\mathcal A(s,u(s))ds+o(1)
\\[2mm]=&\frac{1}{(\sum_{p=0}^\ell\beta_p)^2}\sum_{i=0}^{\ell}\sum_{0\leq \eta,\hat\eta\leq i}\beta_\eta \beta_{\hat\eta} \mathcal A(\tau_\eta,u(\tau_\eta))^\top \int_{\tau_i}^{\tau_{i+1}} \Big(Z^{-1}(t)\Big)^\top
\\[2mm]&\Big(\nabla_x^2f[t](p(t),e_\xi(t),e_\zeta(t))\Big)_{\xi,\zeta=1}^n Z(t)^{-1}dt
  \mathcal A(\tau_{\hat\eta},u(\tau_{\hat\eta}))+o(1),
\end{array}\end{align}
where the term $o(1)$ satisfies $\lim_{\epsilon\to 0^+}o(1)=0$.
Similarly we have
\begin{align}\label{n733}\begin{array}{ll}
&\frac{1}{r_\epsilon^2}\int_0^T \vec X_{v^\epsilon}(t)^\top[R(\tilde p(t),e_i(t),f[t],e_\nu(t))]_{i,\nu=1}^n\vec X_{v^\epsilon}(t)dt
\\=&\frac{1}{(\sum_{p=0}^\ell\beta_p)^2}\sum_{i=0}^{\ell}\sum_{0\leq \eta,\hat\eta\leq i}\beta_\eta \beta_{\hat\eta} \mathcal A(\tau_\eta,u(\tau_\eta))^\top \int_{\tau_i}^{\tau_{i+1}} \Big(Z^{-1}(t)\Big)^\top
\\&[R(\tilde p(t),e_\xi(t),f[t],e_\nu(t))]_{ \xi,\nu=1}^nZ(t)^{-1}dt
  \mathcal A(\tau_{\hat\eta},u(\tau_{\hat\eta}))+o(1).
\end{array}\end{align}

Recalling (\ref{n719}), (\ref{n901}), (\ref{n715}), (\ref{n725}) and  (\ref{n716}), we obtain
\begin{align}\begin{array}{ll}
&\frac{2}{r_\epsilon^2}\int_0^T \vec p(t)[ A(t,v^\epsilon(t))-A(t,\bar u(t))]\vec X_{v^\epsilon}(t)dt
\\=&\frac{2}{r_\epsilon^2}\sum_{i=0}^\ell\hat\beta_i^\epsilon\{\sum_{\eta=0}^\ell\int_{F_\eta^\epsilon}\vec p(t)[A(t,u_\eta(t))-A(t,\bar u(t))]
Z^{-1}(t)
\\&
\int_0^tZ(s)[\vec f(s,u_i(s))-\vec f(s,\bar u(s))]ds
\}dt+o(1)
\\=&\frac{2}{r_\epsilon^2}\sum_{i=0}^\ell\hat\beta_i^\epsilon\sum_{\eta=0}^\ell\hat\beta_\eta^\epsilon\int_0^T\vec p(t)[A(t,u_\eta(t))-A(t,\bar u(t))]Z^{-1}(t)
\\&\int_0^tZ(s)[\vec f(s,u_i(s))-\vec f(s,\bar u(s))]ds\,dt+o(1)
\\=&\frac{2}{r_\epsilon^2}\sum_{\eta=0}^\ell \sum_{0\leq i\leq \eta}^\ell\hat\beta_i^\epsilon\hat\beta_\eta^\epsilon\int_{\tau_\eta-r_\epsilon}^{\tau_\eta}\vec p(t)[A(t,u(t))-A(t,\bar u(t))]Z^{-1}(t)
\\&\cdot\int_{[0,t]\cap [\tau_i-r_\epsilon,\tau_i]}\mathcal A(s,u(s))dsdt
+o(1)
\\=&\frac{2}{r_\epsilon^2}\sum_{\eta=0}^\ell \hat\beta_\eta^\epsilon\int_{\tau_\eta-r_\epsilon}^{\tau_\eta}\vec p(t)[A(t,u(t))-A(t,\bar u(t))]Z^{-1}(t)
\\&\left\{\sum_{0\leq i<\eta}\hat\beta_i^\epsilon\int_{\tau_i-r_\epsilon}^{\tau_i}\mathcal A(s,u(s))ds+\hat\beta_\eta^\epsilon\mathcal A(\tau_\eta,u(\tau_\eta))(t-\tau_\eta+r_\epsilon)\right.
\\\label{n721}&\left.+\hat\beta_\eta^\epsilon\int_{\tau_\eta-r_\epsilon}^t [\mathcal A(s,u(s))-\mathcal A(\tau_\eta,u(\tau_\eta))]ds\right\}dt+o(1)
\\=&\frac{2}{(\sum_{p=0}^\ell\beta_p)^2}\sum_{\eta=0}^\ell\beta_\eta\vec p(\tau_\eta)[A(\tau_\eta,u(\tau_\eta))-A(\tau_\eta,\bar u(\tau_\eta))]Z^{-1}(\tau_\eta)
\\ &\cdot\sum_{0\leq i<\eta}\beta_i\mathcal A(\tau_i,u(\tau_i))+\frac{2}{r_\epsilon^2}\sum_{\eta=0}^\ell(\hat\beta_\eta^\epsilon)^2 \vec p(\tau_\eta)[A(\tau_\eta,u(\tau_\eta))-A(\tau_\eta,\bar u(\tau_\eta))]
\\&\cdot Z^{-1}(\tau_\eta)\mathcal A(\tau_\eta,u(\tau_\eta))\int_{\tau_\eta-r_\epsilon}^{
\tau_\eta}(t-\tau_\eta+r_\epsilon)dt
\\&+\frac{2}{r_\epsilon^2}\sum_{\eta=0}^\ell (\hat\beta_\eta^\epsilon)^2\int_{\tau_\eta-r_\epsilon}^{\tau_\eta}\left\{\vec p(t)[A(t,u(t))-A(t,\bar u(t))]Z^{-1}(t)\right.
\\&\left.-\vec p(\tau_\eta)[A(\tau_\eta,u(\tau_\eta))-A(\tau_\eta,\bar u(\tau_\eta))]Z^{-1}(\tau_\eta)\right\}\mathcal A(\tau
_\eta,u(\tau_\eta))(t-\tau_\eta+r_\epsilon)dt
\\&+\frac{2}{r_\epsilon^2}\sum_{\eta=0}^\ell(\hat\beta_\eta^\epsilon)^2\int_{\tau_\eta-r_\epsilon}^{\tau_\eta}\vec p(t)[A(t,u(t))-A(t,\bar u(t))]Z^{-1}(t)
\\&\int_{\tau_\eta-r_\epsilon}^t[ \mathcal A(s,u(s))-\mathcal A(\tau_\eta,u(\tau_\eta))]dsdt+o(1).
\end{array}\end{align}
 Set
 \begin{align*}
A_i^\epsilon=&\{\tau;\,|\tau-\tau_i|\leq r_\epsilon, |\mathcal A(\tau,u(\tau))-\mathcal A(\tau_i,u(\tau_i))|>\epsilon\},
\\  B_i^\epsilon=&\left\{\tau;  |\tau-\tau_i|\leq r_\epsilon, \left|Z(\tau)[A(\tau,u(\tau))-A(\tau,\bar u(\tau))]Z^{-1}(\tau)-\right.\right.
\\&\left.\left.Z(\tau_i)[A(\tau_i,u(\tau_i))-A(\tau_i,\bar u(\tau_i))]Z^{-1}(\tau_i)\right|>\epsilon\right\},
\end{align*}
   for $i=0,\cdots,\ell$. Recalling (\ref{n725}),  (\ref{n720}) and (\ref{n903}), we have
\begin{align*}
&\frac{1}{r_\epsilon^2}\int_{\tau_\eta-r_\epsilon}^{\tau_\eta}\vec p(t)[A(t,u(t))-A(t,\bar u(t))]Z^{-1}(t)\int_{\tau_\eta-r_\epsilon}^t[ \mathcal A(s,u(s))-\mathcal A(\tau_\eta,u(\tau_\eta))]dsdt
\\=&\frac{1}{r_\epsilon^2}\int_{\tau_\eta-r_\epsilon}^{\tau_\eta}\vec p(t)[A(t,u(t))-A(t,\bar u(t))]Z^{-1}(t)
\\&\left(\int_{[\tau_\eta-r_\epsilon,t]\cap (A^\epsilon_\eta)^c}[\mathcal A(s,u(s))-\mathcal A(\tau_\eta,u(\tau_\eta))]ds\right.
\\&\left.+\int_{[\tau_\eta-r_\epsilon,t]\cap A^\epsilon_\eta}[\mathcal A(s,u(s))-\mathcal A(\tau_\eta,u(\tau_\eta))]ds\right)dt
\\=&o(1),
\end{align*}
and
\begin{align*}
&\frac{1}{r_\epsilon^2}\int_{\tau_\eta-r_\epsilon}^{\tau_\eta}\left\{\vec p(t)[A(t,u(t))-A(t,\bar u(t))]Z^{-1}(t)\right.
\\&\left.-\vec p(\tau_\eta)[A(\tau_\eta,u(\tau_\eta))-A(\tau_\eta,\bar u(\tau_\eta))]Z^{-1}(\tau_\eta)\right\}\mathcal A(\tau_\eta,u(\tau_\eta))(t-\tau_\eta+r_\epsilon)dt
\\=&\frac{1}{r_\epsilon^2}\int_{[\tau_\eta-r_\epsilon,\tau_\eta]\cap (B_\eta^\epsilon)^c}\left\{\vec p(t)[A(t,u(t))-A(t,\bar u(t))]Z^{-1}(t)\right.
\\&\left.-\vec p(\tau_\eta)[A(\tau_\eta,u(\tau_\eta))-A(\tau_\eta,\bar u(\tau_\eta))]Z^{-1}(\tau_\eta)\right\}\mathcal A(\tau_\eta,u(\tau_\eta))(t-\tau_\eta+r_\epsilon)dt
\\&+\frac{1}{r_\epsilon^2}\int_{[\tau_\eta-r_\epsilon,\tau_\eta]\cap B_\eta^\epsilon}\left\{\vec p(t)[A(t,u(t))-A(t,\bar u(t))]Z^{-1}(t)\right.
\\&\left.-\vec p(\tau_\eta)[A(\tau_\eta,u(\tau_\eta))-A(\tau_\eta,\bar u(\tau_\eta))]Z^{-1}(\tau_\eta)\right\}\mathcal A(\tau_\eta,u(\tau_\eta))(t-\tau_\eta+r_\epsilon)dt
\\=&o(1).
\end{align*}
Consequently (\ref{n721}) is reduced to
\begin{align*}
&\frac{2}{r_\epsilon^2}\int_0^T \vec p(t)[ A(t,v^\epsilon(t))-A(t,\bar u(t))]\vec X_{v^\epsilon}(t)dt
\\=& \frac{2}{(\sum_{p=0}^\ell\beta_p)^2}\sum_{\eta=0}^\ell\beta_\eta \vec p(\tau_\eta)[A(\tau_\eta,u(\tau_\eta))-A(\tau_\eta,\bar u(\tau_\eta))]Z^{-1}(\tau_\eta)\sum_{0\leq i<\eta}\beta_i\mathcal A(\tau_i,u(\tau_i))
\\&+\frac{1}{(\sum_{p=0}^\ell\beta_p)^2}\sum_{\eta=0}^\ell(\beta_\eta)^2 \vec p(\tau_\eta)[A(\tau_\eta,u(\tau_\eta))-A(\tau_\eta,\bar u(\tau_\eta))]Z^{-1}(\tau_\eta)\mathcal A(\tau_\eta,u(\tau_\eta))+o(1).
\end{align*}
We devide (\ref{n761}) by $\frac{1}{r_\epsilon^2}$, and  insert (\ref{n730}), (\ref{n733}) and the above relation into it. As $\epsilon$ approaches to $0^+$,  we obtain (\ref{n740}) by using (\ref{n735}) and (\ref{n738}). $\Box$

\medskip
Then, we shall prove Theorem \ref{n790}.
\\

{\it Proof of Theorem \ref{n790}.}\quad By \cite[Theorem I.7.6, p.150]{war72}, the set valued map $U(\cdot)$ is measurable. We obtain from Castaing' theorem \cite[Theorem I.7.8, p.152]{war72} that, there exist measurable selections $\hat u_1(\cdot), \hat u_2(\cdot), \cdots,$ of $U(\cdot)$ such that $\{\hat u_1(t), \hat u_2(t), \cdots\}$ is dense in $U(t)$ for all $t\in[0,T]$. Let $\mathcal T\subset[0,T]$ be the set such that $\{\mathcal A(\cdot,
\hat u_\eta(\cdot))\}_{\eta\geq 1}$ and $\left\{  Z(\cdot)\Big(A(\cdot,\hat u_\eta(\cdot))-A(\cdot,\bar u(\cdot))\Big)\right.
\left.Z^{-1}(\cdot)\right\}_{\eta\geq 1}$ are approximately continuous over it. It is obvious that $|\mathcal  T|=T$.

Fix any  $\tau_0, \tau_1,\cdots, \tau_\ell\subset \mathcal  T$ with $0
<\tau_0<\cdots<\tau_\ell<T$ and $\ell\geq k+j$,  and any $r_i\in U(\tau_i)$ $(i=0,\cdots,\ell)$ and $\beta_0,\cdots,\beta_\ell\in(0,+\infty)$ satisfying (\ref{n791}) and (\ref{n792}),
Then, fix any small $\epsilon>0$.  There exist $u_0^\epsilon(\cdot),\cdots,u_\ell^\epsilon(\cdot)\in\{\hat u_\eta(\cdot)\}_{\eta\geq 1}$ such that
\begin{align*}
&|\mathcal A(\tau_i,r_i)-\mathcal A(\tau_i,u_i^\epsilon(\tau_i))|<\epsilon,
\\& |Z(\tau_i)\Big(A(\tau_i,r_i)-A(\tau_i,u_i^\epsilon(\tau_i))\Big)Z^{-1}(\tau_i)|<\epsilon,
\end{align*}
for $i=0,1,\cdots,\ell$, and
\begin{equation}\label{n803}0^{j+k}\in Int \,co\{(\nabla\Phi_1^\top,\cdots,\nabla\Phi_j^\top,\nabla\Psi^\top)^\top\mathcal A(\tau_i,u_i^\epsilon(\tau_i))\}_{i=0}^\ell.\end{equation}
Applying \cite[Lemma]{War78} to (\ref{n791}), we obtain that there exist $\beta_0^\epsilon,\cdots,\beta_\ell^\epsilon\in(0,+\infty)$ such that
\begin{align}\label{n796}
&\lim_{\epsilon\to0^+}\beta_\eta^\epsilon=\beta_\eta,\; \eta=0,1,\cdots,\ell;
\\ \label{n797}&\nabla\Phi_i\sum_{\eta=0}^\ell\beta_\eta^\epsilon\mathcal A(\tau_\eta,u_\eta^\epsilon(\tau_\eta))=0,\;i=0,1,\cdots,\ell,
\\ \label{n798}& \nabla\Psi \sum_{\eta=0}^\ell\beta_\eta^\epsilon\mathcal A(\tau_\eta,u_\eta^\epsilon(\tau_\eta))=0.
\end{align}
Set
$$u^\epsilon(t)= \sum_{i=0}^\ell I_{[\tau_i-\delta,\tau_i+\delta)}(t)u_i^\epsilon(t)+I_{(\cup_{i=0}^\ell[\tau_i-\delta,\tau_i+\delta])^c}(t)\bar u(t),\;\forall \,t\in[0,T],
$$
 where $\delta\in(0,\frac{1}{2}\min\{\tau_0,\tau_1-\tau_0,\cdots,\tau_\ell-\tau_{\ell-1},T-\tau_\ell\})$.
Then,
\begin{align*}
&u^\epsilon(t)\in U(t),\;\forall\, t\,\in[0,T],
\\
&u^\epsilon(\tau_i)=u_i^\epsilon(\tau_i),\;i=0,1,\cdots,\ell.
\end{align*}
and
$\mathcal A(\cdot,u^\epsilon(\cdot))$ and $Z(\cdot)\Big(A(\cdot,u^\epsilon(\cdot))-A(\cdot,\bar u(\cdot))\Big)Z^{-1}(\cdot)$ are approximately continuous at $\tau_0,\cdots,\tau_\ell$. Consequently, (\ref{n800}) holds with $u(\cdot)$ and $(\beta_0,\cdots,\beta_\ell)$ replaced respectively  by $u^\epsilon(\cdot)$ and $(\beta_0^\epsilon,\cdots,\beta_\ell^\epsilon)$. Recall (\ref{n803}). By Theorem \ref{n770}, we obtain
(\ref{n740}) with $u(\cdot)$ and  $(\beta_0,\cdots,\beta_\ell)$ replaced respectively  by $u^\epsilon(\cdot)$ and $(\beta_0^\epsilon,\cdots,\beta_\ell^\epsilon)$, and we obtain (\ref{n810})
when $\epsilon$ approaches to $0^+$. $\Box$

\setcounter{equation}{0}
\section{Appendix }\label{ap}
\def\theequation{5.\arabic{equation}}
%\hskip\parindent
\subsection{Exponential map}\label{a2}
For this part, we refer the readers to \cite[Chapter 3]{c} and \cite[Chapter 3]{wsy}.

A differentiable curve $\gamma(t)$ on $ M$ with $t\in[0,\alpha)$ (for some $\alpha>0$) is called a geodesic if it satisfies
$$
\nabla_{\dot{\gamma}(t)}\dot{\gamma}(t)=0,\quad t\in[0,\alpha).
$$
Fix $x\in  M$. For any $v\in T_x M$, there exists a unique geodesic $\gamma_v(\cdot)$ satisfying $\gamma_v(0)=x$ and $\dot{\gamma}_v(0)=v$. Let $[0, \ell_v)$ be  the maximal  interval on which $\gamma_v(\cdot)$ is defined. Let $O_x\subset T_x  M$ be the set of vectors $v$ such that $\ell_v>1$. Then one can define  the  exponential  map as follows
$$\exp_x: O_x\to  M,\quad  \exp_x v=\gamma_v(1).$$ It has been shown that $O_x$ is a neighborhood of the origin $O\in T_x M$, and $\exp_x$ maps straight line
segments  in $T_x  M$ passing through the origin $O\in T_x M$  to geodesic
segments  in $ M$ passing through $x$. For any $v\in T_x M$, the differential of $\exp_x$ at $v$ is a linear map, denoted by
$$
d\exp_{x}|_v: \;\ T_vT_x M
\to T_{\exp_xv} M,$$
 where $T_vT_x M $ denotes the tangent space of the manifold $T_x M$ at the point $v\in T_x M$.

Given an $\epsilon>0$, write
 \begin{equation}\label{409}
B(O,\epsilon)\equiv\{v\in T_x M;\ \  |v|< \epsilon\}
\ \
\textrm{and} \ \
B_x(\epsilon)\equiv\{y\in  M;\ \  \rho(x,y)<\epsilon\}.
\end{equation}
We call
$
i(x)\equiv\sup\{\epsilon>0;\ \  \hbox{The map }\exp_x:B(O,\epsilon)\to B_x(\epsilon) \textrm{ is diffeomorphic}\}
$
the injectivity radius  at the point $x$ (e.g., \cite[p. 142]{p1}).

We list the following property of the exponential map, which can be found in many books on Riemannian geometry (e.g. the proof of \cite[Proposition 2.9, p. 65]{c}).
\begin{lem}\label{317} For any $x\in  M$, the map $\exp_x$ is a local
diffeomorphism, whose  differential at the origin $O\in T_x M$ satisfies
\begin{equation}\begin{array}{c}\label{203}
d\exp_x|_O=d\exp_x^{-1}\Big|_x=\textrm{the identity operator on} \,\,T_x M.
\end{array}\end{equation}
Furthermore, for any $y\in  M$ with $\rho(x,y)<i(x)$, there exists a unique shortest piecewise smooth curve which is also a geodesic in $M$, connecting $x$ and $y$.
\end{lem}

\subsection{Parallel translation and tensors}\label{a1}
For the details of this part, we refer the readers to \cite[Chapter I and Chapter III ]{kn}, \cite[Chapter 2]{p1}, \cite[Chapter 1]{wsy} and \cite[Chapter 1]{h1}.

For any $x\in M$ and $r,s\in\mathbb N$, a multilinear map
 $$
F: \ \ \underbrace{T_x^*  M\times\cdots\times T_x^*  M}_{r\;\mbox{times}}\times \underbrace{T_x M\times\cdots\times T_x M}_{s\;\mbox{times}}\to \mathbb R
$$
is called a tensor of order $(r,s)$ at $x$. Denote by $\T_s^r(x)$ the tensor space of type $(r,s)$  at $x$.
A smooth tensor field $\T$
of type $(r,s)$ on  $M$ is  a smooth assignment of a tensor $ \T(x)\in T_s^r(x)$ to each point $x$ of $ M$.
The norm of $\T$ at $x\in M$ is defined as follows:
 \begin{equation}\label{270}
\begin{array}{r}|\T(x)|=\sup\big\{\T(x)(Y_1,\cdots,Y_r,\lambda_1,\cdots,\lambda_s);\ \ Y_j\in T^*_x M,\lambda_l\in T_x M,\\[3mm]
 |Y_j|\leq 1,\, |\lambda_l|\leq 1, j=1,\cdots,r, l=1,\cdots,s\big\},\quad x\in M.\end{array}
\end{equation}
Denote by $\T_s^r(M)$ the set of all tensor fields of type $(r,s)$ over $M$.

Let $\gamma: [0,\ell]\to  M$ ($l>0$) be a differentiable  curve  with $\gamma(0)=x\in M$ and $\gamma(\ell)=y\in M$.
Given a vector $v\in T_xM$, there exists a unique vector field $X$ along $\gamma$ satisfying
\begin{equation}\label{n7}
\nabla_{\dot{\gamma}(s)}X=0,\qquad\forall\; s\in[0,\ell], \quad X(\gamma(0))=v.
\end{equation}
The mapping $T_x M\ni v\mapsto X(\gamma(\ell))\in T_{y} M$ is a linear isometry between $T_x M$ and $T_{y} M
$. We call this map the parallel translation along the curve $\gamma$, and denote it by $L^{\gamma}_{xy}v$. The
parallel translation along the curve $\gamma$ enjoys the following property:
\begin{equation}\label{r1}
\langle L_{xy}^\gamma v, L_{xy}^\gamma w\rangle=\langle v,w\rangle,\quad\forall v,w\in T_xM.
\end{equation}

For any $\eta\in T_x^*M$, we define $L_{xy}^\gamma\eta\in T_y^*M$ by
\begin{equation}\label{r50}L_{xy}^\gamma\eta(X)=\eta((L_{xy}^\gamma)^{-1}X),\;\forall  X\in T_y M.
\end{equation}
One can extend the parallel translation of a vector  at $x\in M$ along the  curve $\gamma$ to a tensor $\T\in\T_s^r(x)$ by
$$
L_{xy}^\gamma \T(v_1,\cdots,v_r,\eta_1,\cdots,\eta_s)=\T((L_{xy}^\gamma)^{-1}v_1,\cdots,(L_{xy}^\gamma)^{-1}v_r,(L_{xy}^\gamma)^{-1}\eta_1,
\cdots,(L_{xy}^\gamma)^{-1}\eta_s),
$$
 for all $v_1,\cdots,v_r\in T_y^*M$ and $\eta_1,\cdots,\eta_s\in T_y M$.

In particular, if $\rho(x,y)<\min\{i(x),i(y)\}$, according to Lemma \ref{317}, there is a unique shortest geodesic  $\gamma$ connecting  $x$ and $y$. In this case, we employ $L_{xy}$ instead of $L_{xy}^{\gamma}$ for abbreviation.

Let $\T$ be a tensor field. Take any $v\in T_xM$. Let $\gamma$ be a smooth curve  such that $\gamma(0)=x$ and $\dot{\gamma}(0)=v$. Then the covariant derivative of a tensor field   (in terms  of parallel translation) is defined as follows (see \cite[p. 42]{h1}):
\begin{align}\label{pt}
\nabla_v\T=\lim_{t\to 0}\frac{1}{t}\Big((L_{x\gamma(t)}^\gamma)^{-1} \T(\gamma(t))-\T(x)\Big).
\end{align}

Denote by $\nabla \T$ the covariant differential of $\T$, which is a tensor field of order $(r,s+1)$, and    is defined by (see \cite[p. 124]{kn})
\begin{equation}\label{100}
\nabla \T(x)(Y_1,\cdots,Y_r,\lambda_1,\cdots,\lambda_s,Z)=\nabla_Z\T(Y_1,\cdots,Y_r,\lambda_1,\cdots,\lambda_s),
\end{equation} for all $Y_1,\cdots,Y_r\in T^*_x  M$ and $\lambda_1,\cdots,\lambda_s,Z\in T_x M$.

In particular, a smooth function $f\in C^\infty( M)$ is a tensor of order $(0,0)$. $\nabla f$ and $\nabla^2f$ are respectively  tensors of order $(0,1)$ and $(0,2)$. We obtain from (\ref{pt}) and the definition of differential of a smooth function that
\begin{align}\label{ced}
\nabla f=df.
\end{align}
We call $\nabla^2f$  the Hessian of the function $f$, which  is a symmetric tensor, and  can be computed by
\begin{equation}\label{223}
\nabla^2f(x)(X,Y)=Y(x) (X f)-(\nabla_{Y(x)}X)f,\qquad x\in  M,\ X,Y\in T M.
\end{equation} For a smooth function $h: M\times M\to \mathbb R$ of two arguments, we denote by $\nabla_i h$ the covariant derivative of $h$ with respect to  the $i^{th}$ argument with $i=1,2$.
The  higher order derivatives of $h$  are defined as follows: For $i,j=1,2$ with $i\neq j$, any $(x_1,x_2)\in M\times M$ and $X,Y,Z\in \mathcal X(M)$,
\begin{equation}\label{108}\begin{array}{ll}
\displaystyle\nabla_i\nabla_jh(x_1,x_2)(X,Y)\equiv Y(x_i)\Big(X(x_j)(h(x_1,x_2))\Big)=Y(x_i)(\langle \nabla_j h(x_1,x_2),X(x_j)\rangle);
\\[3mm] \displaystyle\nabla_i^2h(x_1,x_2)(X,Y)\equiv Y(x_i)\Big(X(x_i)h(x_1,x_2)\Big)-\nabla_{Y(x_i)}Xh(x_1,x_2);
\\[3mm] \displaystyle\nabla_i^2\nabla_j h(x_1,x_2)(X,Y,Z)\equiv\nabla_i^2(\langle X(x_j),\nabla_j h(x_1,x_2)\rangle)(Y,Z);
\\[3mm] \displaystyle\nabla_i\nabla_j^2h(x_1,x_2)(X,Y,Z)\equiv Z(x_i)\Big(\nabla_j^2h(x_1,x_2)(X,Y)\Big).
\end{array}\end{equation}

\subsection{Useful lemmas}\label{a3}

\begin{lem}\label{17}(\cite[Lemma 2.2]{cdz})\,\,
For any $x,y\in  M$ with $\rho(x,y)<\min\{i(x),i(y)\}$, $X,X_1,X_2\in T_x M$ and $Y\in T_y  M$, it holds that
\begin{eqnarray}
&&\label{80}|\displaystyle\exp_x^{-1}y|=|\exp_y^{-1}x|=\rho(x,y), \qquad \nabla_{X_1}L_{x\cdot}X=0,
\\[2mm] &&\displaystyle\label{81}\nabla_1\rho^2(x,y)=-2\widetilde
{\exp_x^{-1}y},\qquad\nabla_2\rho^2(x,y)=-2\widetilde{\exp_y^{-1}x},
\\[2mm] &&\label{18}\displaystyle
L_{xy}exp_x^{-1}y=-exp_y^{-1}x, \qquad L_{xy}d_1\rho^2(x,y)=-d_1\rho^2(y,x),
\\[2mm] &&\label{92}\displaystyle
\nabla_1\nabla_2\rho^2(x,y)(Y,X)=-2\langle d\exp_y^{-1}|_xX,Y\rangle,
\\[2mm] &&\label{56}\displaystyle\langle d\exp_x^{-1}\Big|_yY,X\rangle=\langle d\exp_y^{-1}\Big|_x X,Y\rangle,
\\[2mm] &&\label{96}\displaystyle\nabla_1\nabla_2\rho^2(x,y)(Y,X)=-\nabla_1^2\rho^2(x,y)(L_{yx}Y,X)-\nabla_1\rho^2(x,y)(\nabla_X L_{y\cdot}Y),
\\[2mm] && \label{472}\displaystyle \nabla_1^2\rho^2(x,x)(X_1,X_2)=\nabla_2^2\rho^2(x,x)(X_1,X_2)=2\langle X_1,X_2\rangle,
\\[2mm] &&\label{260} \displaystyle\nabla_i\nabla_j^2\rho^2(x,x)=\nabla_i^2\nabla_j\rho^2(x,x)= \nabla_i^3\rho^2(x,x)=0,\quad i,j=1,2,\,\,i\neq j,
\end{eqnarray} where the notions $\nabla_1\nabla_2\rho^2$, $\nabla_i^2\nabla_j\rho^2$ and  $\nabla_i\nabla_j^2\rho^2$ with $i,j=1,2$ and $i\neq j$ are defined in (\ref{108}),  $\nabla_i^2\rho^2$ is the Hessian of $\rho^2$ with respect to the $i^{th}$ argument, $\nabla_i^3\rho^2$ is the covariant derivative of the Hessian $\nabla_i^2\rho^2(x,x)$ with respect to  the $i^{th}$ argument (see (\ref{100})),    $d_i$ stands for the exterior derivative of a function on $M\times M$ with respect to the $i^{th}$ argument for $i=1,2$, and $\widetilde
{\exp_x^{-1}y}$ is the dual covector of $\exp_x^{-1}y$.
\end{lem}

\medskip

Denote by $[X,Y]\equiv XY-YX$ the Lie bracket of vector fields $X$ and $Y$.
Denote by  $R$  the curvature tensor (of $(M,g)$), which is a correspondence that associates to every pair $X,Y\in \mathcal X(M)$ a mapping $R(X,Y): \mathcal X(M)\to \mathcal X(M)$ given by $$R(X,Y)Z=\nabla_X\nabla_Y Z-\nabla_Y\nabla_X Z-\nabla_{[X,Y]}Z,\qquad \forall\; Z\in \mathcal X(M).$$
We write
$$
R(X,Y,Z,W)=\langle R(X,Y)Z,W\rangle,\qquad \forall\; X,Y,Z,W\in \mathcal X(M).
$$

\begin{lem}\label{425} (\cite[Lemma 4.1]{cdz}) Let $\T$ be a tensor field on $ M$. Then, the following two conditions are equivalent:
\begin{description}
\item[$(i)$] There exists a positive constant $L$ such that
$
|\nabla \T|\leq L$;
\item[$(ii)$] There exists a positive constant $L$ such that
$
|L_{x_1x_2} \T(x_1)- \T(x_2)|\leq L\rho(x_1,x_2)$,
 for all $x_1, x_2\in  M$ with $\rho(x_1,x_2)<\min\{i(x_1),i(x_2)\}$.
\end{description}
\end{lem}

\end{document}